\documentclass{amsart}
\usepackage{amssymb,amsmath,amsfonts,euscript,mathrsfs,amsthm}
\usepackage[dvips]{graphicx}

\usepackage[dvips]{color}

\input xy
\xyoption{all}
\CompileMatrices

\DeclareMathOperator{\Id}{Id}

\def\cA{{\mathscr A}}

\def\cC{{\mathscr C}}

\def\cF{{\mathscr F}}

\def\cH{{\mathscr H}}

\def\cK{{\mathscr K}}

\def\cL{{\mathscr L}}

\def\cT{{\mathscr T}}
\def\cR{{\mathscr R}}

\def\cU{{\mathscr U}}
\def\cV{{\mathscr V}}
\def\cW{{\mathscr W}}

\def\fh{{\mathfrak h}}
\def\fg{{\mathfrak g}}

\def\fl{{\mathfrak l}}

\def\fs{{\mathfrak s}}

\def\Hom{\operatorname{Hom}}

\def\End{\operatorname{End}}

\def\res{\operatorname{res}}

\def\id{\operatorname{id}}

\def\vac{|0\rangle}

\newtheorem{thm}{Theorem}[section]
\numberwithin{equation}{thm}

\newtheorem{prop}[thm]{Proposition}
\newtheorem{lem}[thm]{Lemma}
\newtheorem{cor}[thm]{Corollary}

\theoremstyle{definition}
\newtheorem{defn}[thm]{Definition}

\newtheorem{ex}[thm]{Example}

\theoremstyle{remark}

\newtheorem{rem}[thm]{Remark}

\newtheorem{nolabel}[thm]{}

\DeclareMathOperator{\lie}{Lie}

\newcommand{\ys}{Y}

\begin{document}
\title{Supersymmetric Vertex algebras}
\author{Reimundo Heluani}
\address{Department of Mathematics,  MIT \\ Cambridge, MA 02139 \\ USA}
\email{heluani@math.mit.edu}

\author{Victor G. Kac}
\address{Department of Mathematics, MIT \\ Cambridge, MA 02139 \\ USA}
\email{kac@math.mit.edu}
\date{}
\thanks{Supported in part by NSF grants DMS-0201017 and DMS-0501395}

\begin{abstract}We define and study the structure of \emph{ SUSY Lie
	conformal and vertex
algebras}. This leads to effective rules for computations with
superfields.  
\end{abstract}
\maketitle
\tableofcontents
\section{Introduction}
\begin{nolabel}
	Vertex algebras were introduced about 20 years ago by Borcherds
\cite{borcherds}. They provide a rigorous definition of the chiral part of
2-dimensional conformal field theory, intensively studied by physicists. Since
then,
they have had important applications to string theory and conformal field
theory, and to mathematics, by providing tools to study the most interesting
representations of infinite dimensional Lie algebras. Since their appearance, they
have been extensively studied in many papers and books (for the latter we
refer to \cite{flm}, \cite{frenkelhuang}, \cite{kac:vertex}, \cite{huang}, \cite{frenkelzvi},
\cite{beilinsondrinfeld}).

%
\label{no:intro.1}
\end{nolabel}
\begin{nolabel}
	The purpose of the present paper is to define and study the structure of
	supersymmetric (SUSY) vertex algebras. This theory encompasses the formalism
	of superfields extensively used by physicists (see eg. \cite{cohn},
	\cite{dolgikh} and references therein). 

	Recall \cite{kac:vertex} that a vertex algebra $(V, \vac, T, Y)$ is a
	vector superspace $V$ (space of states) with an even vector $\vac$ (vacuum
	vector), an even endomorphism $T$ (translation operator), and a parity
	preserving bilinear product with values in Laurent series in an
	indeterminate $z$ over $V$:
	\begin{equation*}
		V \otimes V \rightarrow V( (z)), \quad a\otimes b \mapsto Y(a,z) b
		= \sum_{n \in \mathbb{Z}} (a_{(n)}b) z^{-1-n},
	\end{equation*}
	subject to the following axioms $(a,\,b \in V)$:
	\begin{itemize}
		\item (vacuum axioms) $Y(a,z)\vac|_{z=0} = a$,
			$T \vac = 0$,
		\item (translation invariance) $[T,Y(a,z)] = \partial_z Y(a,z)$, 
		\item (locality) $(z-w)^N [Y(a,z), Y(b,w)]= 0$ for some $N \in
			\mathbb{Z}_+$. 
	\end{itemize}
	The vertex algebras which arise in conformal field theory carry a
	\emph{conformal vector} $\nu \in V$, so that the coefficients of the
	corresponding field 
	\begin{equation}
		Y(\nu,z) \equiv L(z) = \sum_{n \in \mathbb{Z}} L_n z^{-2-n}
		\label{eq:intro.1.2}
	\end{equation}
	satisfy the Virasoro commutation relations
	\begin{equation}
		[L_m, L_n] = (m-n) L_{m+n} + \delta_{m, -n} \frac{m^3- m}{12} c,
		\label{eq:intro.1.3}
	\end{equation}
	where $c \in \mathbb{C}$ is the \emph{central charge}, and also the
	following two properties hold:
	\begin{itemize}
		\item \begin{equation} L_{-1} = T, \label{eq:intro.pedo} 
			\end{equation} 
		\item $L_0$ is diagonalizable on $V$ with eigenvalues bounded
			below.
	\end{itemize}
	Sometimes, such a vertex algebra $V$ (called \emph{conformal}) admits a
	``supersymmetry'',  namely, V carries a
	\emph{superconformal vector $\tau$}, such that the corresponding field
	\begin{equation}
		Y(\tau,z) \equiv G(z) = \sum_{n \in \frac{1}{2} + \mathbb{Z}} G_n
		z^{- 3/2 - n}
		\label{eq:intro.1.4}
	\end{equation}
	satisfies the following properties:
	\begin{enumerate}
		\item $\nu = \tfrac{1}{2} G_{-1/2} \tau$ is a conformal vector with
			central charge $c$;
		\item the operators $G_n$ $(n \in 1/2 + \mathbb{Z})$ along with the
			Virasoro operators $L_n$ $(n \in \mathbb{Z})$, appearing in
			$L(z) = Y(\nu, z)$, form the Neveu-Schwarz algebra with
			central charge $c$, namely, along with the Virasoro
			relations (\ref{eq:intro.1.3}), the following commutation
			relations hold:
			\begin{equation}
				[G_m,L_n] = \left( m - \frac{n}{2} \right) G_{m+n},
				\qquad [G_m, G_n] = 2 L_{m+n} +
				\frac{c}{3} \left( m^2 - \frac{1}{4} \right)
				\delta_{m, -n}.
				\label{eq:intro.1.5}
			\end{equation}
	\end{enumerate}
	Then $V$ is called called an $N=1$ superconformal vertex algebra.
	In particular, we then can obtain an enhanced translation invariance property as
	follows. Let $S = G_{-1/2}$, let $\theta$ be an odd indeterminate,
	commuting with $z$: $\theta^2 = 0$, $\theta z = z \theta$, and to each $a
	\in V$ associate a \emph{superfield}:
	\begin{equation}
		Y(a,z,\theta) = Y(a,z) + \theta Y(Sa,z).
		\label{eq:intro.1.6}
	\end{equation}
	Then one can show, using the so called commutator formula, that 
	\begin{equation}
		[S, Y(a,z,\theta)] = (\partial_\theta - \theta \partial_z)
		(a,z,\theta).
		\label{eq:intro.1.7}
	\end{equation}
	Since, by the second formula in (\ref{eq:intro.1.5}) and
	(\ref{eq:intro.pedo}), we have 
	\begin{equation}
		S^2 = T,
		\label{eq:intro.1.8}
	\end{equation}
	formula (\ref{eq:intro.1.7}) implies the usual translation invariance
	$	[T, Y(a,z,\theta)] = \partial_z Y(a,z,\theta)$.
	
	This leads to the following definition of an $N_K = 1$ \emph{SUSY vertex
	algebra} $(V$, $\vac$, $S$, $Y)$, where $S$ is an odd endomorphism of the space
	of states $V$, and $Y$ is a parity preserving bilinear product with values
	in $V( (z))[\theta]$:
	\begin{equation*}
		V \otimes V \rightarrow V( (z))[\theta], \quad a\otimes b \mapsto
		Y(a,z,\theta)b = \sum_{\stackrel{n \in \mathbb{Z}}{i = 0,1}}
		\theta^{1-i} z^{-1-n} a_{(n|i)} b,
	\end{equation*}
	subject to the following axioms:
	\begin{itemize}
		\item (vacuum axioms) 
			$Y(a,z,\theta)\vac|_{z = 0, \theta = 0} = a$, $S \vac =
			0$,
		\item (translation invariance) $[S, Y(a,z,\theta)] =
			(\partial_\theta - \theta \partial_z) Y(a,z,\theta)$,
		\item (locality) $(z-w)^N [Y(a,z,\theta), Y(b,w,\zeta)] = 0$ for
			some $N \in \mathbb{Z}_+$. 
	\end{itemize}
	This SUSY vertex algebra is called \emph{conformal} if there exists an odd
	vector $\tau \in V$, such that 
	\begin{equation*}
		Y(\tau,z,\theta) = G(z) + 2 \theta L(z),
	\end{equation*}
	and such that the coefficients of the expansions of these fields as in (\ref{eq:intro.1.2}) and
	(\ref{eq:intro.1.4}), 
	satisfy the commutation relations (\ref{eq:intro.1.3}) and
	(\ref{eq:intro.1.5}) of the Neveu-Schwarz algebra, $S = G_{-1/2}$ and $L_0$
	is diagonalizable with eigenvalues bounded below.

	Of course, a SUSY vertex algebra $(V, \vac, S, Y(a,z,\theta))$ gives rise
	to an ordinary vertex algebra $(V, \vac, T=S^2, Y(a,z) = Y(a,z,0))$, and a
	superconformal vector $\tau$ gives rise to a conformal vector $\nu =
	\tfrac{1}{2} S \tau$. 
	However, computations with superfields are much more effective than with
	the ordinary fields.
		\label{no:intro.1.1}
\end{nolabel}
\begin{nolabel}
	Now, the Neveu-Schwarz is the simplest, after Virasoro, among the
	superconformal Lie algebras. A
	\emph{superconformal Lie algebra} is a pair $(\cL, \cF)$, where $\cL$ is a
	Lie superalgebra and
	\begin{equation*}
		\cF = \left\{ a^i(z) = \sum_{n \in \mathbb{Z}} a_n^i z^{-1-n}
		\right\}_{i \in I}
	\end{equation*}
	is a finite family of formal distributions whose coefficients $a_n^i$
	span $\cL$, satisfying the following properties:
	\begin{enumerate}
		\item $\cF$ contains a Virasoro formal distribution $L(z) = \sum_{n
			\in \mathbb{Z}} L_n z^{-2-n}$, such that $[L_{-1}, a^i(z)] =
			\partial_z a^i(z)$;
		\item the formal distributions $a^i(z)$ are pairwise local:
			$(z-w)^N [a^i(z), a^j(w)]= 0$ for some $N \in
			\mathbb{Z}_+$, or, equivalently:
			\begin{equation*}
				[a^i(z), a^j(w)] = \sum_{s = 0}^{N-1} c^{ij}_s(w)
				\partial_w^s \delta(z,w) 
			\end{equation*}
			for some formal distributions $c^{ij}_s(z)$;
		\item $c^{ij}_s (z) \in \mathbb{C}[ \partial_z] \cF$;
		\item $\cL$ is a non-split central extension of an almost simple Lie
			superalgebra $\tilde{\cL}$ (i.e. all non-zero ideals of
			$\tilde{\cL}$ contain its derived subalgebra
			$\tilde{\cL}'$).
	\end{enumerate}
\label{no:intro.2}
\end{nolabel}
\begin{nolabel}
		A complete list of centerless superconformal Lie algebras consists of four
	series $(N \in \mathbb{Z}_+)$: $W(1|N)$, $S(1| N+2; a)$, $\tilde{S}(1|
	N+2)$, $K(1|N)$, and two exceptions $K(1|4)'$ and $CK(1|6)$ (see
	\cite{kacfattori}). Here $W(1|N)$ denotes the Lie superalgebra of all
	derivations of $\mathbb{C}[z,z^{-1}, \theta^1, \dots, \theta^{N}]$, where
	$z$ is an even indeterminate and $\theta^i$'s are odd anticommuting
	indeterminates, commuting with $z$ (in particular, $W(1|0)$ is the centerless
	Virasoro algebra).  The Lie superalgebras $S(1|N;a)$ and $\tilde{S}(1|N)$
	(resp. $K(1|N)$) are subalgebras of $W(1|N)$ consisting of vector fields,
	annihilating certain supervolume forms (resp. preserving the super-contact
	form\footnote{Here and further we use the sign convention of
	 \cite{deligne2},which differs
	 from the one used in \cite{kacleur} and\cite{kacfattori}. The main difference is that the de Rham
	differential $d$ is an \emph{even} derivation.} $dz + \sum_{i=1}^N \theta^i d\theta^i$, up to multiplication by a
	function), $K(1|4)'$ is the derived algebra of $K(1|4)$, and $CK(1|6)$ is a
	certain subalgebra of $K(1|6)$. The only algebras on the list that admit a
	central extension are $W(1|N)$ with $N \leq 1$, $S(1|2;a)$,
	$\tilde{S}(1|2)$, $K(1|N)$ with $N \leq 4$ and $K(1|4)'$ (see
	\cite{kacleur} and \cite{kacfattori}). Note that the Neveu-Schwarz algebra
	is a central extension of $K(1|1)$. 

	The subalgebra $\cL_\leq = \mathrm{span } \{ a^i_n | i \in I, \, n \geq 0\}$
	is called the annihilation subalgebra of $\cL$. The reason for this name
	comes from the fact that the space 
	\begin{equation}
		V(\cL) = U(\cL)/U(\cL) \cL_\leq,
		\label{eq:intro.2.3}
	\end{equation}
	where $U(\cL)$ denotes the universal enveloping superalgebra of $\cL$, carries a
	canonical structure of a vertex algebra, called the universal enveloping
	vertex algebra of $\cL$, where $\vac$ is the image of $1$, $T$ is induced
	by $L_{-1}$, and the image of $\cF$ is contained in the space of fields of
	$V(\cL)$, so that $\cL_\leq \vac = 0$. 

	In all cases the subalgebra $\cL_\leq$ consists of all vector fields of $\cL$
	which are regular at the origin. Denote by $\cL_<$ the subalgebra of
	$\cL_\leq$, consisting of vector fields, which vanish at the origin. Then in
	all cases we can find a complementary subalgebra $\cL_{\mathrm{tr}}$ (the
	``translation'' subalgebra)  
	\begin{equation}
		\cL_\leq = \cL_{\mathrm{tr}} \oplus \cL_<.
		\label{eq:intro.2.4}
	\end{equation}
	In the cases $W$ and $S$ one can choose $\cL_{\mathrm{tr}} = \mathrm{span}_{\mathbb{C}}
	\{\partial_z, \partial_{\theta^1}, \dots, \partial_{\theta^N}\}$. However,
	in the remaining cases the simplest choice is 
	\begin{equation*}
		\cL_{\mathrm{tr}} = \mathrm{span}_\mathbb{C} \{\partial_z, \partial_{\theta^1}
		- \theta^1 \partial_z, \dots, \partial_{\theta^N} - \theta^N
		\partial_z \}.
	\end{equation*}
	Denote these Lie superalgebras by $W_\mathrm{tr}$ and $K_\mathrm{tr}$ respectively, and by
	$\cH_W$  (resp. $\cH_K$) the universal enveloping superalgebra of
	$W_\mathrm{tr}$ (resp. $K_\mathrm{tr}$). It is an associative algebra on
	one even
	generator $T$ and $N$ odd generators $S^i$, subject to the relations:
	\begin{equation*}
		[T, S^i] = 0, \quad [S^i, S^j] = 0 (\text{resp. } = 2 \delta_{i,j}
		T).
	\end{equation*}
	Note that the operators $T$ and $S$ of an $N_K=1$ SUSY vertex algebra
	define a representation of the associative superalgebra $\cH_K$ for $N=1$
	(see \ref{eq:intro.1.8}). 

	This leads to the definition of an $N_W$ (resp. $N_K$) $=N$ SUSY vertex
	algebra as a $\cH_W$ (resp. $\cH_K$)-module $V$ with the vacuum vector
	$\vac$ and a parity preserving bilinear product with values in $V(
	(z))[\theta^i, \dots, \theta^N]$, 
	\begin{equation*}
		a \otimes b \mapsto Y(a,z,\theta^1, \dots, \theta^N) =
		\sum_{n \in \mathbb{Z}, J} \theta^{J^c} z^{-1-n} a_{(n|J)}b,
	\end{equation*}
	where $J$ runs over all ordered subsets of $\{1, \dots, N\}$, $J^c$ denotes the
	ordered complement of $J$ in $\{1, \dots, N\}$ and $\theta^J = \theta^{j_1}
	\dots \theta^{j_k}$ for $J = \{j_1, \dots, j_k\}$, subject to the vacuum,
	translation invariance and locality axioms. The vacuum and locality axioms
	generalize in the obvious way; translation invariance means the following:
	\begin{equation*}
		\begin{aligned}
			{[}T, Y(a,z,\theta)] &= \partial_z Y(a,z,\theta), \\
			[S^i, Y(a,z,\theta)] &= \partial_{\theta^i}
			Y(a,z,\theta) \quad \text{in the } N_W = N
			\text{ case},\\
			{[}S^i, Y(a,z,\theta)] &= (\partial_{\theta^i} - \theta^i
			\partial_z) Y(a,z,\theta) \quad \text{in the } N_K=N
			\text{ case}.
		\end{aligned}
	\end{equation*}
\end{nolabel}
\begin{nolabel}
	Given a superconformal Lie algebra $\cL$, and a decomposition
	(\ref{eq:intro.2.4}) of its annihilation subalgebra, such that $\cL_{\mathrm{tr}} \simeq
	W_\mathrm{tr}$ (resp. $ \simeq K_\mathrm{tr}$), we define an \emph{$\cL$-conformal $N_W$ (resp.
	$N_K$) $=N$ SUSY vertex algebra} $V$ by the property that $V$ carries a
	representation of $\cL$ such that the representation of $U(\cL_\mathrm{tr})$ coincides
	with that of $\cH_W$ (resp. $\cH_K$), the formal distributions of $\cL$ are
	represented by fields of $V$, and $L_0$ is diagonalizable with 
	eigenvalues bounded below. The ``minimal'' example of an
	$\cL$-superconformal vertex algebra is $V(\cL)$, given by
	(\ref{eq:intro.2.3}), or $V^c(\cL) = V(\cL)/(C - c 1) V(\cL)$ in the case
	$\cL$ has a central element $C$, called the \emph{SUSY Virasoro} vertex
	algebra associated to $\cL$.
	\label{no:intro.2aa}
\end{nolabel}
\begin{nolabel}
	We develop the structure theory of SUSY vertex algebras along the lines of
	\cite{kac:vertex}. 

	First, we develop the calculus of formal superdistributions. As usual, the
	key role is played by the formal super delta-function
	\begin{equation}
		\delta(Z,W) = (\theta^1 - \zeta^1)\dots(\theta^N -\zeta^N) (i_{z,w}
		- i_{w,z}) \frac{1}{z-w},
		\label{eq:intro.3.1}
	\end{equation}
	where $Z = (z,\theta^1, \dots, \theta^N)$, $W = (w, \zeta^1, \dots,
	\zeta^N)$, and $i_{z,w}$ signifies the expansion in the domain $|z|> |w|$.
	The key result is the decomposition formula of a local formal
	superdistribution:
	\begin{equation}
		a(Z,W) = \sum_{j \geq 0, J} \partial_W^{(j|J)} \delta(Z,W) c_{j|J}
		(W),
		\label{eq:intro.3.2}
	\end{equation}
	where 
	\begin{equation}
		c_{j|J}(W) = \res_Z (Z-W)^{j|J} a(Z,W),
		\label{eq:intro.3.3}
	\end{equation}
	and $\res_Z$ stands for the coefficient of $\theta^1 \dots \theta^N
	z^{-1}$. There are actually two cases, $W$ and $K$, of this formula,
	depending on the choice of $\partial_W^{j|J}$ in (\ref{eq:intro.3.2}) and
	the respective choice of $(Z-W)^{j|J}$ in (\ref{eq:intro.3.3}). 

	We let for $J = \{j_1 < \dots < j_s\}$:
	\begin{equation*}
		\begin{aligned}
			\partial_W^{j|J} &= \partial_w^j \partial_{\theta^{j_1}}
			\dots \partial_{\theta^{j_s}} && \text{ in case } W, \\
			\partial_W^{j|J} &= \partial_w^j
			(\partial_{\theta^{j_1}} + \theta^{j_1} \partial_w) \dots
			(\partial_{\theta^{j_s}} + \theta^{j_s} \partial_w)&& \text{
			in case } K,
		\end{aligned}
	\end{equation*}
	and in both cases we let $\partial_W^{(j|J)} = (-1)^{|J|(|J|+1)/2}
	\partial^{j|J}/j!$. (The second choice will be denoted by $D_W^{j|J}$.)
	Respectively, we let 
	\begin{equation*}
		\begin{aligned}
			(Z-W)^{j|J} &= (z-w)^j (\theta^{j_1} - \zeta^{j_1}) \dots
			(\theta^{j_s} - \zeta^{j_s}) && \text{  in case } W, \\
			(Z-W)^{j|J} &= \left( z-w-\sum_{i = 1}^N \theta^i \zeta^i
			\right)^j (\theta^{j_1} - \zeta^{j_1}) \dots
			(\theta^{j_s} - \zeta^{j_s}) &&\text{  in case } K.
		\end{aligned}
	\end{equation*}
\end{nolabel}
\begin{nolabel}
	Next, the formal Fourier transform is defined by 
	\begin{equation*}
		\cF^\Lambda_{Z,W} a(Z,W) = \res_Z \exp ( (Z-W)\Lambda)
		a(Z,W),
	\end{equation*}
	where $\Lambda = (\lambda, \chi^1, \dots, \chi^N)$, $\lambda$ is an even
	indeterminate, $\chi^i$'s are odd indeterminates, commuting with $\lambda$
	and satisfying the following commutation relations:
	\begin{equation*}
		\begin{aligned}
			\chi^i \chi^j + \chi^j \chi^i &= 0 &&\text{ in case } W,\\
			\chi^i \chi^j + \chi^j \chi^i &= - 2 \delta_{i,j} \lambda
			&&\text{ in case } K.
		\end{aligned}
	\end{equation*}
	Defining the \emph{$\Lambda$-bracket}
	\begin{equation*}
		[a(W)_\Lambda b(W)] = \cF^\Lambda_{Z,W} [a(Z),b(W)]
	\end{equation*}
	of formal superdistributions, we arrive at the notion of the $N_W$
	(resp. $N_K$)$=N$ SUSY Lie conformal algebra $\cR$, which, as usual,
	encodes the singular part of the operator product expansion (OPE) of a
	local pair of superfields. The case $N=0$ is that of an ordinary Lie
	conformal superalgebra \cite{kac:vertex}.

	The new phenomenon for $N > 0$ is that the $\Lambda$-bracket has parity $N
	\mod 2$. Consequently, the bracket $[a_\Lambda b]|_{\Lambda = 0}$ induces
	on $\cR / (T\cR + \sum_i S^i \cR)$ a structure of a Lie superalgebra of
	parity $N \mod 2$. 

	On the other hand, the structure of a SUSY Lie conformal algebra is an
	important part of the structure of a SUSY vertex algebra. As in the case of
	vertex algebras, the only missing ingredient is the normally ordered
	product of superfields, which is defined as usual:
	\begin{equation*}
		:a(Z)b(Z): = a_+(Z) b(Z) + (-1)^{p(a)p(b)} b(Z) a_-(Z),
	\end{equation*}
	where for a superfield $a(Z) = \sum_{n \in \mathbb{Z}, J} \theta^{J^c}
	z^{-1-n} a_{(n|J)}$ we let \[a_-(Z) = \sum_{n\geq 0} \theta^{J^c} z^{-1-n}
	a_{(n|J)}, \quad a_+(Z) = a(Z) - a_-(Z).\] We prove that the non-commutative
	Wick formula, which allows to compute the singular part of the OPE for
	normally ordered products, generalizes to the SUSY $N_W=N$ and $N_K=N$
	cases almost verbatim:
	\begin{equation}
		[a_\Lambda :bc:] = :[a_\Lambda b]c: + (-1)^{p(a) + N) b}
		:b[a_\Lambda c]:  + \int_0^\Lambda [ [a_\Lambda
		b]_\Gamma c] d\Gamma.
		\label{eq:intro.3.10}
	\end{equation}
	Furthermore, we show that the uniqueness and existence theorems, as well as
	all basic identities for vertex algebras, extend to the SUSY case with
	minor modification of signs. In particular, we show that a $N_W$ (resp. $N_K$)
	$=N$ SUSY vertex algebra is a $N_W$ (resp. $N_K$)$=N$ SUSY Lie conformal
	algebra with the $\Lambda$ bracket:
	\begin{equation*}
		[a_\Lambda b] = \res_Z e^{\Lambda Z} Y(a,Z)b,
	\end{equation*}
	together with a unital, ``quasicommuitative'' and ``quasiassociative''
	differential superalgebra structure $: :$, which are related by formula
	(\ref{eq:intro.3.10}). This is an equivalent definition of a SUSY vertex
	algebra, analogous to the one given in \cite{kacbakalov2} for vertex
	algebras. Removing ``quantum corrections'', we obtain the definition of a
	SUSY Poisson vertex algebra.
	\label{no:intro.3}
\end{nolabel}
\begin{nolabel}
	As in the vertex algebra case, most of the basic examples of SUSY vertex algebras
	are constructed starting with a SUSY Lie conformal algebra $\cR$, and
	taking its universal enveloping vertex algebra $V(\cR)$ or $V^c(\cR)$. This
	can be defined as the universal enveloping vertex algebra $V(\cL)$ of the
	corresponding formal superdistribution Lie algebra $V(\cL)$ (resp.
	$V^c(\cL)$), defined in the same way as in (\ref{eq:intro.2.3}). By abuse
	of terminology, we often say that $V(\cR)$ is a SUSY vertex algebra
	generated by $\cR$. 

	The simplest example of an $N_K = N$ SUSY vertex algebra is the well-known
	boson-fermion system, generated by one superfield $\Psi$ of parity $N \mod
	2$, subject to the following $\Lambda$-bracket:
	\begin{equation*}
			{[}\Psi_\Lambda \Psi] = \Lambda^{1|N} \text{ for even
			} N, \quad
		[\Psi_\Lambda \Psi] = \Lambda^{0|N} \text{ for odd } N.
	\end{equation*}
	Viewed as an ordinary vertex algebra, this SUSY vertex algebra in the case
	$N=1$ is generated by one boson and one fermion (hence the name
	``boson-fermion'' system). 

	Another example is the SUSY vertex algebra, generated by the current SUSY
	Lie conformal algebra, associated to a Lie superalgebra $\fg$ with an
	invariant bilinear form $(\, ,\,)$. For $N$ even, the corresponding SUSY
	$N_W = N$ or $N_K = N$ $\Lambda$-bracket is defined as
	\begin{equation*}
		[a_\Lambda b] = [a, b] + \lambda (a,b), \quad a,\,b \in \fg.
	\end{equation*}
	For $N$ odd, we consider the vector superspace $\Pi \fg$ with reversed
	parity, and define for $\bar{a}, \bar{b} \in \Pi\fg$:
	\begin{equation*}
		[\bar{a}_\Lambda\bar{b}] = (-1)^{p(a)} \left( \overline{[a,b]} +
		(a,b) \sum_{i=1}^N \chi^i \right).
	\end{equation*}
	In the case $N_K = 1$, using the normally ordered products of the above
	superfields, one reproduces the construction of the $N_K=1$ super Virasoro
	field from \cite{kactodorov}.

	Unfortunately, we do not know how to construct a SUSY lattice vertex
	algebra. 

	A less routine example is the following. Let $\fg$ be a Lie algebra and $F$
	a $\fg$-module. We associate to this data an $N_K=1$ SUSY Lie conformal
	algebra $\cR(\fg, F)$, which is a free $\cH_K$-module over the vector
	superspace:
	\begin{equation*}
		 (F \oplus \fg^*) \oplus (\fg \oplus \Pi \fg^*), 
	\end{equation*}
	where $(F \oplus \fg^*)$ and $(\fg \oplus \Pi \fg^*)$ are the even and odd
	parts  respectively \and $\Pi$ is the change of parity functor, with the following non-zero $\Lambda$-brackets
	($X, Y \in \fg$, $f \in F$, $\alpha \in \fg^*$, $\bar{\alpha} \in \Pi
	\fg^*$):
	\begin{equation*}
		\begin{aligned}
			{[}X_\Lambda Y] &= [X,Y], & [X_\Lambda f] &= X f, \\
			[X_\Lambda \alpha] &= X \alpha + \lambda <\alpha, X>, &
			[X_\Lambda \bar{\alpha}] &= \overline{X\alpha} + \chi
			<\alpha, X>.
		\end{aligned}
	\end{equation*}
	The corresponding SUSY universal enveloping vertex algebra is denoted by
	$V(\fg, F)$. The SUSY $N_K=1$ vertex algebra $V(\fg, F)$ in the case when
	$\fg$ is the Lie algebra of vector fields on a manifold $M$, $F$ is the
	space of functions on $M$ and $\fg^*$ is the space of differential
	$1-$forms on $M$, is used in \cite{heluani2} to construct the chiral de
	Rham complex \cite{malikov}, as a sheaf of $N_K=1$ SUSY vertex algebras on
	$M$, and study its SUSY properties.
	\label{no:intro.4}
\end{nolabel}
\begin{nolabel}
	In the subsequent paper \cite{heluani4} the formalism of SUSY vertex
	algebras is applied to associate to any (strongly)conformal SUSY $N_W$ (resp
	$N_K$)$=N$ vertex algebra $V$ a vector bundle $\cV_X$ on any supercurve
	(resp. superconformal curve) $X$, along the lines of \cite{frenkelzvi},
	where this is done for ordinary curves.
	\label{no:intro.5}
\end{nolabel}
%
\section{Preliminaries}\label{sec:intro}
\label{sub:super_recolection}
In this section we recall some notation and basic results on vertex algebras. We
also give the first examples of SUSY vertex algebras constructed via
ordinary vertex algebras. 
 The reader is referred to \cite{kac:vertex} for an introduction to the
vertex algebra theory. 
\begin{defn}
	Let $\cA$ be a 
	 Lie superalgebra. An \emph{$\cA$-valued formal distribution} is a formal
	expression of the form:
	\begin{equation*}
		B(z) = \sum_{n \in \mathbb{Z}} B_{(n)} z^{-1-n}
	\end{equation*}
	where $B_{(n)} \in \cA$ have the same parity for all $n \in
	\mathbb{Z}$; this parity is called the \emph{parity of $B(z)$}.
	The coefficients $B_{(n)}$ are called the
	 \emph{Fourier modes} of $B(z)$, and $z$ is an indeterminate. 
        A pair of formal distributions $B(z), \,
	C(w)$  is \emph{local} if 
	\begin{equation*}
		(z-w)^N [B(z), C(w)] = 0 \qquad \text{for some } N \in \mathbb{Z}_+.
	\end{equation*}
	If $\cA= \End(V)$, where $V$ is a vector superspace, with the usual
	superbracket,  we say that
	$B(z)$ is a \emph{field} if, for every $v \in V$, $B_{(n)}v = 0$ for
	large enough $n$.
\end{defn}
\begin{nolabel}
	Let $V$ be a vertex algebra\footnote{We will denote a vertex algebra by its
	space of states $V$ when there is no possible confusion.} (see
	\ref{no:intro.1.1}).
	The map $Y$ is called the \emph{state-field correspondence} and we
	will use this map to identify a vector $a \in V$ with its
	corresponding field $Y(a,z)$.

\end{nolabel}
\begin{nolabel}
	Given a vertex algebra $V$ we denote
	\begin{equation}
		\begin{aligned}
		  a_{(n)}b &= a_{(n)} \bigl(b\bigr), \qquad
			[a_\lambda b] = \sum_{k \geq 0}
			\frac{\lambda^k}{k!} a_{(k)}b, \qquad
			:ab: = a_{(-1)}b. 
		\end{aligned}
		\label{eqn:nth-product}
	\end{equation}
	The first operation is called the $n$-th product, the second is called
	the $\lambda$-bracket and the third the normally ordered product.
\end{nolabel}
\begin{nolabel}
	For each $n \in \mathbb{Z}$, define the $n$-th product of
	$\End(V)$-valued fields $A(z)$ and $B(z)$ as follows. Denote by
	$i_{z,w}$ the \emph{expansion in the domain $|z| > |w|$}: 
	\begin{equation}
		i_{z,w}  z^m w^n (z-w)^k = z^{m+k} w^n i_{z,w} \left(
		1 - \frac{w}{z} \right)^k = \sum_{j \geq 0} (-1)^j
		\binom{k}{j}
		z^{m+k-j} w^{n+j}.
		\label{eq:concrete_expansion}
	\end{equation}
	Define
	\begin{multline}
		A(w)_{(n)}B(w) = \res_z \left( i_{z,w} (z-w)^n A(z)B(w) - \right.
		\\ \left. -
		i_{w,z} (z-w)^n (-1)^{p(A)p(B)} B(w) A(z)\right),
		\label{}
	\end{multline}
	where $p(A)$ denotes the parity of $A(w)$. It can
	be shown that the following $n$-th product identity holds (cf. \cite[Prop. 4.4]{kac:vertex})
	\begin{equation}
		Y(a_{(n)}b,z) = Y(a,z)_{(n)}Y(b,z) \qquad \forall n \in
		\mathbb{Z},
		\label{eqn:nth-product-identity}
	\end{equation}
	hence,
	\begin{equation}
		Y(Ta,z) = \partial_z Y(a,z).
		\label{eqn:nth-product-identityb}
	\end{equation}
\end{nolabel}
\begin{nolabel}
	 In a vertex algebra $V$ we have the following commutator formulas \cite[p 112]{kac:vertex}
	 \begin{equation*}
		 \begin{aligned}
			 {[}a_{(m)}, b_{(n)}] &= \sum_{j \geq 0} 
			 \binom{m}{j} \left(
			 a_{(j)}b \right)_{(m+n-j)}, \\
			 {[}a_{(m)}, Y(b,w)] &= \sum_{j \geq 0}
			 \left(\frac{\partial_w^j w^m}{j!}\right)
			 Y(a_{(j)}b,w). 
		 \end{aligned}
	 \end{equation*}
	 This formula shows that the space of Fourier modes of all fields of a
	 vertex algebra is closed under the Lie bracket, and, moreover, the
	 commutation relations are expressed in terms of $j$-th products.
\end{nolabel}
\begin{defn}
	A \emph{Lie conformal algebra } is a super
	$\mathbb{C}[\partial]$-module $\cR$ equipped with a parity
	preserving bilinear map
	\begin{equation*}
		{[}\,_\lambda \, ] : \cR \otimes \cR \rightarrow
		\mathbb{C}[\lambda] \otimes \cR, \qquad a \otimes b \mapsto
		[a_\lambda b],
		\label{}
	\end{equation*}
	satisfying the following axioms:
	\begin{itemize}
		\item Sesquilinearity:
			\begin{equation*}
					{[}\partial a_\lambda b] = -
					\lambda [a_\lambda b],  \quad
					{[}a_\lambda \partial b] =
					(\partial + \lambda) [a_\lambda
					b].
				\label{}
			\end{equation*}
		\item Skew-commutativity:
			\begin{equation*}
				{[}b_\lambda a] = - (-1)^{p(a)p(b)}
				[a_{-\partial - \lambda} b].
				\label{}
			\end{equation*}
		\item Jacobi identity:
			\begin{equation*}
				[a_\lambda[b_\mu c]] = [ [ a_\lambda
				b]_{\lambda + \mu}c] + (-1)^{p(a)p(b)}
				[b_\mu[a_\lambda c] ] ,
				\label{}
			\end{equation*}
	\end{itemize}
	for all $a, \, b, \, c \, \in \cR$. Here $\lambda$ and $\mu$ are commuting
	indeterminates.

	Given a Lie conformal algebra $\cR$, we can associate to it a vertex
	algebra $V(\cR)$ (cf. \cite{kac:vertex}, \cite{kacbakalov2}) called the \emph{universal
	enveloping vertex algebra of $\cR$}, as defined in the introduction. If $\cR$ is generated by
	some vectors $\left\{ a_i \right\}$ as a
	$\mathbb{C}[\partial]$-module, we say that $V(\cR)$ is generated
	by the same vectors. If $C \in \cR$ is a central element such
	that $\partial C = 0$, given any complex number $c$, we denote by
	$V^c(\cR)$ the quotient of $V(\cR)$ by the ideal 
	$(C-c)V(\cR)$. 

	One can show \cite{kac:vertex} that a vertex algebra $V$ is
	canonically a Lie conformal algebra with
	the $\lambda$-bracket defined in (\ref{eqn:nth-product}) and
	$\partial = T$.
\end{defn}
\begin{ex}
	The \emph{Virasoro} vertex algebra $\mathrm{Vir^c}$ is generated by an even field
	$L$ satisfying:
	\begin{equation}
		[L_{\lambda}L] = (\partial + 2\lambda)L +
		\frac{c}{12} \lambda^3.
		\label{eqn:vir}
	\end{equation}
	The complex number $c$ is called the \emph{central charge}. Expanding the
	corresponding 
	field as in (\ref{eq:intro.1.2}) we obtain the familiar commutation
	relations of the Virasoro algebra (\ref{eq:intro.1.3}). 
\end{ex}
\begin{nolabel}
	Let $\nu \in V$ be a conformal vector (see \ref{no:intro.1.1}) and let $L(z)$
	be the corresponding
	Virasoro field. A vector $a \in V$ satisfying $[L_\lambda a] = (T + \Delta
	\lambda ) a + O(\lambda^2)$ is said to have \emph{conformal weight}
	$\Delta$. If, moreover, $a$ satisfies $[L_\lambda a] = (T + \Delta \lambda)
	a$ we say that $a$ is \emph{primary}. 
	\label{no:primary_def}
\end{nolabel}
\begin{ex} \label{ex:2.11}
	The Neveu Schwarz ($\mathrm{NS}$) vertex algebra is  generated by an even Virasoro
	field $L$ (satisfying (\ref{eqn:vir})) and an odd primary field $G$ of
	conformal weight $3/2$,
	satisfying the commutation
	relation:
	\begin{equation*}
		{[}G_\lambda G] = 2 L + \frac{\lambda^2}{3} c .
		\label{}
	\end{equation*}

	If we expand the corresponding fields as in (\ref{eq:intro.1.2}) and
	(\ref{eq:intro.1.4}) we obtain the commutation relations
	(\ref{eq:intro.1.5}) and (\ref{eq:intro.1.3}) of the Neveu-Schwarz algebra. 

	As we have seen in the introduction, given a vertex algebra with an $N=1$
	superconformal vector $\tau$, we obtain an $N_K=1$ SUSY vertex algebra (see
	also \ref{no:k.vertex.defin} for a definition) by defining the superfields
	(\ref{eq:intro.1.6}). In particular, the Neveu-Schwarz algebra gives rise
	to such an $N_K=1$ SUSY vertex algebra.
\end{ex}
Below we give some examples of vertex algebras with an $N=1$ superconformal vector. By the
construction in \ref{no:intro.1.1},
they are automatically $N_K=1$ SUSY vertex algebras.
\begin{ex} \cite[Ex. 5.9a]{kac:vertex} \label{ex:2.16}
	Let $V$ be the universal enveloping vertex algebra of the Lie conformal
	algebra generated by an even vector (free boson)
	 $\alpha$ and an odd vector (free fermion) $\varphi$, namely
	\begin{equation*}
		\begin{aligned}
			{[}\alpha_\lambda \alpha] = \lambda, \qquad
			{[}\varphi_\lambda \varphi] = 1, \qquad
			{[}\alpha_\lambda \varphi] = 0.
		\end{aligned}
		\label{}
	\end{equation*}
	Then $V$ is a (simple) vertex algebra with a family of $N=1$
	superconformal vectors 
	\begin{equation*}
		\tau = (\alpha_{(-1)} \varphi_{(-1)} + m
		\varphi_{(-2)})\vac, \qquad m \in \mathbb{C},
		\label{}
	\end{equation*}
	of central charge $c = \tfrac{3}{2} - 3 m^2$. 
\end{ex}
\begin{ex} \cite{kactodorov} \cite[Thm 5.9]{kac:vertex} \label{ex:g-super}
	Let $\fg$ be a finite dimensional Lie algebra with a
	non-degenerate invariant symmetric bilinear form $(\, ,
	)$, normalized by the condition $(\alpha,\alpha) = 2$ for a long root
	$\alpha$,  and let $h^\vee$ be the dual Coxeter number. We construct a vertex algebra $V^k(\fg_{\mathrm{super}})$
	generated by the usual currents $a,b \,  \in \fg$, satisfying:
	\begin{equation*}
		{[}a_\lambda b] = [a, b] + (k + h^\vee) \lambda (a, b),
		\label{}
	\end{equation*}
	and the odd super currents $\bar{a} \in  \Pi\fg$ (as before $\Pi$ stands
	for reversal of parity), satisfying:
	\begin{equation*}
		\begin{aligned}
			{[}a_{\lambda} \bar{b}] = \overline{[a,b]}, \qquad
			{[}\bar{a}_\lambda \bar{b}] = (k+h^\vee) (a, b).
	\end{aligned} 
		\label{}
	\end{equation*}
	Let $\left\{ a^i \right\}$ and $\left\{ b^i \right\}$  be dual
	bases of $\fg$. Provided that $k \neq -h^\vee$ the vertex algebra
	$V^k(\fg_{\mathrm{super}})$ admits an $N=1$ superconformal vector
	\begin{equation*}
		\tau = \frac{1}{k + h^\vee}\left( \sum_{i} a^i_{(-1)}
		\bar{b}^i_{(-1)} + \frac{1}{3(k + h^\vee)} \sum_{i,j,r} ([a^i, a^j], a^r) \bar{b}^i_{(-1)}
		\bar{b}^j_{(-1)} \bar{b}^r_{(-1)}\right) \vac,
		\label{}
	\end{equation*}
	of central charge
	\begin{equation*}
		c_k = \frac{k \mathrm{dim} \fg}{k + h^\vee} +
		\frac{1}{2} \mathrm{dim} \fg.
	\end{equation*}
	This is known as the \emph{Kac-Todorov} construction. The formulas in
	\cite{kac:vertex} should be corrected as above.
\end{ex}
\begin{ex} \cite[Thm 5.10]{kac:vertex} \label{ex:n=2}
	The $N=2$ vertex algebra is generated by a Virasoro field $L$ of central
	charge $c$, an
	even field $J$, primary of conformal weight $1$, and two odd fields
	$G^\pm$, primary of conformal weight $3/2$. The remaining commutation
	relations are:
	\begin{equation*}
	\begin{gathered}
			{[}J_\lambda J] = \frac{c}{3} \lambda, \quad 
			{[}G^\pm_\lambda G^\pm] = 0, \quad
			{[}J_\lambda G^\pm] = \pm G^\pm,\\
			{[}G^+_\lambda G^-] = L + \frac{1}{2} \partial J
			+ \lambda J + \frac{c}{6} \lambda^2.
		\label{}
	\end{gathered}
	\end{equation*}
	This vertex algebra contains an $N=1$ superconformal vector:
	\begin{equation*}
		\tau = G^+_{(-1)} \vac + G^-_{(-1)} \vac.
		\label{}
	\end{equation*}
	Also, this vertex algebra admits a $\mathbb{Z}/2 \mathbb{Z} \times
	\mathbb{C}^*$ family of automorphisms. The generator of
	$\mathbb{Z}/2 \mathbb{Z}$ is given by $L \mapsto L$, $J \mapsto
	-J$ and $G^\pm \mapsto G^\mp$. The $\mathbb{C}^*$ family is given
	by $G^+ \mapsto \mu G^+$ and $G^- \mapsto  \mu^{-1} G^-$.
	Applying these automorphisms, we get a family of $N=1$
	superconformal
	vectors. 
	By expanding the corresponding fields
	\begin{equation*}
			L(z) = \sum_{n \in \mathbb{Z}} L_n z^{-2-n}, \qquad
			G^\pm(z) = \sum_{n \in 1/2 + \mathbb{Z}} G^\pm_n
			z^{-3/2 
			-n}, \qquad
			J(z) = \sum_{n \in \mathbb{Z}} J_n z^{-1-n},
	\end{equation*}
	we get the commutation relations of the Virasoro operators $L_n$, and the following
	remaining commutation relations
	\begin{gather*}
			{[}J_m, J_n] = \frac{m}{3} \delta_{m,-n} c, \qquad
			{[}J_m, G^\pm_n] = \pm G^\pm_{m+n}, \qquad
			[G^\pm_m, L_n] = \left( m -\frac{n}{2}  \right)
			G^\pm_{m+n}, \\
			[L_m, J_n] = - n J_{m+n},  \qquad
			{[}G^+_m, G^-_n] = L_{m+n} +
			\frac{m-n}{2} J_{m+n} + \frac{c}{6} \left( m^2 -
			\frac{1}{4}
			\right) \delta_{m,-n}.  
	\end{gather*}
	
	\begin{subequations} 
	Sometimes it is convenient to introduce a different set of generating
	fields for this vertex algebra. We define $\tilde{L} = L - 1/2 \partial
	J$. This is a Virasoro field with central charge zero,
	namely
   		${[}\tilde{L}_\lambda \tilde{L}] = (\partial + 2\lambda)
		\tilde{L}$.
	With respect to this Virasoro element, $G^+$ is
	primary of conformal weight $2$ and $G^-$ is primary of conformal
	weight $1$; $J$ has conformal weight $1$ but is no longer
	a primary field. To summarize the commutation relations, we write
	\begin{equation}
		\begin{aligned}
			Q(z) = G^+(z) = \sum_{n \in \mathbb{Z}} Q_n
			z^{-2-n},& \qquad
			H(z) = G^-(z) = \sum_{n \in \mathbb{Z}} H_n
			z^{-1-n}, \\
			\tilde{L}(z)&= \sum_{n \in \mathbb{Z}} T_n
			z^{-2-n}. \\
		\end{aligned}
		\label{eq:n=2expansion}
	\end{equation}
		The corresponding $\lambda$-brackets of these fields are
		given by:
	\begin{equation}
		\begin{aligned}
			{[}\tilde{L}_\lambda \tilde{L}] &= (\partial + 2\lambda)
			\tilde{L}, & 
			{[}\tilde{L}_\lambda J] &= (\partial + \lambda)J -
			\frac{\lambda^2}{6}c, & 
			{[}\tilde{L}_\lambda Q] &= (\partial + 2\lambda) Q, \\
			{[}\tilde{L}_\lambda H] &= (\partial + \lambda) H, & 
			[H_\lambda Q] &= \tilde{L} - \lambda J +
			\frac{c}{6} \lambda^2. & &
		\end{aligned}
		\label{eq:agregado.w.n=2}
	\end{equation}
	The commutation relations between the Fourier coefficients are:
	\begin{equation}\begin{aligned}
			{[}T_m, T_n] &= (m-n) T_{m+n}, & [Q_m,Q_n] &= [H_m, H_n] = 0, \\ 
			{[} T_m, H_n] &= - n H_{m+n}, &
			[T_m , J_n] &= -n J_{m+n} - m (m+1)
			\frac{c}{12}\delta_{m, -n},\\  
			{[} T_m, Q_n] &= (m-n) Q_{m+n}, &
			[H_m, Q_n] &= T_{m+n} - m J_{m+n} + m
			(m-1)\frac{c}{6} \delta_{m, -n}  			
		\end{aligned}
		\label{eq:n=2newcommute}
	\end{equation}
\end{subequations}
\end{ex}
\begin{nolabel} \label{defn:n=2}
	In the subsequent sections, we will study in detail the structure
	theory of SUSY vertex algebras. Here we introduce some basic notation,
	used in the examples further on. 

	Let $V$ be a vector superspace over $\mathbb{C}$. Let $z$ be an even indeterminate and
	$\theta^1, \dots, \theta^N$ be odd anticommuting indeterminates which commute with $z$. For an ordered subset $I=(i_1,
	\dots, i_k) \subset \{1,\dots, N\}$, we will write
	$\theta^I =
	\theta^{i_1} \dots \theta^{i_k}$ and let $N\setminus I$ be the
	ordered complement of $I$ in $\{1,\dots, N\}$. 
	
	An $\End(V)$-valued
	superfield is a expression of the form:
	\begin{equation}
		A(z,\theta^1,\dots,\theta^N) = \sum_{(n|I):n \in
		\mathbb{Z}}  \theta^{N\setminus I}
		A_{(n|I)} z^{-1-n}
		\label{eq:definition_fourier}
	\end{equation}
	where $I$ runs over all ordered subsets of the set $\{1, \dots, N\}$, $A_{(n|I)} \in
	\End(V)$, and for each $I$ and $v \in V$
	we have $A_{(n|I)} v= 0$ for $n$ large enough. We will usually write $A(z,\theta)$ or
	simply $A(Z)$ for this field, 
	where $Z = (z, \theta^1, \dots, \theta^N)$. 
\end{nolabel}
\begin{rem}
	Define  an $N=2$ superconformal vertex algebra as a vertex
	algebra with a vector $\tau$ and two operators $S^1, S^2$ satisfying
	\begin{equation*}
		{[T}, S^i] = 0, \qquad [S^i, S^j] = 2 \delta_{i,j} T,
	\end{equation*}
	such that the corresponding fields $J(z) = -i Y(\tau, z)$, $L(z) =
	\tfrac{1}{2} Y(S^2 S^1 \tau, z)$ and
	\begin{equation}
		\begin{aligned}
			G^{(1)}(z) &\equiv G^+(z) + G^-(z) = -Y(S^2 \tau, z), \\
			G^{(2)}(z) & \equiv i \left(G^+(z) - G^-(z)\right) = Y(S^1
			\tau,
			z),\\
		\end{aligned}
		\label{eq:conf.def.n=2.2}
	\end{equation}
	satisfy the $\lambda$-brackets of Example
	\ref{ex:n=2}, $L_{-1} = T$, $G^{(i)}_{-1/2}  =S^i$, and $L_0$ is diagonalizable with eigenvalues bounded
	below. Then we obtain an $N_K = 2$ SUSY vertex algebra (see
	\ref{no:k.vertex.defin} for a definition) by letting 
	\begin{equation*}
		\ys(a,Z) = Y(a,z) + \theta^1 Y(S^1 a,z) + \theta^2 Y(S^2a, z) +  \theta^2
		\theta^1 Y(S^1 S^2 a, z).
		\label{n=2.k.superfields}
	\end{equation*}
	Similarly, given a vertex algebra with two vectors $\nu, \tau$ and an
	odd operator $S$ such that $[T,S] = 0$, $S^2 = 0$ and the
	associated fields:
	\begin{equation*}
		\begin{aligned}
			J(z) &= - Y(\tau,z), & H(z) &= Y(\nu, z) \\
			Q(z) &= Y(S \tau, z), & \tilde{L}(z) &= Y(S \nu, z) -
			 \partial_z J(z)
		\end{aligned}
	\end{equation*}
	satisfy the commutation relations
	(\ref{eq:agregado.w.n=2}), $T_{-1} = T$, $Q_{-1}= S$, $T_0$ is diagonalizable with eigenvalues
	bounded below, and $J_0$ is diagonalizable, we obtain an $N_W = 1$ SUSY
	vertex algebra (see \ref{defn:w_defin} for a definition)
	by letting:
	\begin{equation*}
		\ys(a,Z) = Y(a,z) + \theta Y(S a, z).
		\label{n=2.w.superfields}
	\end{equation*}\label{rem:n=2.construction}
\end{rem}
\begin{ex} \label{ex:n=2p}
	Following the previous remark, we can give the $N=2$ vertex algebra, as
	defined in Example \ref{ex:n=2}, the structure of an $N_K=2$ SUSY vertex
	algebra by letting $S^1 =
	(G^+_{(0)} + G^-_{(0)})$ and $S^2 = i (G^+_{(0)} - G^-_{(0)})$.
	Also we check directly that letting 
	\begin{equation}
		\tau =  \sqrt{-1} J_{(-1)} \vac,
		\label{eq:defintau}
	\end{equation}
	we get:
	\begin{equation}
	  \ys(\tau, z, \theta^i) = \sqrt{-1} J(z) +\theta^1
	  G^{(2)}(z) - \theta^2  G^{(1)}(z) + 2 \theta^1 \theta^2 L(z)
		\label{eq:tau.k.2.field}
	\end{equation}
	where $G^{(1)}(z) = G^+(z) + G^-(z)$ and $G^{(2)}(z) = i (G^+(z) -
	G^-(z))$. It follows that $[S^i, S^j] = 2 \delta_{ij} T$,  $\tau_{(0|0)} = 2 T$, $\tau_{(0|1)}= -
	S^1$ and $\tau_{(0|2)} = - S^2$ (cf. \ref{defn:k.conf.vertex.def}
	below).
	
	 We note that $G^{(i)}$ are primary of
	conformal weight $3/2$, and $J$ is primary of conformal weight $1$. The other
	commutation relations between the generating fields $L,\,J,\, G^{(i)}$
	$(i=1,2)$ are 
	\begin{equation*}
		\begin{aligned}
			{[}{G^{(i)}}_\lambda G^{(i)}] &= 2 L +
			\frac{c \lambda^2}{3}, &
			[{G^{(1)}}_\lambda G^{(2)}] &= -i \left( \partial
			+ 2\lambda
			\right) J, \\
			[J_\lambda G^{(1)}] &= - i G^{(2)}, &
			[J_\lambda G^{(2)}] &= i G^{(1)},
		\end{aligned}
	\end{equation*}
	or, equivalently,
	\begin{equation*}
		\begin{aligned}
			{[}G^{(i)}_m, G^{(i)}_n] &= 2 L_{m+n} + \left( m^2 -
		\frac{1}{4}
		\right) \frac{c}{3} \delta_{m, -n}, &
		{[}G^{(1)}_m, G^{(2)}_n] &= i \left( n - m \right) J_{m+n},\\
		[J_{m}, G^{(1)}_n] &= -i G^{(2)}_{m+n}, &
		[J_m, G^{(2)}_n] &= i G^{(1)}_{m+n}.
		\end{aligned}
	\end{equation*}

	Similarly, we can view the $N=2$ vertex algebra of Example \ref{ex:n=2} as
	an $N_W=1$ SUSY vertex algebra as follows. We will use the
	generating fields $\tilde{L}, Q, H,$ and $J$ with the commutation
	relations (\ref{eq:n=2newcommute}). Define the superfields:
	\begin{equation*}
		\ys(a, z, \theta) = Y(a,z) + \theta Y(Q_{-1}a, z),
	\end{equation*}
	and  let $T = T_{-1}$,
	$S=S^1=Q_{-1}$, so that $T$ and $S$ commute and $S^2 = 0$.  
	
	Note that defining the vectors $\nu = H_{(-1)}
	\vac$ and $\tau = - J_{(-1)} \vac$ we have in particular
	\begin{equation*}
		\begin{aligned}
			\ys(\nu, z, \theta) &= H(z) + \theta (\tilde{L}(z) + \partial_z J(z)), \\
		\ys(\tau, z, \theta) &= - J(z) + \theta Q(z).
	\end{aligned}
	\end{equation*}
	Therefore, if we consider the Fourier modes as defined in
	(\ref{eq:definition_fourier}), we have 
	\begin{equation*}
		\nu_{(0,0)} = T, \quad  \tau_{(0,0)} =  S.
	\end{equation*}
	Moreover, it is
	easy to see that the field $\tilde{L}(z) + \partial_z J(z)$ is also a Virasoro
	field and the conformal weights of the generating fields $\tilde{L},H,Q,J$ are
	positive with respect to this Virasoro field as well. It follows that
	the operator
	$\nu_{(1,0)}$ acts diagonally with non-negative eigenvalues (cf. Definition
	\ref{defn:w.conformal.definition} below).
\end{ex}
\begin{ex} \cite[Ex. 5.9d]{kac:vertex}\label{ex:2.23} 
	Consider the vertex algebra generated by a pair of free charged
	bosons $\alpha^\pm$ and a pair of free charged fermions
	$\varphi^\pm$ where the only non-trivial commutation relations are:
	\begin{equation*}
			{[}{\alpha^\pm}_\lambda \alpha^\mp] = \lambda, \qquad
			[{\varphi^\pm}_\lambda \varphi^\mp] = 1.
		\label{}
	\end{equation*}
	This vertex algebra contains the following family of $N=2$ vertex
	subalgebras ($\mu \in \mathbb{C}$):
	\begin{equation*}
		\begin{aligned}
			G^\pm &= :\alpha^\pm \varphi^\pm: \pm m \partial
			\varphi^\pm, \quad
			J = :\varphi^+ \varphi^-: -m (\alpha^+ + \alpha^-), \\
			L &= :\alpha^+ \alpha^-: + \frac{1}{2} : \partial
			\varphi^+ \varphi^-:  + \frac{1}{2} :\partial \varphi^-
			\varphi^+: - \frac{m}{2} \partial (\alpha^+ -
			\alpha^-).
		\end{aligned}
		\label{}
	\end{equation*}
	The vector $\tau$ given by (\ref{eq:defintau}) provides this vertex algebra
	with the structure of an
	$N_K=2$ SUSY vertex algebra, by letting $T=L_{-1}$ and $S^i
	= G^{(i)}_{-1/2}$ (see (\ref{eq:conf.def.n=2.2})). As in Example
	\ref{ex:n=2p}, we can view this vertex algebra as an $N_W=1$
	SUSY vertex algebra. 
\end{ex}
\begin{ex}
	An example of an  $N_W=N$ SUSY vertex algebra for each
	$N$ can be constructed as follows. Denote by $A$ the set of all ordered monomials
	$\theta^{i_1} \dots \theta^{i_s}$ and consider the
	superalgebra $\mathbb{C}[t, t^{-1},\theta^1, \dots, \theta^N]$
	where $t$ is even and $\theta^i$ are odd indeterminates. Let
	$W(1|N)$ be the Lie superalgebra of derivations of this superalgebra, and define the following
	collection of $W(1|N)$-valued formal
	distributions:
	\begin{equation*}
		\cF = \left\{ a^j(z) = \sum_{n \in \mathbb{Z}} (t^n a \partial_j)
		z^{-1-n} \Bigl| \; a \in A,\, j = 0,1,\dots ,N \right\},
		\label{}
	\end{equation*}
	where $\partial_j = \partial_{\theta^j}$ if $j > 0$ and
	$\partial_0 = \partial_t$. The pair $(W(1|N), \cF)$ is a 
	\emph{formal distribution Lie superalgebra}. The corresponding
	Lie conformal superalgebra is the free
	$\mathbb{C}[\partial]$-module $\cW_N$ generated by the vectors $a^j$, with $a
	\in A$ and $j = 0, \dots, N$, and the following $\lambda$-brackets (cf.
	\cite{kac:vertex}, \cite{kacfattori}):
	\begin{equation*}
		\begin{aligned}
			{[}{a^i}_\lambda b^j] &= (a \partial_i b)^j +
			(-1)^{p(a)}
			( (\partial_j a)b)^i, \qquad i, j \geq 1, \\
			{[}{a^i}_\lambda b^0] &= (a \partial_i b)^0 -
			(-1)^{p(b)}
			(ab)^i \lambda, \qquad
			[{a^0}_\lambda b^0] = - \partial (ab)^0  - 2
			(ab)^0 \lambda.
		\end{aligned}
		\label{}
	\end{equation*}
	Let $V(\cW_N)$ be the associated universal
	enveloping vertex algebra. The field 
	\begin{equation*}
		L(z) = -1^0(z) + \sum_{i=1}^n \partial_z (\theta^i)^i(z),
		\label{}
	\end{equation*}
	is a Virasoro field, and the elements
	$(\theta^i)^j$ are primary of conformal weight $1$, while the
	elements $-1^i$ are primary of conformal weight $2$. We will need
	later its Fourier modes, which are given by:
	\begin{equation*}
		L_n = - t^{n+1} \partial_t - (n+1) \sum t^{n} \theta^i
		\partial_{\theta^i}.
		\label{}
	\end{equation*}
	We define
	$T = L_{-1}=-\partial_t$. In
	order to be consistent with previous notation we define the
	fields $Q^i(z)=-1^i(z)$ and write down their Fourier modes which
	are
	\begin{equation*}
		Q^i_n = - t^{n+1} \partial_{\theta^i}.
		\label{}
	\end{equation*}
	In particular, we define
	$S^i = Q^i_{-1}$ for $i \geq 1$ and note that $(S^i)^2=0$ and $[T,S^i]=0$.

	In order to construct an $N_W=N$ SUSY vertex algebra from the vertex
	algebra $V(\cW_N)$ we proceed as before, 
	defining the superfields
	\begin{equation}
		\ys(a, z, \theta^1, \dots, \theta^N) = \sum_{I}
		(-1)^{\frac{I(I-1)}{2}} \theta^I Y( S^{i_1} \dots
		S^{i_s} a, z),
		\label{eq:definicionlarga}
	\end{equation}
	where the summation is taken over all ordered subsets $I=(i_1, \dots, i_s)$
	of $\{1, \dots, N\}$. It is
	straightforward to check that the data $(V(\cW_N), T, S^i, \vac, \ys)$ is indeed
	an $N_W=N$ SUSY vertex algebra. Moreover, this is a $W(1|N)$-conformal
	$N_W=N$ SUSY vertex algebra, as defined in \ref{no:intro.2aa}. We shall return to this example in
	\ref{ex:w_n.conformal}. \label{ex:w_n.series.aa}
\end{ex}
%
\begin{ex}
	We can similarly construct an $N_K=N$ SUSY vertex algebra $V(\cK_N)$ for any
	$N$. For this we define a subalgebra $K(1|N)$ of $W(1|N)$, of
	those differential operators preserving the form $\omega = dt +
	\sum \theta^i d\theta^i$ up to multiplication by a function (recall
	that we consider $d$ to be an even derivation, as in \cite{deligne2},
	but not in \cite{kacleur} and \cite{kacfattori}). It consists of differential
	operators of the form: 
	\begin{equation}
		D^f = f \partial_0 + \frac{1}{2} (-1)^{p(f)} \sum_{i = 1}^N (\theta^i
		\partial_0 + \partial_i)(f)(\theta_i \partial_0 + \partial_i),
		\label{eq:fattorikac1}
	\end{equation}
	for $f \in \mathbb{C}[t, t^{-1}, \theta^1, \dots, \theta^N]$. These operators satisfy 
	\begin{equation*}
		[D^f, D^g] = D^{\{f,g\}}, 
	\end{equation*}
	where 
	\begin{equation*}
		\{f,g\} = \left( f - \frac{1}{2}\sum_{i = 1}^N \theta^i \partial_i
		f\right) \partial_0 g - \partial_0 f \left( g -
		1\frac{1}{2} \sum_{i = 1}^N \theta^i \partial_i g \right) +
		(-1)^{f} \frac{1}{2} \sum_{i = 1}^N \partial_i f \partial_i g.
	\end{equation*}

	In
	particular $K(1|N)$ contains the operators
	\begin{equation*}
		\begin{aligned}
			L_n &= - t^{n+1} \partial_t - \frac{n+1}{2} t^n \sum
			\theta^i \partial_{\theta^i}, \qquad n \in \mathbb{Z},\\
			G^{(i)}_{n} &= -t^{n+1/2}(\partial_{\theta^i} - \theta^i
			\partial_t) + \left( n + \frac{1}{2} \right)
			t^{n-1/2} \theta^i \sum \theta^j
			\partial_{\theta^j}, \qquad n \in \frac{1}{2} +
			\mathbb{Z}.
		\end{aligned}
		\label{}
	\end{equation*}
	It is easy to see that the operators $L_n$ span a centerless Virasoro Lie
	algebra. 

	As in the $W(1|N)$ case, we construct the corresponding Lie conformal superalgebra as
	follows. It is the free $\mathbb{C}[\partial]$-module $\cK_N$ generated by
	vectors $a \in A$,
         with the following $\lambda$-brackets
	\cite{kacfattori}
	\begin{equation*}
		[a_\lambda b] = \left( \left( \frac{r}{2} - 1 \right) \partial (ab) + (-1)^r
		\frac{1}{2} \sum_{i=1}^n \partial_i a \partial_i b \right) + \lambda \left(
		\frac{r+s}{2} - 2 \right) ab,
		\label{fattorikac2}
	\end{equation*}
	where $a = \theta^{i_i} \dots \theta^{i_r}$, $b = \theta^{j_1} \dots
	\theta^{j_s}$. 
	
	We denote by $V(\cK_N)$ its universal enveloping vertex
	algebra, and we define the operators
	$T = L_{-1}$ and $S^i = G^{(i)}_{-1/2}$. Now we define the
	state-field correspondence as in (\ref{eq:definicionlarga}):
	\begin{equation*}
		\ys(a, z, \theta) = \sum_{I} (-1)^{\frac{I(I-1)}{2}} \theta^I Y(S^I a, z).
		\label{}
	\end{equation*}
	All the properties of an $N_K=N$ SUSY vertex algebra are
	straightforward to check as in the previous cases. Moreover, this is a
	$K(1|N)$-conformal $N_K=N$ SUSY vertex algebra. We will return to this
	example in
	\ref{ex:k.n.series}. \label{ex:k_series.begin} 
\end{ex}
\numberwithin{thm}{subsection}
\section{Structure theory of $N_W=N$ SUSY vertex algebras}\label{sec:structure}
In this section we develop the structure theory of SUSY Lie conformal algebras and SUSY vertex
algebras along the lines of \cite{kac:vertex} (see also \cite{kacdesole2} for a better
exposition). Proofs are rather straightforward adaptations of those in the vertex algebra
case, the only difficulty being the problem of signs.
\subsection{Formal distribution calculus} \label{sub:formal}
\begin{nolabel} 
In what follows we fix the ground field to be the complex numbers $\mathbb{C}$ and
$N$ to be a non-negative integer.
Let $\theta^1, \dots, \theta^N$ be Grassmann variables and $I=\{i_1, \dots, i_k\}$
be an \emph{ordered $k$-tuple: $1 \leq i_1 <\dots < i_k \leq N$}. We will denote
\begin{equation*}
\theta^I = \theta^{i_1} \dots \theta^{i_k}, \qquad \theta^N = \theta^1 \dots
\theta^N.
\end{equation*}
For an element $a$ in a vector superspace we will
denote $(-1)^a=(-1)^{p(a)}$, where $p(a) \in \mathbb{Z}/2 \mathbb{Z}$ is the parity of $a$,
and, given a $k$-tuple $I$ as above, we will let
$(-1)^I = (-1)^k$.
Given two disjoint ordered tuples $I$ and $J$, we define $\sigma(I,J) = \pm 1$ by 
\begin{equation*}
	 \theta^I \theta^J = \sigma(I,J)\theta^{I \cup J},
\end{equation*}
and we define $\sigma(I,J)$ to be zero if $I \cap J \neq \emptyset$. Also, unless
noted otherwise, all ``union'' symbols ``$\cup$'' will denote disjoint
unions\footnote{This will not be true in section \ref{sub:k.formal} where we
analyze $N_K=n$ SUSY vertex algebras}. It
follows easily, by looking at $\theta^I \theta^J \theta^K$,  that for three mutually disjoint tuples, $I,\, J$ and $K$ we have:
\begin{equation}
	\sigma(I,J)\sigma(I\cup J, K) = \sigma(I,J\cup K)\sigma(J,K), \quad
	\sigma(I,J) = (-1)^{IJ} \sigma(J,I).
	\label{eq:not.2b}
\end{equation}
Here and further $(-1)^{IJ}$ stands for $(-1)^{(\sharp I)(\sharp J)}$. 

We will denote by $N\setminus I$ the ordered complement of $I$ in $\{1, \dots, N\}$
and define $\sigma(I) := \sigma(I, N\setminus I)$. 
It follows from the definitions that
$	\theta^I \theta^{N\setminus I} = \sigma(I) \theta^N $.
\label{no:not}
\end{nolabel}
\begin{nolabel}
	Let $Z=(z,\theta^1,\dots,\theta^N)$ and $W=(w,\zeta^1,\dots, \zeta^N)$
	denote two sets of coordinates in the formal superdisk $D=D^{1|N}$. As
	before, all $\theta^i$ and $\zeta^j$ anticommute.
	
	Let $\mathbb{C}[ [z]]$ be the algebra of formal power series in $z$; its
	elements are are series $\sum_{n \geq 0} a_n z^n$ with $a_n \in \mathbb{C}$. 
	The superalgebra of regular functions
	in $D$ is defined as  $\mathbb{C}[ [Z]] := \mathbb{C}[ [z]] \otimes
	\mathbb{C}[\theta^1,\dots, \theta^N]$. Similarly, we define the superalgebra $\mathbb{C}[ [Z, W]] :=
	\mathbb{C}[ [z, w]] \otimes \mathbb{C}[\theta^1, \dots, \theta^N, \zeta^1, \dots,
	\zeta^N]$. 
	 
	For any $\mathbb{C}$-algebra $R$, we denote by $R( (z))$ the algebra of $R$-valued formal
	Laurent series, its elements are series of the form $\sum_{n \in \mathbb{Z}} a_n z^n$
	such that $a_n \in R$ and there exists $N_0 \in \mathbb{Z}$ such that $a_n = 0$ for all
	$n \leq N_0$. If $R$ is a field, so is $R( (z))$. We denote $R( (Z)) := R( (z))
	\otimes_\mathbb{C} \mathbb{C}[\theta^1, \dots, \theta^N]$. Denote also by $\mathbb{C}(
	(Z))( (W))$ the superalgebra $R( (W))$ where $R = \mathbb{C}( (Z))$; its elements are
	Laurent series in $W$ whose coefficients are Laurent series in $Z$. Similarly we have the
	superalgebra $\mathbb{C}( (W))( (Z))$. 

	Denote by $\mathbb{C}( (z,w))$ the field of fractions of $\mathbb{C}[ [z,w]]$ and put
	$\mathbb{C}( (Z, W)) := \mathbb{C}( (z,w)) \otimes_{\mathbb{C}} \mathbb{C}[\theta^1,
	\dots, \theta^N, \zeta^1, \dots, \zeta^N]$. One may think of this
	superalgebra as the algebra of
	meromorphic functions in the formal superdisk $D^{2|2 N}$.  Given such a meromorphic
	function, we can ``expand it near
	the $w$ axis'', to obtain an element of $\mathbb{C}( (Z))( (W))$. Indeed,  $\mathbb{C}[ [z,w]]$ embeds
	naturally in $\mathbb{C}( (z))( (w))$ and $\mathbb{C}( (w))( (z))$ respectively. Since
	$\mathbb{C}( (z,w))$ is the ring of fractions of $\mathbb{C}[ [z,w]]$ and $\mathbb{C}(
	(z))( (w))$ and $\mathbb{C}( (w))( (z))$ are fields, these embedding induce respective
	algebra embeddings
	\begin{equation*}
		\mathbb{C}( (z))( (w)) \overset{i_{z,w}}{\hookleftarrow} \mathbb{C}( (z,w))
		\overset{i_{w,z}}{\hookrightarrow} \mathbb{C}( (w))( (z)).
		\label{}
	\end{equation*}
	(A concrete example is given by (\ref{eq:concrete_expansion})).
	
	Tensoring with the corresponding Grassmann superalgebras, we obtain
	superalgebra embeddings
	\begin{equation*}
		\mathbb{C}( (Z))( (W)) \overset{i_{z,w}}{\hookleftarrow} \mathbb{C}( (Z,W))
		\overset{i_{w,z}}{\hookrightarrow} \mathbb{C}( (W))( (Z)).
		\label{}
	\end{equation*}

	Let $\cU$ be a vector superspace. An $\cU$-valued \emph{formal
	distribution} is an expression of the form 
	\begin{equation*}
		a(Z) = \sum_{(n|I): n \in \mathbb{Z}} Z^{n|I} a_{n|I}, \qquad
		a_{n|I} \in \cU.
	\end{equation*}
	The space of such distributions will be denoted $\cU[ [Z, Z^{-1}]]$. 
	We denote by
	$\mathbb{C}[Z, Z^{-1}] := \mathbb{C}[z,z^{-1}] \otimes \mathbb{C}[\theta^1, \dots, \theta^N]$ the
	superalgebra of Laurent polynomials. A $\cU$-valued formal distribution is canonically
	a linear functional
	$\mathbb{C}[Z, Z^{-1}]
	\rightarrow \cU$. To see this, we define the \emph{super residue} as the
	coefficient of $Z^{-1|N}$:
	\begin{equation*}
		\res_Z a(Z) = a_{-1|N}.
	\end{equation*}
	This clearly satisfies 
	\begin{equation}
		\res_Z \partial_z a(Z) = \res_Z \partial_\theta
	a(Z) = 0.
	\label{eq:nseporque}
	\end{equation}
	Given a $\cU$-valued formal distribution $a(Z)$ we obtain a linear map
	$\mathbb{C}[Z, Z^{-1}] \rightarrow \cU$ given by
	\begin{equation*}
		f(Z) \mapsto \res_Z a(Z) f(Z).
	\end{equation*}
	Conversely, every formal distribution arises in this way. Indeed we have:
	\begin{equation}
		\res_Z Z^{n|I} a(Z) = \sigma(I) a_{-1-n|N\setminus
		I}.
		\label{eq:distrib.4}
	\end{equation}
	Therefore the formal distribution $a(Z)$ can be written as 
	\begin{equation*}
		a(Z) = \sum_{(n|I): n \in \mathbb{Z}} Z^{-1-n|N\setminus I}
		a_{(n|I)},
	\end{equation*}
	where 
	\begin{equation*}
		a_{(n|I)} = \sigma(I) \res_Z Z^{n|I} a(Z).
	\end{equation*}
	\label{no:distrib.6}

	We can similarly define $\cU$-valued formal distributions in two variables, as
	expressions of the
	form 
	\begin{equation*}
		a(Z, W) = \sum_{(j|J), (k|K)} Z^{j|J}W^{k|K} a_{j|J, k|K},  \qquad a_{j|J,
		k|K} \in \cU.
		\label{}
	\end{equation*}
	The space of such formal distributions will be denoted $\cU[ [Z, Z^{-1}, W, W^{-1}] ]$. 

	Note that in the
	case $\cU = \mathbb{C}$, both algebras $\mathbb{C}( (Z))( (W))$ and $\mathbb{C}( (W))(
	(Z))$ are embedded in $\mathbb{C}[ [ Z, Z^{-1}, W, W^{-1}]]$. We
	will denote by $i_{z, w}$ and $i_{w, z}$ 
	the corresponding embeddings of $\mathbb{C}( (Z, W))$ in $\mathbb{C}[ [Z, Z^{-1}, W,
	W^{-1}]]$. When $f(Z,W)$ is a Laurent polynomial (that is a polynomial in $z,\,
	z^{-1},\, w, \, w^{-1}$ and the odd variables) then the embeddings $i_{z,w} f$ and
	$i_{w, z} f$ coincide on $\mathbb{C}[ [Z, Z^{-1}, W, W^{-1}]]$. Indeed, it is immediate to see
	that 
	\begin{equation}
		\mathbb{C}( (Z))( (W)) \cap \mathbb{C}( (W))( (Z)) = \mathbb{C}[ [Z, W]][
		z^{-1}, w^{-1}],
		\label{eq:intersection_added}
	\end{equation}
	where the intersection is taken in $\mathbb{C}[ [Z, Z^{-1}, W, W^{-1}]]$. The images under
	these embeddings 
	are different for other
	functions, as we will see below (cf.
	\ref{prop:delta_function.1}). 

	A $\cU$-valued formal distribution in two
	variables is called \emph{local} if there exists a non-negative integer $n$ such that 
	\begin{equation*}
		(z-w)^n a(Z,W) = 0.
	\end{equation*}
\end{nolabel}
\begin{nolabel}
	Note that the differential operators $\partial_z, \partial_{\theta^i}$ and $\partial_w,
	\partial_{\zeta^i}$ act in the usual way on the
	spaces $\mathbb{C}( (Z, W))$, $\mathbb{C}[ [Z, Z^{-1}, W, W^{-1}]]$. For $j \in
	\mathbb{Z}_+$ and $J= (j_1, \dots, j_k)$ we will denote 
	\[\partial_Z^{j|J} = \partial_z^j \partial_{\theta^{j_1}} \dots
	\partial_{\theta^{j_k}}.\] We define
	\begin{equation*}
		\partial_Z^{(j|J)} := \frac{(-1)^{\frac{J(J+1)}{2}}}{j!}
		\partial_Z^{j|J}, \qquad
		Z^{(j|J)} := \frac{(-1)^{\frac{J(J+1)}{2}}}{j!} Z^{j|J}.
		\label{}
	\end{equation*}
	One checks easily that the embeddings $i_{z,w}$ and $i_{w,z}$ defined above, commute
	with the action of the differential operators $\partial_Z^{j|J}$ and
	$\partial_W^{j|J}$. 
	
	We will denote 
	\begin{equation}
	\begin{gathered}
		Z-W
	= (z-w, \theta^1 - \zeta^1, \dots, \theta^N - \zeta^N), \qquad Z^{n|I} = z^n
	\theta^I, \\ (Z-W)^{j|J} = (z-w)^j \prod_{i \in J} (\theta^i - \zeta^i),
	\quad \partial_W = (\partial_w, \partial_{\zeta^1}, \dots,
	\partial_{\zeta^N}),
	\label{eq:w.intro.recontra}
	\end{gathered}
\end{equation}
	and for any $\cU$-valued formal distribution $f(Z)$, we define its
	\emph{Taylor expansion} as:
	\begin{equation}
		f(Z) = e^{(Z-W) \partial_W} f(W),
		\label{eq:taylor_expansion_defn}
	\end{equation}
	where \[(Z-W) \partial_W = (z-w) \partial_w + \sum_i (\theta^i - \zeta^i)
	\partial_{\zeta^i}.\] Expanding the exponential in (\ref{eq:taylor_expansion_defn}) we
	obtain:
	\begin{equation*}
		f(Z) = \sum_{(j|J): j \geq 0} (-1)^J (Z-W)^{j|J} \partial_W^{(j|J)} f(W).
		\label{}
	\end{equation*}
	\label{no:Taylor_1_def}
\end{nolabel}
\begin{rem}
	In the definition of formal distributions and super residues, we can replace $\mathbb{C}$ by any
	commutative superalgebra $\cA$, and $\cU$ by any $\cA$-module. We see
	immediately that the residue map is of parity $N \mod 2$, that is, for
	$\chi \in \cA$, and $u(Z)$ an $\cU$-valued distribution, we have:
	\begin{equation*}
		\res_Z \chi u(Z) = (-1)^{\chi N} \res_Z u(Z).
	\end{equation*}
	On the other hand, this residue map is a morphism of right $\cA$-modules,
	namely:
	\begin{equation*}
		\res_Z u(Z) \chi = \Bigl( \res_Z u(Z) \Bigr) \chi.
	\end{equation*}
	\label{rem:parity_of_res}
\end{rem}
\begin{prop}
	There exists a unique $\mathbb{C}$-valued formal distribution $\delta(Z,W)$ such that
	for every
	function $f \in \cU[Z Z^{-1}]$ we have $\res_Z \delta(Z,W) f(Z) = f(W)$. 
	\label{prop:delta_function.1}
\end{prop}
\begin{proof}
	For uniqueness, let $\delta$ and $\delta'$ be two such distributions, then
	$\beta = \delta - \delta'$ satisfies $\res_Z \beta(Z,W) f(Z) = 0$ for all
	functions $f(Z)$. Decomposing $\beta(Z,W) = \sum
	\beta_{n|I,m|J} W^{m|J} Z^{n|I} $, and multiplying by $Z^{k|L}$ we see that 
	$\beta_{-1-k|N-L, m|J} = 0$ for all $m|J$, hence $\beta = 0$.
	Existence will be proved below.
\end{proof}
\begin{nolabel}
	We define the formal $\delta$-function as the $\mathbb{C}$-valued formal
	distribution in two variables, given by
	\begin{equation}
		\delta(Z,W) = (i_{z,w} - i_{w,z}) (Z-W)^{-1|N} = (i_{z,w} -
		i_{w,z}) \frac{(\theta - \zeta)^N}{z-w} 
		\label{eq:delta.1}
	\end{equation}
	It follows that
	\begin{equation*}
		\partial_w^{(n)} \delta(Z,W) := \frac{1}{n!}\partial_w^n \delta(Z,W) =  (i_{z,w} - i_{w,z}) (Z-W)^{-1-n|N}.
	\end{equation*}
	This distribution has the following properties:
	\begin{enumerate}
		\item $(Z-W)^{m|J} \partial_W^{n|I} \delta (Z,W) = 0$ whenever $m >
			n$ or $J \supsetneq I$, 
		\item $(Z-W)^{j|J} \partial_{W}^{(n|I)} \delta(Z,W) =
			\sigma(I \setminus J,J) 
			\partial_W^{(n-j| I\setminus J)} \delta(Z,W)$ if $n \geq j
			$ and $I \supset J$, 
		\item $\delta(Z,W) = (-1)^N \delta(W,Z)$,
		\item $\partial_Z^{j|J} \delta(Z,W) = (-1)^{j+N+J} \partial_W^{j|J}
			\delta(W,Z)$,
		\item $\delta(Z,W) a(Z) = \delta(Z,W) a(W)$, where $a(Z)$ is any
			formal distribution,
		\item $\res_Z \delta(Z,W) a(Z) = a(W)$, where $a(Z)$ is any formal
			distribution,
		\item $\exp \left((Z-W)\Lambda \right) \partial_W^{n|I}
			\delta(Z,W) =
			(\Lambda + \partial_W)^{n|I} \delta(Z,W)$, where $\Lambda
			= (\lambda, \chi^1, \dots, \chi^N)$, $\chi^i$ are odd
			anticommuting variables, $\lambda$ is even, $\lambda$
			commutes with $\chi^i$, and we  write 
			\begin{equation*}
				\begin{aligned}
					(Z-W)\Lambda = (z-w)\lambda + \sum_i
					(\theta^i - \zeta^i) \chi^i, \qquad
					(\Lambda + \partial_W) = (\lambda +
					\partial_w, \chi_i +
					\partial_{\theta^i}).
				\end{aligned}
			\end{equation*}
	\end{enumerate}
	\begin{proof}
		Writing $\partial_\zeta^I = \partial_{\zeta^{i_1}} \dots
		\partial_{\zeta^{i_k}}$ we have 
		\begin{equation*}
			(Z-W)^{m|J} \partial_W^{n|I} \delta (Z,W) = (z-w)^m n!
			(i_{z,w} - i_{w,z}) (z-w)^{-1-n} (\theta - \zeta)^J
			\partial_\zeta^I (\theta - \zeta)^N
		\end{equation*}
		Now this clearly vanishes if $m \geq 1+n$ since then the two
		embeddings $i_{z,w}$ and $i_{w,z}$ coincide on the regular function
		$(z-w)^{m-n-1}$. The other factor is clearly zero if $J \supsetneq
		I$ since for every $j \in J \setminus I$ we have a factor
		$(\theta^j
		- \zeta^j)$ in $\partial_\zeta^I (\theta - \zeta)^N$. This proves
		(1). 
		
		In order to prove (2) we write:
		\begin{multline}
			(Z-W)^{j|J} \partial_{W}^{n|I} \delta(Z,W) = n!(i_{z,w} -
			i_{w,z}) (z-w)^{j-1-n} (\theta - \zeta)^J \partial_\zeta^I
			(\theta-\zeta)^N \\= \frac{n!}{(n-j)!} (n-j)! (i_{z,w} - i_{w,z})
			(z-w)^{-1-(n-j)} (-1)^{J} \sigma(J,
			I\setminus J) (\theta - \zeta)^J \partial_\theta^{J}
			\partial_\zeta^{I\setminus J} (\theta-\zeta)^N  \\ =
			\frac{n!}{(n-j)!} 
			(-1)^{\frac{J(J+1)}{2}} \sigma(J,
			I \setminus J) 
			 (n-j)! (i_{z,w} - i_{w,z}) (z-w)^{-1-(n-j)}
			 \partial_\zeta^{I\setminus J} (\theta-\zeta)^N  \\
			 =(-1)^{\frac{J(J+1)}{2}}
			\sigma(J, I \setminus J) \frac{n!}{(n-j)!}
			\partial_W^{n-j| I\setminus J} \delta(Z,W).
		\end{multline}
		(3) is obvious and (4) follows from (3) easily. In order to prove
		(5) wee see that from (1) we have $\delta z = \delta w$ therefore
		we get $\delta(Z,W) z^k = \delta(Z,W) w^k$. On the other hand, also from (1)
		it follows that $\delta(Z,W) \theta^i = \delta(Z,W) \zeta^i$. Hence
		$\delta(Z,W)
		\theta^I = \delta(Z,W) \zeta^I$ and we have proved that $\delta(Z,W) Z^{n|I}
		= \delta(Z,W) W^{n|I}$. The result follows by linearity now.

		(6) follows by taking residue in (5). To prove (7) we first expand
		the exponential in power series:
		\begin{equation}
			\exp\left( (Z-W)\Lambda \right) \partial^{n|I}_W
			\delta(Z,W) = \sum_{(j|J): j \geq 0} (-1)^J (Z-W)^{(j|J)}
			\Lambda^{j|J}
			\partial^{n|I}_W \delta(Z,W).
			\label{eq:delta.4}
		\end{equation}
		Now using (2) we see that this is:
		\begin{equation}
			\sum_{(j|J): j \geq 0} \binom{n}{j} \Lambda^{j|J} \sigma(J, I \setminus J)
			\partial_W^{n-j|I \setminus J} \delta(Z,W).
			\label{eq:delta.5}
		\end{equation}
		On the other hand we can expand the right hand side of (7) as:
		\begin{multline*}
			(\Lambda + \partial_W)^{n|I} = (\lambda + \partial_w)^n
			(\chi + \partial_\zeta)^I  = \sum_{(j|J): j \geq 0} \binom{n}{j} \lambda^j
			\partial_w^{n-j} \sigma(J, I\setminus J) \chi^J
			\partial_\zeta^{I\setminus J}  \\
			= \sum_{(j|J): j \geq 0} \binom{n}{j} \Lambda^{j|J} \sigma(J, I\setminus J)
			\partial_W^{n-j| I \setminus J}.
		\end{multline*}
		Comparing with (\ref{eq:delta.5}) we get the result.
	\end{proof}
	\label{no:delta}
\end{nolabel}
\begin{lem}
	Let $a(Z,W)$ be a local formal distribution in two variables. Then $a(Z,W)$ can be
	uniquely decomposed as
	\begin{equation}
		a(Z,W) = \sum_{(j|J):j \geq 0} 	\left( \partial_W^{(j|J)}
		\delta(Z,W) \right) c_{j|J}(W), 
		\label{eq:decomp.1}
	\end{equation}
	where the sum is finite. The coefficients $c_{j|J}$ are given by
	\begin{equation}
		c_{j|J}(W) = \res_Z (Z-W)^{j|J} a(Z,W).
		\label{eq:decomp.2}
	\end{equation}
	\label{lem:decomp}
\end{lem}
\begin{proof}
 	First we note that if $a(Z,W)$ is local then the sum on the right hand side is
	finite.
	Let $b(Z,W)$ the difference between the right hand side and the left hand
	side of (\ref{eq:decomp.1}). We find:
	\begin{equation*}
		\begin{aligned}
			\res_Z (Z-W)^{k|K} b(Z,W) &= \res_Z (Z-W)^{k|K} a(Z,W) -
			\\ 
			& \quad - \res_Z \sum_{(j|J): j \geq 0}  
			  (Z-W)^{k|K}
			\left( \partial_W^{(j|J)} \delta(Z,W) \right) c_{j|J}(W) \\
			& = c_{k|K}(W) - \res_Z \left(
			\partial_W^{(j-k|J\setminus K)}
			\delta(Z,W)\right)c_{(j|J)}(W)\\
			&= c_{k|K}(W) - \res_Z \delta(Z,W)c_{k|K}(W) = 0,
		\end{aligned}
	\end{equation*}
	where in the second line we have used (2) of \ref{no:delta}. It follows
	 that $b(Z,W)$ has no negative powers of $z$. Moreover, $b(Z,W)$ is local,
	since $a(Z,W)$ is, and the right hand side of (\ref{eq:decomp.1}) is
	local by (1) of \ref{no:delta}. We can write then 
	\begin{equation*}
	  b(Z,W) = \sum_{(j|J): j \geq 0} Z^{j|J} b_{j|J}(W),
	\end{equation*}
	and since $(z-w)^n b(Z,W) = 0$ we obtain:
	\begin{equation*}
	  \sum_{\stackrel{(j|J)}{j \geq k \geq 0}} \binom{n}{k} Z^{j|J} w^{n-k}
	  b_{j-k|J}(W) = 0,
	\end{equation*}
 	which easily shows that $b(Z,W)=0$. 	
	 Uniqueness is clear by taking
	residues on both sides of (\ref{eq:decomp.1}).
\end{proof}
\begin{nolabel}
	Let $a(Z,W)$ be a
	formal distribution in two variables. We define its formal Fourier transform by:
	\begin{equation}
		\cF^\Lambda_{Z,W} a(Z,W) = \res_Z \exp \left( (Z-W)\Lambda \right)
		a(Z,W),
		\label{eq:fourier.1}
	\end{equation}
	where $\Lambda=(\lambda, \chi^1, \dots, \chi^N)$, $\lambda$ is an even
	variable, and $\chi^i$ are odd anticommuting
	variables, commuting with $\lambda$. 

	Expanding this exponential we have (recall (\ref{eq:delta.4})) :
	\begin{multline}
			\cF^\Lambda_{Z,W} a(Z,W) = \res_Z \sum_{(j|J): j \geq 0}
			\Lambda^{j|J}
			(Z-W)^{(j|J)} a(Z,W) \\
			= \sum_{(j|J): j \geq 0} (-1)^{JN} \Lambda^{j|J}
			\res_Z (Z-W)^{(j|J)} a(Z,W) 
			= \sum_{(j|J): j \geq 0} (-1)^{JN} 
			\Lambda^{(j|J)} c_{j|J}(W) 
		\label{eq:fourier.2}
	\end{multline}
	where $c_{j|J}$ are defined by (\ref{eq:decomp.2}) and we write, as before
	\begin{equation*}
		\Lambda^{j|J} = \lambda^j \chi^{j_1} \dots \chi^{j_k}, \quad \Lambda^{(j|J)} := \frac{(-1)^{\frac{J(J+1)}{2}}}{j!} \Lambda^{j|J}.
		\label{}
	\end{equation*}
	\label{no:fourier}
\end{nolabel}
\begin{prop}
	The formal Fourier transform satisfies the following properties:
	\begin{enumerate}
		\item sesquilinearity:
			\begin{equation*}
				\begin{aligned}
					\cF^\Lambda_{Z,W} \partial_z a(Z,W) &=
					-\lambda \cF^\Lambda_{Z,W} a(Z,W) =
					[\partial_w, \cF^\Lambda_{Z,W}] a(Z,W), \\
					\cF^\Lambda_{Z,W} \partial_{\theta^i}
					a(Z,W) &= - (-1)^N \chi^i \cF^\Lambda_{Z,W}
					a(Z,W) =
					(-1)^N [\partial_{\zeta^i},
					\cF^\Lambda_{Z,W}]a(Z,W).
				\end{aligned}
			\end{equation*}	
		\item For any local formal distribution $a(Z,W)$ we have:
			\begin{equation}
				\begin{aligned}
					(-1)^N \cF^\Lambda_{Z,W}a(W,Z) &=
					\cF^{-\Lambda-\partial_W}_{Z,W} a(Z,W),\\
					&= \cF^\Gamma_{Z,W} a(Z,W)|_{\Gamma =
					-\Lambda - \partial_W},
				\end{aligned}	
				\label{eq:fourierprop.2}
			\end{equation}
			where $- \lambda - \partial_W = (-\lambda - \partial_w, -
			\chi^i - \partial_{\zeta^i})$. 
		\item For any formal distribution in three variables $a(Z,W,X)$ we have
			\begin{equation*}
				\cF^\Lambda_{Z,W} \cF^\Gamma_{X,W} a(Z,W,X) =
				(-1)^N
				\cF^{\Lambda + \Gamma}_{X,W} \cF^\Lambda_{Z,X}
				a(Z,W,X),
			\end{equation*}
			where $\Gamma = (\gamma, \eta^1, \dots, \eta^N)$, with
			$\eta^i$ odd anticommutative variables and $\gamma$ is
			even and commutes with $\eta^i$, $\Lambda + \Gamma$ is
			the sum $(\lambda + \gamma, \chi^i + \eta^i)$, and
			the superalgebra
			$\mathbb{C}[\Lambda, \Gamma]$ is commutative.
	\end{enumerate}
	\label{prop:fourier}
\end{prop}
\begin{proof}
	The proof of the first equality of both lines of (1) follows from
	(\ref{eq:nseporque}). For the first equality of the second line 
	we have
	\begin{multline*}
			\cF^\Lambda_{Z,W} \partial_{\theta^i} a(Z,W) = \res_Z \exp
			\left( (Z-W)\Lambda \right) \partial_{\theta^i} a(Z,W) \\
			= - \res_Z \left( \partial_{\theta^i} \exp\left(
			(Z-W)\Lambda
			\right) \right) a(Z,W) 
			= - \res_Z \chi^i \exp\left( (Z-W)\Lambda \right) a(Z,W) \\
			= - (-1)^N \chi^i \cF^\Lambda_{Z,W} a(Z,W).
	\end{multline*}
	For the second equality of the second line of (1) we have:
	\begin{multline*}
			{[}\partial_{\zeta^i}, \cF^\Lambda_{Z,W}] a(Z,W) =
			(-1)^N \left( \res_Z \partial_{\zeta^i} \exp\left( (Z-W)
			\Lambda
			\right) a(Z,W) - \right. \\  \quad \left. - \exp\left( (Z-W)
			\Lambda 
			\right) \partial_{\zeta^i} a(Z,W) \right) = 
			= (-1)^N \res_Z \left( \partial_{\zeta^i} \exp \left(
			(Z-W)\Lambda
			\right) \right) a(Z,W) \\ 
			= - \chi^i \cF^\Lambda_{Z,W} a(Z,W).
	\end{multline*}
	To prove (2) it is enough, by Lemma \ref{lem:decomp}, to prove the
	statement when $a(Z,W) = \left(\partial^{j|J}_W \delta(Z,W) \right) c(W)$.
	In this case we have:
	\begin{equation*}
		\begin{aligned}
			\cF^\Lambda_{Z,W} a(W,Z) &= \cF^\Lambda_{Z,W} \left(
			\partial_Z^{j|J} \delta(W,Z)
			\right) c(Z) 
			= \cF^\Lambda_{Z,W} (-1)^{j+J+N}  \partial_W^{j|J}
			\delta(Z,W) c(Z).
		\end{aligned}
	\end{equation*}
	Now using (7) in \ref{no:delta} we can express the last expression
	above as:
	\begin{multline*}
		(-1)^{j+J+N} \res_Z (\Lambda + \partial_W)^{j|J} \delta(Z,W) c(Z)
		=\\=
		(-1)^{j + J + JN+N} (\Lambda + \partial_W)^{j|J} \res_Z \delta(Z,W)
		c(Z)  = (-1)^{j+J+JN+N} (\Lambda + \partial_W)^{j|J} c(W).
	\end{multline*}
	On the other hand we have 
	\begin{multline*}
		\cF^\Gamma_{Z,W} a(Z,W)|_{\Gamma=-\Lambda-\partial_W}  = (-1)^{JN}
		(-\Lambda - \partial_W)^{j|J} c(W)  \\ =(-1)^{j+J+JN} (\Lambda +
		\partial_W)^{j|J} c(W).
	\end{multline*}
	The proof of (3) is straightforward:
	\begin{multline}
		\cF^\Lambda_{Z,W} \cF^\Gamma_{X,W} = \res_Z \exp\left( 
		(Z-W) \Lambda 
		\right) \res_X \exp \left( (X-W)\Gamma \right) = \\ = \res_Z \res_X
		\exp\left( (Z-W)\Lambda + (X-W)\Gamma \right) = \\ = (-1)^{N} \res_X
		\res_Z \exp \left( (Z-X)\Lambda  + (X-W)(\Lambda+\Gamma) \right) =
		\\ = (-1)^N \res_X \exp\left( (X-W)(\Lambda + \Gamma) \right) \res_Z
		\exp\left( (Z-X)\Lambda \right) = (-1)^N \cF^{\Lambda +
		\Gamma}_{X,W} \cF^\Lambda_{Z,X}.
		\label{eq:fourierproof.6}
	\end{multline}
	The sign $(-1)^N$ appears when we commute the residue maps (recall that they have parity $N
	\mod 2$). 
\end{proof}
\subsection{$N_W=N$ SUSY Lie conformal algebras} \label{sub:liesuper}
\begin{nolabel}
	Let $\fg$ be a Lie superalgebra. A pair of $\fg$-valued formal
	distributions $a(Z),
	b(Z)$ is called \emph{local} if the distribution $[a(Z),b(W)]$ is local. By
	the decomposition Lemma \ref{lem:decomp} we have for such a pair:
	\begin{equation*}
		[a(Z),b(W)] = \sum_{(j|J): j \geq 0}
		\left(\partial_W^{(j|J)} \delta(Z,W) \right) c_{j|J}(W).
	\end{equation*}
	where 
	\begin{equation*}
	  c_{j|J} (W) = \res_Z (Z-W)^{j|J} [a(Z),b(W)].
	\end{equation*}

	We define $a(W)_{(j|J)}b(W) = c_{j|J}(W)$ and we call this
	operation the $(j|J)$-product. Let us also define
	the $\Lambda$-bracket of two $\fg$-valued formal distributions by 
	\begin{equation}
		[a_\Lambda b](W) = \cF^\Lambda_{Z,W} [a(Z),b(W)],
		\label{eq:lie.2}
	\end{equation}
	where $\cF_{Z,W}^\Lambda$ is the formal Fourier transform defined in
	\ref{no:fourier}.  
	It follows from the definitions and from (\ref{eq:fourier.2}) that
	\begin{equation*}
		[a_\Lambda b] = \sum_{(j|J):j \geq 0}
		(-1)^{JN} \Lambda^{(j|J)} a_{(j|J)}b.
	\end{equation*}
	\label{no:lie}
	Note also that the $\Lambda$-bracket has parity $N \mod 2$ (this
	follows from the fact that the residue map has parity $N \mod 2$).

	A pair $(\fg, \cR)$ consisting of a Lie superalgebra $\fg$ and a family
	$\cR$ of
	pairwise local $\fg$-valued formal distributions $a(Z)$, whose coefficients span $\fg$, stable under all $j|J$-th products and
	under the derivations $\partial_z$ and $\partial_{\theta^i}$ is called an $N_W = N$
	\emph{formal distribution Lie superalgebra}.
\end{nolabel}
\begin{prop}
	The $\Lambda$-bracket defined in (\ref{eq:lie.2}) satisfies the following
	properties:
	\begin{enumerate}
		\item Sesquilinearity for a pair $(a(Z), b(W))$:
			\begin{xalignat}{2}
				[\partial_z a_\Lambda b] &= -\lambda [a_\Lambda b]
				& [a_\Lambda \partial_w b] &= (\partial_w +
				\lambda) [a_\Lambda b]   \label{eq:lieprop.a} \\
				[\partial_{\theta^i} a_\Lambda b] &= -(-1)^N \chi^i
				[a_\Lambda b] & [a_\Lambda \partial_{\zeta^i} b] &=
				(-1)^{a+N} (\partial_{\zeta^i} + \chi^i) [a_\Lambda b]
				\label{eq:lieprop.1}
			\end{xalignat}
		\item Skew-symmetry for a local pair $(a(Z),b(W))$:
			\begin{equation*}
				[b_\Lambda a] = - (-1)^{ab+N} [a_{-\Lambda -
				\partial_W} b].
			\end{equation*}
		\item Jacobi identity for a triple $(a(Z), b(X), c(W))$:
			\begin{equation*}
				[a_\Lambda[b_\Gamma c]] = (-1)^{aN +N} [ [a_\Lambda
				b]_{\Lambda
				+ \Gamma} c] + (-1)^{(a+N)(b+N)} [b_\Gamma
				[a_\Lambda c]],
			\end{equation*}
	\end{enumerate}
	where $\Gamma = (\gamma, \eta^1, \dots, \eta^N)$ and the superalgebra
	$\mathbb{C}[\Lambda, \Gamma]$ is commutative.
	\label{prop:lie}
\end{prop}
\begin{proof}
	In order to prove the first equation in (\ref{eq:lieprop.1}) we use
	Proposition \ref{prop:fourier} (1):
	\begin{multline*}
			{[}\partial_{\theta^i} a_\Lambda b] = \cF^\Lambda_{Z,W}
			[\partial_{\theta^i} a(Z), b(W)] 
			= \cF^\Lambda_{Z,W} \partial_{\theta^i} [a(Z),b(W)] \\
			= - (-1)^N \chi^i \cF^\Lambda_{Z,W} [a(Z),b(W)] 
			= - (-1)^N \chi^i [a_\Lambda b].
	\end{multline*}
	For the second equation we have by Proposition \ref{prop:fourier} (1):
	\begin{multline*}
			{[}a_\Lambda \partial_{\zeta^i} b] = \cF^\Lambda_{Z,W}
			[a(Z), \partial_{\zeta^i} b(W)] 
			= \cF^\Lambda_{Z,W} (-1)^a \partial_{\zeta^i} [a(Z),b(W)]
			\\
			= (-1)^a \left( [\cF^\Lambda_{Z,W} , \partial_{\zeta^i}] +
			(-1)^N \partial_{\zeta^i} \cF^\Lambda_{Z,W}\right)
			[a(Z),b(W)] \\
			= (-1)^{a+N} (\chi^i + \partial_{\zeta^i})
			\cF^{\Lambda}_{Z,W} [a(Z),b(W)]
			= (-1)^{a+N} \left(  \chi^i +  \partial_{\zeta^i}
			\right) [a_\Lambda b].
	\end{multline*}
	Skew-symmetry follows from the skew-symmetry property of the Fourier
	transform (\ref{eq:fourierprop.2}) as follows:
	\begin{multline*}
			{[}b_\Lambda a] = \cF^\Lambda_{Z,W} [b(Z),a(W)] 
			= -(-1)^{ab} \cF^\Lambda_{Z,W} [a(W),b(Z)] \\
			= - (-1)^{ab+N} \cF^{-\Lambda - \partial_W}_{Z,W}
			[a(Z),b(W)] 
			= - (-1)^{ab+N} [a_{-\Lambda - \partial_W} b].
	\end{multline*}
	Finally, to prove the Jacobi identity we use Proposition \ref{prop:fourier}
	(3)
	\begin{equation*}
		\begin{aligned}
			{[}a_\Lambda[b_\Gamma c] ] &= \cF^\Lambda_{Z,W} [a(Z),
			\cF^\Gamma_{X,W} [b(X),c(W)]] \\
			&= (-1)^{aN} \cF^\Lambda_{Z,W} \cF^\Gamma_{X,W}
			[a(Z),[b(X),c(W)]]  \\
			&= (-1)^{a N} \cF^\Lambda_{Z,W} \cF^\Gamma_{X,W} [ [a(Z),
			b(X)], c(W)] + \\ & \quad +  (-1)^{ab+aN} \cF^\Lambda_{Z,W}
			\cF^\Gamma_{X,W} [b(X), [a(Z), c(W)]] \\
			&= (-1)^{a N + N} \cF^{\Lambda+ \Gamma}_{X,W}
			\cF^\Lambda_{Z,X} [ [a(Z), b(X)], c(W)] + \\ &\quad
			(-1)^{ab+aN+bN+N}
			\cF^\Gamma_{X,W} [b(X), \cF^\Lambda_{Z,W} [a(Z), c(W)]] \\
			&= (-1)^{a N+ N} [[a_\Lambda b]_{\Gamma+\Lambda} c] +
			(-1)^{(a+N)(b+N)} [b_\Gamma [a_\Lambda c]].
		\end{aligned}
	\end{equation*}
\end{proof}
\begin{defn}
	Let $\mathbb{C}[T,S]:=\mathbb{C}[T,S^1,\dots, S^N]$ be the 
	commutative superalgebra freely 
	generated by an even element $T$ and $N$ odd elements $S^i$. A \emph{
	$N_W = N$ SUSY Lie conformal algebra} is a $\mathbb{Z}/2 \mathbb{Z}$-graded
	$\mathbb{C}[T,S]$-module $\cR$ 
	with a $\mathbb{C}$-bilinear operation $[\,_\Lambda \,] : \cR \otimes_\mathbb{C} \cR
	\rightarrow 
	\mathbb{C}[\Lambda] \otimes_\mathbb{C} \cR$ of parity $N \mod 2$ satisfying the
	following three axioms:
	\begin{enumerate}
		\item Sesquilinearity:
			\begin{equation*}
				\begin{aligned}
					{[T} a_\Lambda b] &= -\lambda [a_\Lambda b] &
				[a_\Lambda T b] &= (T + \lambda) [a_\Lambda b] \\
				[S^i a_\Lambda b] &= -(-1)^N \chi^i [a_\Lambda b] &
				[a_\Lambda S^i b] &= (-1)^{a+N} \left(S^i +
				\chi^i \right) [a_\Lambda b] 
				 	\end{aligned}
				\end{equation*}
		\item Skew-symmetry:
			\begin{equation*}
				[b_\Lambda a] = - (-1)^{ab + N} [b_{-\Lambda -
				\nabla} a],
			\end{equation*}
			where $\nabla = (T, S^1, \dots, S^N)$, the $\Lambda$-bracket in the RHS
			means compute first the $\Gamma$
			bracket and then let $\Gamma = - \Lambda - \nabla$.
		\item Jacobi identity:
			\begin{equation}
				[a_\Lambda [b_\Gamma c]] = (-1)^{aN + N} [
				[a_\Lambda b]_{\Gamma + \Lambda} c] +
				(-1)^{(a+N)(b+N)} [b_\Gamma [a_\Lambda c]].
				\label{eq:conformal_def.3}
			\end{equation}
	\end{enumerate}
	We will drop the adjective SUSY when no confusion may arise.
	\label{defn:conformal_def}
\end{defn}
\begin{rem}
	Even though in this case the situation is simple, it is instructive to
	realize the $\Lambda$ bracket as a morphism of
	$\mathbb{C}[\Lambda]$-modules. Consider the co-commutative Hopf superalgebra $\cH =
	\mathbb{C}[\Lambda]$ with commultiplication $\Delta\lambda = \lambda
	\otimes 1 + 1 \otimes \lambda$, $\Delta\chi^i= \chi^i \otimes 1 + 1 \otimes
	\chi^i$. Note 
	that $\mathbb{C}[\nabla] \simeq \cH$. Consider $\cH$ as a $\cH$-module
	with the adjoint action (which is trivial in this case, given that $\cH$
	is super-commutative). Then we may think of $\cH \otimes \cR$ as an $\cH$
	module, the action is given by $h \mapsto \Delta h$. Similarly $\cR \otimes \cR$ is an $\cH$-module. The $\Lambda$-bracket is
	then a $\cH$-module homomorphism of degree $(-1)^N$. Namely, let $\phi$
	denote the morphism $\cR \otimes \cR \rightarrow \cH \otimes \cR$ which
	is given by the $\Lambda$-bracket. Then for every $h \in \cH$ we have
	\begin{equation*}
		\phi h - (-1)^{hN} h \phi = 0,
	\end{equation*}
	as elements in $\Hom (\cR \otimes \cR, \cH \otimes \cR)$.
	Similarly, the Jacobi identity is an identity in 
	\begin{equation*}
		\Hom (\cR \otimes \cR \otimes \cR, \cH \otimes \cH \otimes \cR).
	\end{equation*}\label{rem:explain.morph}
	We will expand on this in Remark \ref{rem:remarketacommute}.
\end{rem}
\begin{rem}
	According to Proposition \ref{prop:lie}, given any $N_W = N$ SUSY formal distribution Lie
	superalgebra $(\fg, \cR)$, the space $\cR$ is a SUSY Lie conformal algebra
	where $T = \partial_w$ and $S^i = \partial_{\zeta^i}$, and the
	$\Lambda$-bracket is defined by (\ref{eq:lie.2}).
	 \label{rem:conformal_from_lie}
\end{rem}
\begin{defn}
	A Lie superalgebra of degree $p \in \mathbb{Z}/2\mathbb{Z}$ is a vector superspace
	$\fh$ with a
	bilinear operation $\{\, , \, \} : \fh \otimes \fh \rightarrow \fh$ of
	parity $p$
	satisfying:
	\begin{enumerate}
		\item Skew-symmetry: $\{ a, b\} = -(-1)^{ab + p} \{b, a\}$.
		\item Jacobi identity: $\bigl\{ a, \{b,c\}\bigr\} = (-1)^{ap+p}
			\bigl\{ \{a,b\}, c\bigr\} + (-1)^{(a+p)(b+p)} \bigl\{ b, \{a,
			c\} \bigr\}$.
	\end{enumerate}
	\label{defn:pseudo.lie.1aa}
\end{defn}
\begin{lem}
	Let $\fh$ be a Lie superalgebra of degree $p \in \mathbb{Z}/2\mathbb{Z}$. Define $\fg$ as a vector
	superspace to be
	$\fh$ if $p=0 \mod 2$ or $\fh$ with the reversed parity if $p =1 \mod 2$. Define the
	bilinear operation $[\, , \,]: \fg \otimes \fg \rightarrow \fg$ by:
	\begin{equation*}
		[a,b] = (-1)^{ap+p} \{a, b\},
	\end{equation*} where the right hand side is computed in $\fh$ and then we
	reverse the parity if $p = \bar{1}$.
	Then $(\fg, [\,,\,])$ is a Lie superalgebra which we will denote as
	$\lie (\fh)$.
	\label{lem:lie.from.pseudoaa}
\end{lem}
\begin{proof}
	We have:
	\begin{equation*}
		[b,a] = (-1)^{bp+p}\{b,a\} =- (-1)^{bp+ab}\{a,b\} = - 
		(-1)^{(a+p)(b+p)}[a,b],
	\end{equation*}
	which is skew-symmetry for the Lie algebra provided the parity in
	$\fg$ is shifted by $p$. To check Jacobi identity we have:
	\begin{multline*}
		[a,[b,c]]  = (-1)^{pb+ap}\{a,\{b,c\}\} = (-1)^{pb+p}\{\{a, b\},c\} +
		(-1)^{ab + p}\{b, \{a,c\}\} = \\ = (-1)^{pb+p+(a+b+p)p+ap}[
		[a,b],c] + (-1)^{ab+p+ap + bp} [b,[a,c]] = \\ =[ [a,b],c] +
		(-1)^{(a+p)(b+p)} [b, [a,c]].
	\end{multline*}
\end{proof}
\begin{lem}
	Let $\cR$ be a $N_W = N$ SUSY Lie conformal algebra. Then $ \cR / \nabla \cR$
	is naturally a Lie superalgebra of degree $N \mod 2$ with bracket
	\begin{equation*}
		\{a + \nabla \cR, b + \nabla \cR\} = [a_\Lambda b]_{\Lambda = 0} + \nabla \cR. 
	\end{equation*}
	\label{lem:Lie_R_1}
\end{lem}
\begin{proof}
	The fact that the bilinear map $\{\, ,\,\}$ is well defined follows from
	sesquilinearity. Skew-symmetry and the Jacobi identity follow from the corresponding
	axioms for the SUSY Lie conformal algebra $\cR$.
\end{proof}
\begin{lem}
	Let $\cR$ be an $N_W=N$ SUSY Lie conformal algebra.  Then $\tilde{\cR}:=\cR \otimes \mathbb{C}[W,
	W^{-1}]$ is an $N_W=N$ SUSY Lie
	conformal algebra with
	$\Lambda$-bracket:
	\begin{equation}
		[a \otimes f_\Lambda b \otimes g] = (-1)^{fb} [a_{\Lambda + \partial_W}b] \otimes
		f(W)g(W')|_{W' = W},
		\label{eq:lambda_bracket_affine}
	\end{equation}
	and with $\tilde{T} = T \otimes \id + \id \otimes \partial_w$ and $\tilde{S^i} = S^i
	\otimes \id + \id \otimes \partial_{\zeta^i}$. 
	\label{lem:affinization_1}
\end{lem}
\begin{proof}
	We prove here skew-symmetry, the other axioms are checked in a similar way:
	\begin{equation*}
		\begin{aligned}
			{[a \otimes f}_\Lambda b \otimes g] &= (-1)^{fb}[a_{\Lambda + \partial_W}b]
			\otimes f(W)g(W')|_{W = W'} \\ 
			&= - (-1)^{ab+N+fb} [b_{-\Lambda - \partial_W - \nabla} a] \otimes
			f(W)g(W')|_{W = W'} \\
			&= - (-1)^{(a+f)(b+g) + N + ga} [b_{-\Lambda - \partial_W - \nabla -
			\partial_{W'}
			+ \partial_{W'}} a] \otimes g(W') f(W)|_{W = W'} \\
			&= - (-1)^{(a+f)(b + g) + N} [b \otimes g_{-\Lambda - \tilde\nabla} a
			\otimes f]
		\end{aligned}
		\label{}
	\end{equation*}
\end{proof}
\begin{nolabel}
	For any $N_W = N$ SUSY Lie conformal algebra $\cR$, we put $L(\cR) = \tilde\cR /
	\tilde\nabla \tilde\cR$ and $\lie(\cR) := \lie (L(\cR))$ (see Lemmas
	\ref{lem:lie.from.pseudoaa} and \ref{lem:Lie_R_1}). For each $a \in
	\mathbb{\cR}$, let $a_{<n|I>} \in L(\cR)$
	be the image of $a \otimes W^{n|I}$. Similarly define $a_{(n|I)} \in \lie(\cR)$ as the
	image of the following element of $L(\cR)$
	\begin{equation*}
		(-1)^{a I} \sigma(I) a_{<n|I>},
		\label{}
	\end{equation*}
	and define, for each $a \in \cR$, the following $\lie(\cR)$-valued formal distribution
	\begin{equation}
		a(Z) = \sum_{j \in \mathbb{Z}, J} Z^{-1-j|N \setminus J} a_{(j|J)} \; \in
		\lie(\cR)[ [Z, Z^{-1}] ].
		\label{eq:distributions_in_lie_R}
	\end{equation}
	Using (\ref{eq:lambda_bracket_affine}) with $f = W^{n|I}$ and $g = W^{k|K}$ and putting
	$\Lambda = 0$ we compute explicitly the Lie bracket (of parity $N \mod 2$) in
	$L(\cR)$:
	\begin{multline}
		\{a_{<n|I>}, b_{<k|K>}\} = \sum_{j \geq 0, J} (-1)^{aJ + b(I- J)}
		\binom{n}{j} \times \\ \times
		 \sigma(J, I \setminus J) \sigma(I\setminus J, K)
		(a_{(j|J)}b)_{<n-j+k|K \cup (I\setminus J)>}.
		\label{eq:commutator_affine.1aa}
	\end{multline}
	It is straightforward to check using Lemma \ref{lem:lie.from.pseudoaa}, that the Lie
	bracket in $\lie(\cR)$ is given by:
	\begin{multline}
		[a_{(n|I)},b_{(k|K)}] = (-1)^{(a+N-I)(N-K)} \sum_{(j|J): j \geq 0}
		(-1)^{(I-J)(N-J)} \binom{n}{j}  \times \\ \times
		\sigma(I)
	 \sigma(J,I\setminus J) \sigma\left(I\setminus J, (N\setminus K) \setminus
		(I\setminus J)\right)
		\bigl(a_{(j|J)}b\bigr)_{(n+k-j|K\cup (I\setminus
		J))}.
		\label{eq:lie_bracket.4aa}
	\end{multline}
	\label{no:formal_from_R_1}
\end{nolabel}
\begin{prop}
	Let $\cR$ be an $N_W=N$ SUSY Lie conformal algebra, $a$, $b$ two vectors in $\cR$, and
	$a(Z)$, $b(W)$ the corresponding $\lie(\cR)$-valued formal distributions defined by
	(\ref{eq:distributions_in_lie_R}). Then 
	\begin{equation}
		{[}a(Z),b(W)] = \sum_{j \geq 0, J}  \left(
		\partial^{(j|J)}_W \delta(Z,W) \right) \left( a_{(j|J)}b \right) (W).
		\label{eq:distributions_are_local}
	\end{equation}
	\label{prop:distributions_are_local}
\end{prop}
\begin{proof}
	First we expand
	\begin{multline}
		\partial^{(j|J)}_W \delta(Z,W) = \sum_{n \in \mathbb{Z}, I} \binom{n}{j}
		(-1)^{I - J} \sigma(J) \sigma(N\setminus I, I \setminus J)
		Z^{-1-n|N\setminus I} W^{n-j|I\setminus J}.
		\label{eq:delta_expanded}
	\end{multline}
	Now using (\ref{eq:lie_bracket.4aa}) we have:
	\begin{multline}
		{[}a(Z), b(W)] = \sum_{\stackrel{n \in \mathbb{Z}, I}{k \in \mathbb{Z}, K}}
		\binom{n}{j} (-1)^{(I-J)(N-J)} \sigma(I)
		\sigma(J,
		I\setminus J)\times \\ \times
		 \sigma(I \setminus J, (N\setminus K)\setminus(I\setminus J))
		Z^{-1-n|N\setminus I}W^{-1-k|N\setminus K} \left( a_{(j|J)}b \right)_{(n+k-j|K
		\cup (I\setminus J)}.
		\label{eq:horrible.1aa}
	\end{multline}
	On the other hand we have
	\begin{equation}
		W^{-1-k|N\setminus K} = \sigma(I\setminus J, (N\setminus K)\setminus
		(I\setminus J)) W^{n-j|I\setminus J} W^{-1-k-n+j|(N\setminus K)\setminus
		(I\setminus J)},
		\label{eq:horrible.2aa}
	\end{equation}
	and, due to (\ref{eq:not.2b}), 
	\begin{equation}
		\sigma(I) \sigma(J, I\setminus J) = (-1)^{(I-J)(N-I)}
		\sigma(N\setminus I, I \setminus J) \sigma(J).
		\label{eq:horrible.3aa}
	\end{equation}
	Now substituting (\ref{eq:horrible.2aa}) in (\ref{eq:horrible.1aa}) and using
	(\ref{eq:horrible.3aa}) we obtain (\ref{eq:distributions_are_local}). 
\end{proof}
\begin{prop}
	Let $\cR$ be an $N_W=N$ SUSY Lie conformal algebra, then the pair $(\lie(\cR), \cR)$ is
	an $N_W=N$ SUSY formal distribution Lie superalgebra.
	\label{prop:formal_distrib_from_conf}
\end{prop}
\begin{proof}
	The fact that the family of distributions (\ref{eq:distributions_in_lie_R}) is closed
	under $(j|J)$-th products and that they are pairwise local follows from Proposition
	\ref{prop:distributions_are_local} since $a_{(j|J)}b = 0$ for $j \gg 0$ in $\cR$. The
	fact that this family is closed under the derivations $\partial_z, \partial_{\theta^i}$
	follows from the following identities which are straightforward to check
	\begin{equation}
		\begin{aligned}
			(Ta)_{(j|J)} &= - j a_{(j-1|J)}, \\
			(S^i a)_{(j|J)} &= \sigma(e_i, N\setminus J) a_{(j| J\setminus e_i)}.
		\end{aligned}
		\label{eq:horrible_alpedo}
	\end{equation}
\end{proof}
\begin{nolabel}
	Note from (\ref{eq:lambda_bracket_affine}) that $(-\partial_w,
	-\partial_{\zeta^i})$  are derivations of the $(0|0)$-th product of
	$\tilde{\cR}$. Since these operators supercommute with $(\tilde{T},
	\tilde{S}^i)$, they induce derivations $(T, S^i)$
	of the Lie superalgebra $\lie(\cR)$, given by the formulas:
	\begin{equation}
		\begin{aligned}
			T(a_{(j|J)}) &= - j a_{(j-1|J)}, \\
			S^i (a_{(j|J)}) &= \begin{cases} \sigma(N \setminus J, e_i)
				a_{(j|J \setminus e_i)} & \text{ if } i \in J, \\
				0 & \text{ if } i \notin J. \end{cases}
		\end{aligned}
		\label{eq:derivations_lie_R}
	\end{equation}
	Note that $\lie(\cR)$ contains a subalgebra $\lie(\cR)_\leq$ spanned by
	vectors $a_{(j|J)}$ with $j \geq 0$. This subalgebra, called the
	\emph{annihilation subalgebra}, is stable under the action of $\nabla = (T,
	S^i)$. 

	Moreover, it is straightforward to check, using
	(\ref{eq:derivations_lie_R}), that the formal distributions
	(\ref{eq:distributions_in_lie_R}) satisfy:
	\begin{equation}
		Ta(Z) = \partial_z a(Z), \qquad S^i a(Z) = \partial_{\theta^i} a(Z)
		\label{eq:lie_R_is_regular}
	\end{equation}
	namely, the $N_W=N$ formal distribution Lie superalgebra $(\lie(\cR),
	\cR)$ is \emph{regular}.
\label{no:regular_lie_R}
\end{nolabel}
\begin{nolabel}
	Recall that we have defined $(j|J)$-th products of formal distributions for $j
	\geq 0$ in \ref{no:lie}. In order to define these products for $j < 0$ we
	let for a formal distribution $a(Z)= \sum
	Z^{j|J} a_{j|J}$:
	\begin{equation*}
			a_+(Z) = \sum_{(j|J): j \geq 0} Z^{j|J} a_{j|J}, \qquad
			a_-(Z) = \sum_{(j|J): j < 0} Z^{j|J} a_{j|J}
	\end{equation*}
	It follows easily from the definitions that
	\begin{equation}
		\begin{aligned}
			a_+(W) = \res_Z i_{z,w} (Z-W)^{-1|N} a(Z), \qquad
			a_-(W) = -\res_Z i_{w,z} (Z-W)^{-1|N} a(Z).
		\end{aligned}
		\label{eq:normal.2}
	\end{equation}
	Indeed, we have
	\begin{equation*}
		i_{z,w} (Z-W)^{-1|N} = \sum_{(m|J): m \geq 0} (-1)^J
		\sigma(J) W^{m|J} Z^{-1-m|N\setminus J}. 
	\end{equation*}
	Hence:
	\begin{multline*}
		\res_Z i_{z,w}(Z-W)^{-1|N} a(Z) = \res_Z
		\sum_{\stackrel{(m|J): m \geq
		0}{(n|I): n \in \mathbb{Z}}} (-1)^J
		\sigma(J) W^{m|J} Z^{-1-m|N\setminus J} Z^{n|I}
		a_{n|I}  \\ 
		= \res_Z \sum_{(m|J): m \geq 0} (-1)^{J} \sigma(J)
		\sigma(N\setminus J, J) W^{m | J} Z^{-1|N} a_{m | J} \\ = 
		\sum_{(m|J): m \geq 0} W^{m|J} a_{m|J} = a_+(W).
	\end{multline*}
	The second equation in (\ref{eq:normal.2}) follows similarly, or by noting
	that it is a consequence of the first equation in (\ref{eq:normal.2}), the
	definition of the $\delta$ function (\ref{eq:delta.1}) and property (5) in
	\ref{no:delta}. Differentiating (\ref{eq:normal.2}) we find:
	\begin{equation}
		\begin{aligned}
			(-1)^{JN} \partial_W^{(j|J)}
			a(W)_+ &= \sigma(J) \res_Z i_{z,w}
			(Z-W)^{-1-j|N\setminus J} a(Z), \\
			(-1)^{JN} \partial_W^{(j|J)}
			a(W)_- &= -\sigma(J) \res_Z i_{w,z}
			(Z-W)^{-1-j|N\setminus J} a(Z).\\
		\end{aligned}
		\label{eq:normal.5}
	\end{equation}
	These equations (\ref{eq:normal.5}) are called the super \emph{Cauchy
	formulae}.
	\label{no:normal}
\end{nolabel}
\begin{defn}
	Let $V$ be a vector superspace. An $\End(V)$-valued formal distribution
	$a(Z)$ is called a \emph{field} if for every vector $v \in V$ we have $a(Z)
	v \in V( (Z))$, i.e. there are finitely many negative powers of $z$ in
	$a(Z)v$. For two such fields we define their \emph{normally ordered
	product} to be 
	\begin{equation}
		:a(Z)b(Z): \; := a_+(Z) b(Z) + (-1)^{ab} b(Z)a_-(Z)
		\label{eq:normal_defin}
	\end{equation}
	\label{defn:field_and_normal}
\end{defn}
\begin{nolabel}
	The normally ordered product of fields is again a well defined field. Indeed, when
	applied to any vector $v \in V$ the first summand in
	(\ref{eq:normal_defin}) clearly has finitely many negative powers of $z$
	since $b(Z) v \in V( (Z))$ and $a_+(Z)$ has only non-negative powers of
	$z$. For the second summand we see that $a_-(Z) v \in V[Z,Z^{-1}]$, namely it is a
	Laurent polynomial with values in $V$, therefore $b(Z) a_-(Z) v \in V(
	(Z))$ as we wanted.
	\label{no:normal_is_field}
\end{nolabel}
\begin{lem}
	\begin{multline}
		:a(W)b(W): = \res_Z \left( i_{z,w} (Z-W)^{-1|N} a(Z) b(W) -
		\right. \\ \left. - (-1)^{ab} i_{w,z} (Z-W)^{-1|N} b(W)a(Z)
		\right).
		\label{eq:normal_lem.1}
	\end{multline}
	\label{lem:normal_lem}
\end{lem}
\begin{proof}
	This is immediate by (\ref{eq:normal.2}). 
\end{proof}
\begin{nolabel}
	Given the last lemma and the Cauchy formulae (\ref{eq:normal.5}) it
	is natural to define
	\begin{equation}
		a(W)_{(-1-j|N\setminus J)}b(W) = \sigma(J)
		(-1)^{JN}
		:\left(\partial_W^{(j|J)}a(W)\right) b(W):.
		\label{eq:negative.1}
	\end{equation}
	Differentiating (\ref{eq:normal_lem.1}) we find:
	\begin{multline*}
		a(W)_{(-1-j|N\setminus J)} b(W) = \res_{Z}\left( \left( i_{z,w}
		(Z-W)^{-1-j|N\setminus J}\right) a(Z)b(W) - \right. \\ \left. - (-1)^{ab}
		\left(i_{w,z} (Z-W)^{-1-j|N\setminus J} \right) b(W) a(Z) \right).
	\end{multline*}
	Similarly, from the definition of the $j|J$-th products for $j \geq 0$ in
	\ref{no:lie} we have:
	\begin{multline}
		\res_Z\left( \left( i_{z,w} (Z-W)^{j|J}\right) a(Z) b(W) - (-1)^{ab} \left(i_{w,z}
		(Z-W)^{j|J}\right) b(W)a(Z) \right) = \\ = \res_Z (Z-W)^{j|J} \left(a(Z)
		b(W) - (-1)^{ab} b(W) a(Z) \right) = \\ = \res_Z (Z-W)^{j|J} [a(Z),
		b(W)] = a(W)_{(j|J)} b(W).
	\end{multline}
	Therefore we have proved that for every $j \in \mathbb{Z}$ and every tuple
	$J$ we have:
	\begin{multline}
		a(W)_{(j|J)} b(W) = \res_Z\left(  \left(i_{z,w} (Z-W)^{j|J}\right) a(Z) b(W) -
		\right. \\ \left. -
		(-1)^{ab} \left(i_{w,z}
		(Z-W)^{j|J} \right) b(W)a(Z) \right).
		\label{eq:negative.4}
	\end{multline}
	\label{no:negative}
\end{nolabel}
\begin{prop}
	The following identities analogous to sesquilinearity for all pairs $j|J$
	are true:
	\begin{equation}
		\begin{aligned}
			\left(\partial_w a(W)\right)_{(j|J)} b(W) &= - j
			a(W)_{(-j-1|J)}b(W) \\
			\partial_w \left( a(W)_{(j|J)} b(W) \right) &= \left(
			\partial_w a(W)
			\right)_{(j|J)} b(W) + a(W)_{(j|J)} \partial_w b(W) \\
			\left(\partial_{\zeta^i} a(W) \right)_{(j|J)} b(W) &=
			 \sigma (J \setminus e_i,e_i) a(W)_{(j|J \setminus
			e_i)}b(W) \\
			\partial_{\zeta^i} \left( a(W)_{(j|J)}b(W) \right) &=
			(-1)^{N-J} \left( \left(
			\partial_{\zeta^i} a(W) \right)_{(j|J)} b(W) \right. + \\
			& \left. 
			\quad  + (-1)^a
			a(W)_{(j|J)} \left( \partial_{\zeta^i} b(W) \right) \right),
		\end{aligned}
		\label{eq:sesquilinearity.1}
	\end{equation}
	where $e_i$ is the tuple consisting of only one element $\{i\}$ and we
	recall that we are defining $\sigma(e_i, J\setminus e_i)$ to be zero if
	$i \not\in J$. 
	\label{prop:sesquilinearity}
\end{prop}
\begin{proof}
	The first two equations are standard and their proof is similar to the last
	two. We will prove the last two equations by using (\ref{eq:negative.4}).
	If $i \not\in J$ the result is obvious. 
	\begin{multline}
			\res_Z i_{z,w} (Z-W)^{j|J} \partial_{\theta^i}a(Z)b(W) =
			-(-1)^{J} \res_Z \left( \partial_{\theta^i} i_{z,w}
			(Z-W)^{j|J} \right) a(Z) b(W) \\ = -(-1)^J \sigma(e_i, J
			\setminus e_i)
			\res_Z i_{z,w}
			(Z-W)^{j|J\setminus e_i} a(Z) b(W).
			\label{eq:sesqui_proof.1}
	\end{multline}
	Similarly we have:
	\begin{multline}
		- (-1)^{(a+1)b} \res_Z i_{w,z} (Z-W)^{j|J} b(W)\partial_{\theta^i} a(Z)
		= \\ = (-1)^{ab + J} \res_Z \left( \partial_{\theta^i} i_{w,z}
		(Z-W)^{j|J} \right) b(W)a(Z) = \\ = (-1)^{ab + J} \sigma(e_i, J
		\setminus e_i) \res_Z i_{w,z} (Z-W)^{j|J\setminus e_i} b(W)a(Z).
		\label{eq:sesqui_proof.2}
	\end{multline}
	Adding (\ref{eq:sesqui_proof.1}) and (\ref{eq:sesqui_proof.2}) and using
	(\ref{eq:not.2b}), we obtain the third equation in
	(\ref{eq:sesquilinearity.1}).

	Finally, to prove the last relation in (\ref{eq:sesquilinearity.1}) we expand:
	\begin{multline}
		\partial_{\zeta^i} \left( a(W)_{(j|J)}b(W) \right) =
		\partial_{\zeta^i} \res_Z \left( i_{z,w} (Z-W)^{j|J} a(Z) b(W) -
		\right. \\ - \left. (-1)^{ab} i_{w,z} (Z-W)^{j|J} b(W)a(Z) \right)
		= \\ (-1)^N  \res_Z \left(- \sigma(e_i,J\setminus e_i) i_{z,w}
		(Z-W)^{j|J \setminus e_i} a(Z)b(W) + \right. \\ \left. + (-1)^{J+a}
		i_{z,w} (Z-W)^{j|J} a(Z)\partial_{\zeta^i}b(W) - \right. \\ \left.
		- (-1)^{ab} \sigma(e_i, J\setminus e_i) i_{w,z} (Z-W)^{j|J
		\setminus e_i} b(W) a(Z) - \right. \\ \left. - (-1)^{ab + J}
		i_{w,z} (Z-W)^{j|J} \partial_{\zeta^i}b(W) a(Z)
		\right) = \\ = - (-1)^N \sigma(e_i, J\setminus e_i)
		a(W)_{(j|J\setminus e_i)} b(W) + (-1)^{N + J+a}
		a(W)_{(j|J)}\partial_{\zeta^i}b(W) = \\ = (-1)^{N-J} \left( \left(
		\partial_{\zeta^i} a(W)
		\right)_{(j|J)} b(W) + (-1)^{a} a(W)_{(j|J)} \partial_{\zeta^i}b(W)
		\right).
		\label{eq:sesqui_proof.4}
	\end{multline}
\end{proof}
\begin{prop}
	The following identity holds for any $(j|J)$ and any three fields $a=a(W)$,
	$b =b(W)$, $c = c(W)$:
	\begin{multline}
		[a_\Lambda (b_{(j|J)}c)] =  \\ =\sum_{(k|K): k \geq 0}
		(-1)^{(a + K + N)(J + N)}  \sigma(J,K)
		\Lambda^{(k|K)} [a_\Lambda b]_{(j+k|J \cup K)} c + \\ +
		(-1)^{(a+N)(b+N-J)}b_{(j|J)}[a_\Lambda c].
		\label{eq:wick.1}
	\end{multline}
	\label{prop:wick.1}
\end{prop}
\begin{proof}
	The left hand side is
	\begin{multline}
		\res_Z \exp\left( (Z-W)\Lambda \right) [a(Z), (b(W)_{(j|J)}c(W))]
		=\\ =
		\res_Z \exp\left( (Z-W) \Lambda \right) \left( [a(Z), \res_X i_{x,w}
		(X-W)^{j|J} b(X)c(W)] - \right. \\ \left. - (-1)^{bc} 
		[a(Z), \res_X i_{w,x} (X-W)^{j|J} c(W) b(X)] \right) = \\
		(-1)^{a(N-J)} \res_Z \res_X \exp\left( (Z-W)\Lambda \right) i_{x,w}
		(X-W)^{j|J} [a(Z), b(X)c(W)] - \\ - (-1)^{bc + a(N-J)} \res_Z
		\res_X \exp \left( (Z-W) \Lambda\right)i_{w,x} (X-W)^{j|J} [a(Z),c(W)b(X)]
		\label{eq:wick.prop.proof.1}
	\end{multline}
	Using the identity $[a,bc] = [a,b]c + (-1)^{ab} b[a,c]$ we can write the
	first term of the RHS of the last equality as:
	\begin{multline}
		(-1)^{a(N-J)} \res_Z \res_X \exp\left( (Z-X+X-W)\Lambda \right)
		\times \\ \times i_{x,w}
		(X-W)^{j|J} [a(Z), b(X)] c(W) + (-1)^{a(N-J+b)} \res_Z \res_X
		\exp\left( (Z-W)\Lambda \right) \times \\ \times i_{x,w}
		(X-W)^{j|J} b(X)[a(Z),c(W)]
		= \\ = (-1)^{a(N-J)+ N + JN} \res_X \exp\left( (X-W)\Lambda \right)
		i_{x,w} (X-W)^{j|J}
		[a_\Lambda b](X) c(W) + \\ + (-1)^{a(N-J + b) + N + JN + bN}\res_X
		i_{x,w} (X-W)^{j|J} b(X)[a_\Lambda c](W) = \\ = (-1)^{(a+N)(N-J)}
		\res_X \sum_{(k|K):k \geq 0} \frac{(-1)^{\frac{K(K+1)}{2}}}{k!}
		\Lambda^{k|K}  \times \\ \times i_{x,w} (X-W)^{k|K} (X-W)^{j|J}
		[a_\Lambda b](X) c(W)
		+  \\ + (-1)^{(a+N)(N-J +b)} \res_X i_{x,w} (X-W)^{j|J}
		b(X)[a_\Lambda c](W) = \\ 
		= (-1)^{(a+N)(N-J)} \res_X \sum_{(k|K): k \geq 0}
		\frac{(-1)^{\frac{K(K+1)}{2}}}{k!} \sigma(K,J) \times \\ \times
		\Lambda^{k|K} i_{x,w} (X-W)^{k+j|K \cup J} [a_\Lambda b](X) c(W) + \\ 
		+ (-1)^{(a+N)(N-J+b)} \res_X i_{x,w} (X-W)^{j|J} b(X)[a_\Lambda
		c](W)
		\label{eq:wick.prop.proof.2}
	\end{multline}
	Similarly the second term in the RHS of the last equality of
	(\ref{eq:wick.prop.proof.1}) can be written
	as:
	\begin{multline}
		- (-1)^{bc + a(N-J)} \res_Z \res_X \exp\left( (Z-W) \Lambda \right)
		\times \\ \times i_{w,x} (X-W)^{j|J} [a(Z),c(W)]b(X)  -
		(-1)^{bc + a(N-J+c)}
		\res_Z \res_X \times \\ \times \exp\left( (Z-W)\Lambda \right) 
		i_{w,x} (X-W)^{j|J}
		c(W)[a(Z),b(X)] = \\ = - (-1)^{bc + a(N-J) + N + JN} \res_X i_{w,x}
		(X-W)^{j|J} [a_\Lambda c](W) b(X) - \\ - (-1)^{bc +
		(a+N)(N-J+c)} \res_X \exp\left( (X-W)\Lambda \right) i_{w,x}
		(X-W)^{j|J} c(W)[a_\Lambda b](X) \\ = - (-1)^{bc+ (a+N)(N-J)}
		\res_X i_{w,x} (X-W)^{j|J} [a_\Lambda c](W) b(X) - \\ - (-1)^{bc +
		(a+N)(N-J+c)} \res_X \sum_{(k|K): k \geq 0}
		\frac{(-1)^{\frac{K(K+1)}{2}}}{k!} \times \\ \times \sigma(K,J)
		\Lambda^{k|K} i_{w,x} (X-W)^{k+j, K \cup J} c(W)[a_\Lambda b](X).
		\label{eq:wick.prop.proof.3}
	\end{multline}
	Now adding (\ref{eq:wick.prop.proof.2}) and (\ref{eq:wick.prop.proof.3}) we
	get (\ref{eq:wick.1}) (recall that the $\Lambda$-bracket has parity $N \,
	\mathrm{mod} 2$).
\end{proof}
\begin{rem}
	If we multiply both sides of (\ref{eq:wick.1}) by
	\begin{equation*}
		\frac{(-1)^{\frac{J(J+1 +2 a)}{2}}}{j!} \Gamma^{j|J},
	\end{equation*}
	and sum over all pairs $(j|J)$ with $j \geq 0$ we obtain the Jacobi
	identity for the $\Lambda$-bracket that we have already proved in Proposition
	\ref{prop:lie}. Therefore, the identities (\ref{eq:wick.1}) for $j \geq 0$
	are equivalent to the Jacobi identity (\ref{eq:conformal_def.3}).

	Next we note that if we replace $b$ by $\partial_w b$ in (\ref{eq:wick.1})
	we obtain the same identity with $j$ replaced by $j-1$ whenever $j \leq
	-1$. 
	Similarly, replacing $b$ by $\partial_{\zeta^i} b$ we obtain the same
	identity with $J$ replaced by $J \setminus e_i$. It follows the identity
	(\ref{eq:wick.1}) is equivalent to the Jacobi identity (\ref{eq:conformal_def.3})
	and (\ref{eq:wick.1}) with
	$(j|J) = (-1|N)$. In
	this case the formula (\ref{eq:wick.1}) looks as follows:
	\begin{equation*}
		[a_\Lambda :b c:] = \sum_{k \geq 0} \frac{\lambda^k}{k!} [a_\Lambda
		b]_{(k-1|N)}c + (-1)^{(a + N)b} :b[a_\Lambda c]:
	\end{equation*}
	Rewriting the sum as the sum of the $k=0$ term and the rest, this becomes:
	\begin{equation}
		[a_\Lambda :bc:] =  :[a_\Lambda b] c: + (-1)^{(a+N)b}
		:b [a_\Lambda c]:  + \int_0^\Lambda [ [a_\Lambda b]_\Gamma c]
		d\Gamma.
		\label{eq:wick.proof.8}
	\end{equation}
	Here the integral $\int_0^\Lambda$ is computed by taking the indefinite integral in the
	even variable $\gamma$ of $\partial_\eta^N$ of the integrand, and then taking the
	difference of the values at the limits.
	This is the super analogue of the \emph{non-commutative Wick formula}
	\cite{kac:vertex}. 
	Thus, the identity (\ref{eq:wick.1}) is equivalent to the
	Jacobi identity plus this non-commutative Wick formula. 
\end{rem}
The following lemma is proved as in the ordinary vertex algebra case
\cite[Lem. 3.2]{kac:vertex}.
\begin{lem}[Dong's Lemma]
	Given three pairwise local formal distributions $a,b,c$, the pair $(a,
	b_{(j|J)}c)$ is 
	local for any $(j|J)$.
	\label{lem:dong}
\end{lem}
\subsection{Identities and existence theorem}
\label{sub:existence}
In this section we define $N_W=N$ SUSY vertex algebras, derive their identities, and prove an existence
theorem as in the non-super case \cite[Thm. 4.5]{kac:vertex}. 
\begin{defn}
	An $N_W=N$ SUSY vertex algebra consists of a vector superspace $V$,
	an even vector $\vac \in V$, $N$ odd operators
	$S^i$ (the odd translation operators), an even operator $T$ (the
	even translation operator), and a parity preserving linear map
	$\ys$ from $V$ to the space of $\End(V)$-valued superfields $a \mapsto \ys(a,Z)$. The following
	axioms must be satisfied:
	\begin{itemize}
		\item Vacuum axioms:
			\begin{equation*}
				\begin{aligned}
					\ys(a,Z) \vac &= a +
					O(Z), \quad
					T\vac = S^i \vac = 0, \qquad i=1,\dots,
					N.\\
				\end{aligned}
			\end{equation*}
		\item Translation invariance
			\begin{equation}
				\begin{aligned}
				{[}S^i, \ys(a,Z)] &=
				\partial_{\theta^i}  \ys
				(A,Z),\quad
				[T, \ys(a, Z)] = \partial_z
				\ys(a, Z).
			\end{aligned}
				\label{eq:w_defin.2}
			\end{equation}
		\item Locality
			\begin{equation*}
				(z-w)^n [\ys(a,Z),
				\ys(b, W)] = 0 \quad \text{for some } n
				\in \mathbb{Z}_+.
				\label{}
			\end{equation*}
			As before, $O(Z)$ is an element of $V[ [Z]]$ which vanishes
			at $Z=0$.
	\end{itemize}
	Morphisms between $N_W=N$ SUSY vertex algebras are linear maps $f :V_1 \rightarrow V_2$
	such that:
	\begin{equation*}
		\begin{aligned}
			f \circ T_1 &= T_2 \circ f,  \quad
			f(\ys_1(a,Z)b) = \ys(f(a),Z)f(b), \quad \forall\, a,b \in V_1  .
		\end{aligned}
		\label{}
	\end{equation*}
	\label{defn:w_defin}
\end{defn}
\begin{nolabel}
	Given a $N_W=n$ SUSY vertex algebra $V$, we can define the $(j|J)$ product of two
	vectors of $V$ as follows. Expand the field $\ys(a,Z)$ for $a \in V$:
	\begin{equation}
		\ys(a,Z) = \sum_{(j|J): j \in \mathbb{Z}}  Z^{-1-j|N\setminus
		J} a_{(j|J)},
		\label{eq:product.1}
	\end{equation}
	and define the $j|J$-product of two vectors in $V$
	as:
	\begin{equation}
		a_{(j|J)} b := a_{(j|J)}\bigl(b\bigr).
		\label{eq:product.3}
	\end{equation}
	This is a $\mathbb{C}$-bilinear product on  $V$ of parity $N-J \mod 2$. 
	We can rewrite the axioms of the vertex algebra in terms of these products.
	For example, 
	the vacuum axioms are equivalent to:
	\begin{xalignat}{2}
		a_{(-1|N)}\vac &= a, & a_{(j|J)}\vac &= 0
		\quad \text{if } j \geq 0 ,\\
		T\vac = S^i \vac &=0. & &
		\label{eq:product.4}
	\end{xalignat}
	and translation invariance is equivalent to:
	\begin{equation}
		\begin{aligned}
			{[}T, a_{(j|J)}] &= -j a_{(j-1|J)} ,\\
			[S^i,a_{(j|J)}] &= \begin{cases} \sigma(N\setminus J, e_i)
			a_{(j|J\setminus e_i)} \quad &\text{if } i \in J ,\\
			0 \qquad &\text{if } i \notin J.\end{cases}
		\end{aligned}
		\label{eq:product.5}
	\end{equation}

	Of course the fact that $\ys(a,Z)$ is a field is equivalent to
	 $a_{(j|J)}b = 0$ for $j \gg 0$, given $a, b \in V$.
	\label{no:product_axioms}
\end{nolabel}
\begin{thm}
	Let $\cU$ be a vector superspace and $V$ a space of pairwise local  $\End(\cU)$-valued
	fields such that $V$ contains the constant field $\Id$, it is invariant
	under the derivations $\partial_z, \partial_{\theta^i}$ and closed under
	all $(j|J)$-th products. Then $V$ is a $N_W=N$ SUSY vertex algebra with
	vacuum vector $\Id$, translation operators  $T a(Z) = \partial_z
	a(Z)$ and $S^i a(Z) =
	\partial_{\theta^i} a(Z)$, and the $(j|J)$ products are given by the RHS of 
	(\ref{eq:negative.4}) multiplied by $\sigma(J)$\footnote{This
	normalization becomes is necessary because of our choice in
	(\ref{eq:negative.1}), see also Theorem \ref{thm:structure.2}}.
	\label{thm:fields.1}
\end{thm}
\begin{proof}
	To check the vacuum axioms we have:
	\begin{equation*}
		\begin{aligned}
			a(Z)_{(j|J)}1 &= \sigma(J) \res_Z (Z-W)^{j|J}
			[a(Z),1] = 0 \quad \text{if } j \geq 0,\\
			a(Z)_{(-1|N)}1 &= :a(Z)1: = a(Z), \quad
			\partial_z 1 = \partial_{\theta^i} 1 = 0.
		\end{aligned}
	\end{equation*}
	To check translation invariance we have:
	\begin{equation*}
		\partial_z(a(Z)_{(j|J)}b(Z)) - a(Z)_{(j|J)}\partial_zb(Z) = 
		(\partial_z a(Z))_{(j|J)}b(Z),
	\end{equation*}
	but this is $-j a(Z)_{(j-1|J)}b(Z)$, according to
	(\ref{eq:sesquilinearity.1}). Therefore we see that the first equation in
	(\ref{eq:product.5}) holds. For the odd translation operators we write (note
	that the parity of $a_{(j|J)}$ is $a+N-J$ since $\ys$ is parity preserving
	and our choice of decomposing the field in (\ref{eq:product.1})):
	\begin{multline*}
			  \sigma(J) \Bigl(\partial_{\theta^i}
			(a(Z)_{(j|J)}b(Z)) - (-1)^{a+N-J}
			a(Z)_{(j|J)}\partial_{\theta^i}b(Z)\Bigr) = \\=
			(-1)^{N-J} \sigma(J)
			(\partial_{\theta^i}a(Z))_{(j|J)}b(Z),
	\end{multline*}
	and again by (\ref{eq:sesquilinearity.1}) we see that this is 
	\begin{multline*}
		- (-1)^N  \sigma(J) \sigma(e_i, J\setminus e_i)
		a(Z)_{(j|J \setminus e_i)}b(Z)  = \sigma(N\setminus J,
		e_i)\sigma(J\setminus e_i)
		a(Z)_{(j|J\setminus e_i)}b(Z),
	\end{multline*}
	proving the second identity in (\ref{eq:product.5}). In order to
	check locality, we expand 
	\begin{equation}
		\begin{aligned}
			\ys(a(W), X)b(W) &= \sum_{(j|J): j \in \mathbb{Z}} \sigma(J)
			X^{-1-j|N\setminus J} a(W)_{(j|J)}b(W) \\
			&= \res_Z \sum_{(j|J): j \in \mathbb{Z}} (-1)^{(N-J)N} \sigma(J)
			X^{-1-j|N \setminus J} \times \\ & \quad \times \Bigl(
			i_{z,w} (Z-W)^{j|J} a(Z)b(W)
			- (-1)^{ab}i_{w,z}(Z-W)^{j|J}b(W)a(Z) \Bigr) \\
			&= \res_Z \sum_{(j|J): j \in \mathbb{Z}} (-1)^{N-J} \sigma(J)
			 \Bigl(
			i_{z,w} (Z-W)^{j|J} X^{-1-j|N \setminus J} \times \\ &
			\quad \times a(Z)b(W)
			- (-1)^{ab}i_{w,z}(Z-W)^{j|J} X^{-1-j|N \setminus J}
			b(W)a(Z) \Bigr). \\
		\end{aligned}
		\label{eq:fields.proof.5}
	\end{equation}
	We note that
	\begin{equation}
		i_{z,w} \sum_{(j|J): j\in \mathbb{Z}}  (-1)^{(N-J)}\sigma(J)
		(Z-W)^{j|J}
		X^{-1-j|N\setminus J} =
		i_{z,w} \delta(Z-W,X).
		\label{eq:fields.proof.6}
	\end{equation}
	Therefore the RHS of (\ref{eq:fields.proof.5}) reads:
	\begin{equation*}
		\res_Z \Bigl( i_{z,w} \delta(Z-W,X)a(Z)b(W) - (-1)^{ab} i_{w,z}
		\delta(Z-W,X) b(W)a(Z) \Bigr).
	\end{equation*}
	With this last equation we can compute then the commutator $[\ys(a(W)),
	\ys(b(W))]c(W)$. Indeed, the product $\ys(a(W),X)\ys(b(W),Y)c(W)$ is
	given by:
	\begin{multline}
		\res_Z \res_U \Bigl( i_{u,w} i_{z,w}
		\delta(U-W,X)\delta(Z-W,Y)a(U)b(Z)c(W) - \\ - (-1)^{bc} i_{u,w} i_{w,z}
		\delta(U-W,X)\delta(Z-W,Y)a(U)c(W)b(Z) -\\- (-1)^{a(b+c)}
		i_{w,u}i_{z,w} \delta(U-W,X)\delta(Z-W,Y) b(Z)c(W)a(U) +
		\\ + (-1)^{a(b+c)+bc} i_{w,u}i_{w,z} \delta(U-W,X)\delta(Z-W,Y)
		c(W)b(Z)a(U) \Bigr),
		\label{}
	\end{multline}
	and we get a similar expression for the product
	$\ys(b(W),Y)\ys(a(W),X)c(W)$. Subtracting we obtain:
	\begin{multline}
		[\ys(a(W),X),\ys(b(W),Y)]c(W) = \\ \res_Z \res_U \Bigl( i_{u,w}i_{z,w}
		\delta(U-W,X) \delta(Z-W,Y) [a(U),b(Z)]c(W) -\\ - (-1)^{(a+b)c}
		i_{w,u}i_{w,z} \delta(U-W,X) \delta(Z-W,Y) c(W)[a(U),b(Z)]
		\Bigr).
		\label{eq:fields.proof.8}
	\end{multline}
	Let $n \in \mathbb{Z}_+$ be such tht $(u-z)^n[a(U), b(Z)] = 0$. Multiplying
	(\ref{eq:fields.proof.8}) by $(x-y)^n$ we obtain that the RHS  
	vanishes. Indeed, using
	\begin{equation*}
		(x-y) = (z-u)-( (z-w) - x) + ( (u-w)- y),
	\end{equation*}
	we see that all terms in the expansion of $(x-y)^n$ vanish when multiplied by
	$\delta$ functions, with the exception of $(z-u)^n$. But this term
	vanishes
	when multiplied by the factors $[a(U),b(Z)]$ in (\ref{eq:fields.proof.8}).
	 Therefore we have proved locality and the theorem.
\end{proof}
\begin{cor} Any identity on elements of an $N_W = N$ SUSY vertex algebra, holds for
	any collection of pairwise local fields. 
\end{cor}
\begin{lem}
	Let $V$ be a vector superspace and let $\vac$ be an even vector of  $V$. Let $a(Z),
	b(Z)$ be two
	$\End(V)$-valued fields such that $a(Z) \vac \in V[ [Z]]$ and $b(Z)
	\vac \in V[ [Z]]$. Then for all $(j|J)$,
	$a(W)_{(j|J)}b(W) \vac \in V[ [W]]$ and the constant term is 
	\begin{equation}
		\sigma(J) a_{(j|J)}b_{(-1|N)}
	\vac.
	\label{eq:constant_term}
	\end{equation}
	\label{lem:constant_term}
\end{lem}
\begin{proof}
	Applying both sides of (\ref{eq:negative.4}) to the vacuum, we see that the second
	term on the RHS of (\ref{eq:negative.4}) vanishes since it contains only positive
	powers of $z$. The first term in the RHS contains only positive powers of $w$ since
	$i_{z,w} (Z-W)^{j|J}$ does and $b(W)\vac \in \mathbb{C}[ [W]]$. Letting $W=0$ we get
	\begin{equation}
		a(W)_{(j|J)}b(W) \vac|_{W = 0} = \res_Z Z^{j|J} a(Z)\left( b_{(-1|N)} \vac \right).
		\label{eq:constant_proof_1a}
	\end{equation}
	It follows from (\ref{eq:distrib.4}) that the RHS of (\ref{eq:constant_proof_1a}) is
	(\ref{eq:constant_term}). 
\end{proof}
The following lemma is straightforward
\begin{lem}
	Let $A$ and $B_1, \dots, B_N$ be linear operators on a vector superspace $\cU$.
	Suppose that $A$ is even and $B_i$ are odd and they pairwise (super) commute, i.e. $A
	B_i = B_i A$, $B_i B_j = - B_j B_i$. 
	Then there
	exists a unique solution $f(Z) \in \cU[ [Z]]$ to the system of differential equations:
	\begin{equation}
		\partial_z f(Z) = A f(Z), \qquad \partial_{\theta^i} f(Z) = B_i f(Z) \; (i = 1,
		\dots, N),
		\label{eq:differential.1}
	\end{equation}
	for any initial condition $f(0)= f_0$.
	\label{lem:differential}
\end{lem}
\begin{proof}
	Using (\ref{eq:differential.1}), the coefficients of $f(Z)$ can be computed by
	induction, given $f_0$.  
\end{proof}
\begin{prop}\renewcommand{\theenumi}{\alph{enumi}}
	Let $V$ be a $N_W=N$ SUSY vertex algebra. Then for every $a,b \in V$:
	\begin{enumerate}
		\item $\ys(a,Z)\vac = \exp( Z\nabla) a.$ 
		\item $\exp(Z \nabla) \ys(a,W) \exp( -Z \nabla ) = i_{w,z} \ys(a,Z
			+W).$
		\item $\ys(a,Z)_{(j|J)}\ys(b,Z) \vac = \sigma(J)
			\ys(a_{(j|J)}b,Z)\vac$,
	\end{enumerate}
	\renewcommand{\theenumi}{\arabic{enumi}}
	where $\nabla = (T,S^1,\dots,S^N)$ and $Z \nabla  = zT+ \sum_i \theta^i S^i$.
	\label{prop:structure_thm_comienzo}
\end{prop}
\begin{proof}
	We note that both sides in (a) and (c) are elements of $V[ [Z]]$ whereas
	both sides of (b) are elements of $\End(V)[ [W, W^{-1}]][ [Z]]$. Note that
	by evaluating at $Z = 0$ we get equalities in all three cases, the only
	non-trivial case is (c), but it follows from Lemma
	\ref{lem:constant_term}. Let us
	denote the right hand side in each case by
	$X(Z)$. It is easy to show that it satisfies the following systems of
	equations respectively:
	\begin{enumerate}
		\item $\partial_z X(Z) = T X(Z)$, and $\partial_{\theta^i} X(Z) =
			S^i X(Z)$.
		\item $\partial_z X(Z) = [T, X(Z)]$ and $\partial_{\theta^i} X(Z) =
			[S^i, X(Z)]$ by the translation axioms.
		\item $\partial_z X(Z) = T X(Z)$ and $\partial_{\theta^i} X(Z) =
			S^i X(Z)$ by the translation axioms (recall that $T \vac = S^i
			\vac = 0$).	
	\end{enumerate}
	In order to apply Lemma \ref{lem:differential}, we have to show that the left hand side of (a), (b) and (c)
	satisfies the same differential equations (1), (2) and (3) respectively;
\begin{enumerate}
	\item It is immediate by the translation invariance and the second of the vacuum
		axioms.
	\item Denoting $Y(Z) =e^{Z \nabla} \ys(a,W) e^{-Z
		\nabla}$, we have:
		\begin{equation*}
				\partial_z Y(Z) = T Y(Z) - Y(Z)T =
				 [T, Y(Z)],
		\end{equation*}
		and similarly:
		\begin{equation*}
			\begin{aligned}
				\partial_{\theta^i} Y(Z) &= S^i Y(Z) + (-1)^{a}
				e^{Z \nabla} \ys(a,W) (-S^i) e^{-Z \nabla} \\
				&= S^i Y(Z) - (-1)^a Y(Z) S^i =
				 [S^i, Y(Z)].
			\end{aligned}
		\end{equation*}
	\item Denote $Y(Z) = \ys(a,Z)_{(j|J)}\ys(b,Z) \vac$ and 
		recall that from Proposition \ref{prop:sesquilinearity}, 
		$\partial_z$ and $\partial_{\theta^i}$ are derivations
		of the $(j|J)$ products. To simplify notation, we will denote $a(Z)
		= \ys(a,Z)$ and $b(Z) = \ys(b,Z)$. We have:
		\begin{multline*}
				S^i Y(W) =  S^i \res_Z \Bigl( i_{z,w} (Z-W)^{j|J}
				a(Z)b(W)\vac - \\  -(-1)^{ab}i_{w,z} (Z-W)^{j|J}
				b(W)a(Z)\vac \Bigr)  
				= \\ =  (-1)^{N + J} \res_Z \Bigl(i_{z,w} (Z-W)^{j|J}
				[S^i,a(Z)]b(W)\vac  \\  + (-1)^a i_{z,w}
				(Z-W)^{j|J}
				a(Z)[S^i, b(W)]\vac - (-1)^{ab}
				i_{w,z} (Z-W)^{j|J}
				[S^i,b(W)]a(Z)\vac -\\ -(-1)^{ab+b} i_{w,z}
				(Z-W)^{j|J}
				b(W) [S^i, a(Z)]\vac \Bigr),
		\end{multline*}
		and, using $S^i \vac = 0$,
		\begin{multline*}
				= (-1)^{N + J} \res_Z \Bigl(i_{z,w} (Z-W)^{j|J}
				(\partial_{\theta^i}a(Z))b(W)\vac + \\  +
				(-1)^a i_{z,w}
				(Z-W)^{j|J}
				a(Z)(\partial_{\zeta^i} b(W))\vac - \\   -
				(-1)^{ab}
				i_{w,z} (Z-W)^{j|J}
				(\partial_{\zeta^i} b(W))a(Z)\vac 
				-(-1)^{ab+b} i_{w,z}
				(Z-W)^{j|J}
				b(W) (\partial_{\theta^i} a(Z))\vac \Bigr) \\
				= (-1)^{N+J} \Bigl((\partial_{\zeta^i}a(W))_{(j|J)}b(W)
				+ (-1)^a
				a(W)_{(j|J)}(\partial_{\zeta^i}b(W))\Bigr)\vac 
				= \\ = \partial_{\zeta^i} \left( a(W)_{(j|J)}b(W)\vac
				\right).
		\end{multline*}
		The proof for $T$ is similar.
\end{enumerate}
\end{proof}
\begin{prop}[Uniqueness]
	Let $V$ be a $N_W=N$ SUSY vertex algebra and let $a(Z)$ be an $\End(V)$-valued
	field such that the pair $(a(Z), \ys(b,Z))$ is local for every $b \in V$, and 
	 $a(Z)\vac = 0$, then $a(Z) = 0$. 
	\label{prop:uniqueness.1}
\end{prop}
\begin{proof}
	By locality there exists  $n \in \mathbb{Z}_+$ such that
	\begin{equation*}
		(z-w)^n a(Z) \ys(b,W)\vac = (-1)^{ab} (z-w)^n \ys(b,W) a(Z) \vac =
		0.
	\end{equation*}
	By Proposition \ref{prop:structure_thm_comienzo} (1), the left hand side is
	$(z-w)^n a(Z) e^{W \nabla}b$. Letting $W = 0$, we get $z^n a(Z)b = 0$, and
	this holds for all $b$, therefore $a(Z) = 0$.
\end{proof}
As a simple corollary of the previous proposition and Proposition
\ref{prop:structure_thm_comienzo} we obtain the following
\begin{thm}
	In an $N_W=N$ SUSY vertex algebra the following identities hold
	\begin{enumerate}
		\item $\ys(a_{(j|J)}b,Z) = \sigma(J)
			\ys(a,Z)_{(j|J)}\ys(b,Z)$ ( $(j|J)$-th
			product identity). 
		\item $\ys(a_{(-1|N)}b,Z) = :\ys(a,Z)\ys(b,Z):$.
		\item $\ys(Ta, Z) = \partial_z \ys(a,Z)$.
		\item $\ys(S^ia, Z) = \partial_{\theta^i}\ys(a,Z)$.
		\item We have the following OPE formula (the sums are finite):
			\begin{equation}
				\begin{aligned}
				{[}\ys(a,Z), \ys(b,W)] &= \sum_{(j|J): j \geq 0}
				\sigma(J)
				(\partial_W^{(j|J)} \delta(Z,W) )
				\ys(a_{(j|J)}b,W) \\
				&= \sum_{(j|J): j \geq 0} (i_{z,w} - i_{w,z})
				(Z-W)^{-1-j|N\setminus
				J} \ys(a_{(j|J)} b, W).
			\end{aligned}
				\label{eq:ope_w.1}
			\end{equation}
	\end{enumerate}
	\label{thm:structure.2}
\end{thm}
\begin{proof}
	(1) is the combined statement of Dong's Lemma \ref{lem:dong}, and Propositions \ref{prop:uniqueness.1} and
	\ref{prop:structure_thm_comienzo} (c). (2) follows from (1) by letting $j|J =
	-1|N$. To prove (3) we write, using (\ref{eq:product.5}), the $(-2|N)$-product
	identity,
	(\ref{eq:negative.1}) and the vacuum axiom:
	\begin{equation*}
		\ys(Ta, Z) = \ys(a_{(-2,N)}\vac, Z) = \ys(a,Z)_{(-2|N)}\Id = 
		:\partial_z \ys(a,Z) \Id: = \partial_z \ys(a,Z).
	\end{equation*}
	(4) follows similarly:
	\begin{multline*}
		\ys(S^ia, Z) = \ys(a_{(-1,N\setminus
		e_i)} \vac, Z) = \\ = - \sigma(N\setminus e_i, e_i)\sigma(e_i,
		N\setminus e_i)(-1)^N
		:\partial_{\theta^i}\ys(a,Z) \Id: =  
		 \partial_{\theta^i} \ys(a,Z).
	\end{multline*}
	Finally (5) follows from (1) and the decomposition Lemma \ref{lem:decomp}
\end{proof}
\begin{cor} Let $e_i = \{i\}$. One has (cf. (\ref{eq:sesquilinearity.1})):
	\begin{equation*}
		\begin{aligned}
			(Ta)_{(j|J)} &= -j a_{(j-1|J)}, \qquad
			(S^ia)_{(j|J)} = \sigma(e_i,N\setminus J)
			a_{(j|J\setminus e_i)},\\
			T(a_{(j|J)}b) &= (Ta)_{(j|J)}b + a_{(j|J)}T(b), \\
			S^i(a_{(j|J)}b) &= (-1)^{N-J}\Bigl( (S^i a)_{(j|J)}b + (-1)^{a}
			a_{(j|J)} S^i b\Bigr).
		\end{aligned}
	\end{equation*}
	\label{cor:tranlsation.relations.1}
\end{cor}
\begin{lem}
	\begin{equation*}
		i_{x,z} \delta(X - Z, W) = i_{w,z} \delta(X, W+Z).
		\label{}
	\end{equation*}
	\label{lem:deltas_are_equal}
\end{lem}
\begin{proof}
	For simplicity let us assume $N=0$, the general result follows easily. Denote:
	\begin{equation*}
		\begin{aligned}
		\psi &= i_{x,z} i_{x-z,w} (x-w-z)^{-1} \in \mathbb{C}[ [x, x^{-1}, z,z^{-1}, w,
		w^{-1}]],\\
		\varphi &= i_{w,z} i_{x, w+z} (x-w-z)^{-1} \in \mathbb{C}[ [x, x^{-1}, z,z^{-1}, w,
		w^{-1}]] 
	\end{aligned}
		\label{}
	\end{equation*}
	It is straightforward to check that both $\psi$ and $\varphi$ are elements of $K[
	[z,w]]$ where $K = \mathbb{C}( (x))$. On the other hand, since both compositions
	$i_{x,z} i_{x-z,w}$ and $i_{w,z} i_{x, w+z}$ commute with multiplication by $x$, $z$
	and $w$, we have
		$(x-w-z) (\psi - \varphi) = 0$,
	hence $\psi = \varphi$, since $K[ [z,w]]$ has no zero divisors.
	Similarly, we have: \[(i_{x,z} i_{w,x-z} - i_{w,z} i_{w+z,x}) (x-w-z)^{-1}
	= 0,\] and the
	lemma follows.
\end{proof}
\begin{nolabel}
	Taking the generating series in Theorem \ref{thm:structure.2}(1) we obtain for the
	left hand side:
	\begin{equation*}
		\sum_{(j|J): j \in \mathbb{Z}} W^{-1-j|N\setminus J} \ys(a_{(j|J)}b,Z) = \ys \Bigl(
		\ys(a,W)b, Z \Bigr).
	\end{equation*}
	On the right hand side we obtain
	\begin{multline*}
		\sum_{(j|J): j \in \mathbb{Z}} W^{-1-j|N\setminus J} \sigma(J) \res_X \Bigl(
		i_{x,z} (X-Z)^{j|J} \ys(a,X) \ys(b,Z) - \\ - (-1)^{ab} i_{z,x}
		(X-Z)^{j|J} b(Z)a(X) \Bigr) = \\ = \res_X \sum_{(j|J): j \in
		\mathbb{Z}} (-1)^{N\setminus J}
		\sigma(J) \Bigl( i_{x,z}
		(X-Z)^{j|J}W^{-1-j|N\setminus J} \ys(a,X) \ys(b,Z) - \\ - (-1)^{ab}
		i_{z,x} (X-Z)^{j|J} W^{-1-j|N\setminus J} \ys(b,Z) \ys(a,X)
		\Bigr).
	\end{multline*}
	But, according to (\ref{eq:fields.proof.6}), this is 
	\begin{equation}
		\res_X \Bigl( i_{x,z} \delta(X-Z,W) \ys(a,X) \ys(b,Z) - (-1)^{ab}
		i_{z,x} \delta(X-Z, W) \ys(b,Z) \ys(a,X) \Bigr).
		\label{eq:quasi-assoc.vertex.2}
	\end{equation}
	Using Lemma \ref{lem:deltas_are_equal}, the first term gives
	\begin{equation}
		i_{w,z} \ys(a,W+Z) \ys(b,Z).
		\label{eq:quasi-assoc.vertex.3}
	\end{equation}
	In order to compute the second term we expand in Taylor series (cf.
	\ref{eq:taylor_expansion_defn})
	\begin{equation*}
		i_{z,x} \delta(X-Z,W) = \sum_{(k|K): k \geq 0} (-1)^K X^{k|K}
		\partial_{-Z}^{(k|K)}
		\delta(-Z,W). 
	\end{equation*}
	Hence the second term in (\ref{eq:quasi-assoc.vertex.2}) reads: 
	\begin{multline*}
		- (-1)^{ab} \res_X  \sum_{(k|K): k \geq 0} (-1)^K X^{k|K} \partial_Z^{(k|K)}
		\delta(-Z,W) \ys(b,Z) X^{-1-n|N\setminus I} a_{(n|I)} = \\ = -
		\res_X \sum_{(k|K): k \geq 0} (-1)^{ab + (N-I)(b+N-K) + K} \sigma(K, N\setminus I)
		 \times \\ \times X^{k-1-n|K \cup(N\setminus I)}
		 \partial_{-Z}^{(k|K)}
		\delta(-Z,W)\ys(b,Z) a_{(n|I)} = \\ = - \sum_{(k|K): k \geq 0}
		(-1)^{(a+N-K)b + N}
		\sigma(K) \partial_{-Z}^{(k|K)} \delta(-Z,W)
		\ys(b,Z) a_{(k|K)}.
	\end{multline*}
	Adding this to (\ref{eq:quasi-assoc.vertex.3}) and changing $Z$ by $-Z$ we
	obtain the important formula
	\begin{multline}
		\ys \Bigl( \ys(a,W)b, -Z \Bigr) = i_{w,z} \ys(a, W-Z) \ys(b, -Z) - \\ -
		\sum_{(k|K): k \geq 0}  (-1)^{(a+N-K)b + N} \sigma(K)
		 \partial_Z^{(k|K)}\delta(Z,W) \ys(b,-Z) a_{(k|K)}.
		\label{eq:quasi-assoc.vertex.6}
	\end{multline}
	Note now that by acting on any vector $c \in V$ and multiplying this last
	equation by a sufficiently high power
	of $(z-w)$ the second term vanishes, therefore we obtain
	\emph{associativity} for the vertex operators, namely:
	\begin{equation}
		(z-w)^n \ys \Bigl( \ys(a,W)b, -Z \Bigr)c = (z-w)^n \ys(a, W-Z)
		\ys(b, -Z)c, \quad n \gg 0.
		\label{eq:quasi-assoc.vertex.7}
	\end{equation}
	\label{no:quas-assoc.vertex}
\end{nolabel}
As in \cite[3.2.3]{frenkelzvi} we obtain an equivalent formulation which is
called the \emph{Cousin} property. Recall the embedding:
\begin{equation*}
	i_{z,w} : \mathbb{C}( (Z, W)) \hookrightarrow \mathbb{C} ( (Z)) ( (W))
\end{equation*}
Given $f \in \mathbb{C}( (Z,W))$, $i_{z,w} f$ is called the \emph{expansion of $f$ in
the domain} $|z| > |w|$. 
\begin{cor}[Cousin property]
	For any $N_W=n$ SUSY vertex algebra $V$ and vectors $a, b, c \in V$,
	the three expressions:
	\begin{equation*}
		\begin{aligned}
			\ys(a,Z) \ys(b, W)c & \in V( (Z) ) ( ( W) ) \\
			(-1)^{ab} \ys(b,W) \ys(a, Z) c & \in V( (W)) ( (Z)) \\
			\ys\left( \ys(a,Z-W)b,W \right)c &\in V( (W) ) ( (Z-W))
		\end{aligned}
	\end{equation*}
	are the expansions, in the domains $|z| > |w|$, $|w| > |z|$ and $|w| > |w-z|$
	respectively, of the same element of 
	$V [ [ Z, W] ] [ z^{-1}, w^{-1}, (z-w)^{-1}]$.
	\label{cor:cousin.w}
\end{cor}
\begin{proof}
	By the locality axiom, there exists $n \in \mathbb{Z}_+$ such that:
	\begin{equation*}
		(z-w)^n \ys(a, Z) \ys(b,W)c = (-1)^{ab}(z-w)^n \ys(b,W)\ys(a,Z)c.
		\label{}
	\end{equation*}
	Since the LHS is an element of $V ( (Z))( (W))$ and the RHS is an element of $V( (W))(
	(Z))$, it follows that they are both equal to some $\varphi \in V[ [Z, W]][z^{-1},
	w^{-1}]$ (cf.
	(\ref{eq:intersection_added})). Since $i_{z,w}$ and $i_{w,z}$ are algebra morphisms, we
	get
	\begin{equation*}
		\ys(a,Z) \ys(b,W)c = i_{z,w} \frac{\varphi}{(z-w)^n}, \quad (-1)^{ab}
		\ys(b,W)\ys(a,Z)c = i_{w,z} \frac{\varphi}{(z-w)^n}.
		\label{}
	\end{equation*}
	The rest of the corollary is proved in a similar way, using (\ref{eq:quasi-assoc.vertex.7}). 
\end{proof}
\begin{thm}[Skew-symmetry]
	In an $N_W=N$ SUSY vertex algebra the following identity, called
	\emph{skew-symmetry}, holds
	\begin{equation}
		\ys(a,Z)b = (-1)^{ab}e^{Z \nabla} \ys(b,-Z)a
		\label{eq:skew-symmetry.w.1}
	\end{equation}
	\label{thm:skew-symmetry.w}
\end{thm}
\begin{proof}
	By the locality axiom we have for $n \gg 0$
	\begin{equation*}
		(z-w)^n\ys(a,Z)\ys(b,W)\vac = (z-w)^n (-1)^{ab}
		\ys(b,W)\ys(a,Z)\vac
	\end{equation*}
	Now by (1) in Proposition \ref{prop:structure_thm_comienzo} we can write
	this as:
	\begin{multline}
			(z-w)^n \ys(a,Z)e^{W\nabla} b = (z-w)(-1)^{ab}
			\ys(b,W)e^{Z\nabla} a \\
			= (z-w)^n (-1)^{ab}
			e^{Z\nabla}e^{-Z\nabla}\ys(b,W)e^{Z\nabla}a 
			= (z-w)^n(-1)^{ab} e^{Z\nabla}i_{w,z}\ys(b,W-Z)a,
		\label{eq:skew.proof.2}
	\end{multline}
	where in the last line we used (2) of Proposition \ref{prop:structure_thm_comienzo}.
	Now both sides in (\ref{eq:skew.proof.2}) are formal power series in $W$.
	Indeed, since $b_{(j|J)}a = 0$ for $j \gg 0$ we see that by making $n$ large
	enough we may assume that there are no negative powers of $w$ in the RHS.
	We can then let $W=0$ in (\ref{eq:skew.proof.2}) and multiply by $z^{-n}$
	to obtain (\ref{eq:skew-symmetry.w.1}). 
\end{proof}
\begin{nolabel}
	Expanding both sides in (\ref{eq:skew-symmetry.w.1}) we have:
	\begin{gather*}
			\sum_{(j|J): j \in \mathbb{Z}} Z^{-1-j|N\setminus J} a_{(j|J)}b =  (-1)^{ab} \Bigl(
			\sum_{(j|J): j \geq 0} 
			\nabla^{(j|J)}Z^{j|J} \Bigr)\times \\\times   \Bigl(
			\sum_{(k|K): k \in \mathbb{Z}}
			(-Z)^{-1-k|N\setminus K} b_{(k|K)}a \Bigr)
			= (-1)^{ab} \sum_{\stackrel{(j|J): j \geq 0}{(k|K): k
			\in \mathbb{Z}}} (-1)^{1+k+N-K}
			\nabla^{(j|J)} \times \\ 
			\quad \times \sigma(J,
			N\setminus K)  Z^{j-1-k|J \cup (N \setminus
			K)} b_{(k|K)}a
	\end{gather*}
	Taking the coefficient of $Z^{-1-n|N\setminus I}$ on both sides we get:
	\begin{multline}
		a_{(n|I)}b = (-1)^{ab} \sum_{j \geq 0, J\cap I = \emptyset} (-1)^{1-n+N+J-I}
		 \times \\ \times (-\nabla)^{(j|J)}
		\sigma(N\setminus(I \cup J),J)  b_{(n+j|I\cup J)}a.
		\label{eq:skew.cor.2}
	\end{multline}
	In particular, when $(n|I)=(-1|N)$ in (\ref{eq:skew.cor.2}), we get:
	\begin{equation*}
		:ab: - (-1)^{ab} :ba: = (-1)^{ab} \sum_{j \geq 1}
		\frac{(-T)^{j}}{j!} \left(
		b_{(-1+j|N)}a\right),
	\end{equation*}
	or, equivalently, after exchanging $a$ and $b$:
	\begin{equation}
		:ab: - (-1)^{ab} :ba: = \int_{-\nabla}^0 [a_\Lambda b] d\Lambda.
		\label{eq:skew.cor.agregado}
	\end{equation}
	Identity (\ref{eq:skew.cor.agregado}) is called the \emph{
	quasi-commutativity} of the normally ordered product.
%
%
	\label{no:skew.cor}
\end{nolabel}
\begin{nolabel}
	Define the following \emph{formal Fourier transform} by 
	\begin{equation*}
		F^\Lambda_Z a(Z) = \res_Z e^{Z \Lambda} a(Z).
		\label{}
	\end{equation*}
	It is a linear map from the space of $\cU$-valued formal distributions in $Z$ to
	$\cU[ [\Lambda]]$. It has the following properties which are immediate to check:
		\begin{align}
			F^\Lambda_Z \partial_z a(Z) &= - \lambda F^\Lambda_Z a(Z),\label{al.1p} \\
			F^\Lambda_Z \partial_{\theta^i} a(Z) &= - (-1)^N \chi^i F^\Lambda_Z
			a(Z), \label{al.2p} \\
			F^\Lambda_Z \left(e^{Z\nabla} a(Z) \right) &= F^{\Lambda + \nabla}_Z
			a(Z) \; \text{if } a(Z) \in \cU( (Z)), \label{al.3p} \\
			F^\Lambda_Z a(-Z) &= - F^{-\Lambda}_Z a(Z),\label{al.4p} \\ 
			F^\Lambda_Z \left( \partial^{(j|J)}_W \delta(Z,W) \right) &= (-1)^{JN}
			e^{W\Lambda} \Lambda^{(j|J)}. \label{al.5p}
		\end{align}
	\label{no:fourier_1variable}
\end{nolabel}
\begin{thm}
	Let $V$ be a $N_W=N$ SUSY vertex algebra. Then  $V$ is a $N_W = N$ SUSY
	Lie conformal 
	algebra with $\Lambda$-bracket:
	\begin{equation}
		[a_\Lambda b] = F^\Lambda_Z\ys(a,Z)b = \sum_{(j|J): j \geq 0}
		(-1)^{JN} \sigma(J)
		 \Lambda^{(j|J)} (a_{(j|J)}b).
		\label{eq:conformal_vertex.1}
	\end{equation}
	\label{thm:conformal_vertex.1}
\end{thm}
\begin{proof}
	The sesquilinearity relations follow from Corollary \ref{cor:tranlsation.relations.1}
	for $j \geq 0$. Applying $F^\Lambda_Z$ to both sides of (\ref{eq:skew-symmetry.w.1})
	and using (\ref{al.3p}) and (\ref{al.4p}) we get the skew-symmetry relation. In order
	to prove the Jacobi identity, apply $F^\Lambda_Z$ to the OPE formula (\ref{eq:ope_w.1})
	applied to $c$, and use (\ref{al.5p}) to obtain (cf. (\ref{eq:wick.1})):
	\begin{equation*}
		[a_\Lambda \ys(b,W)c] = (-1)^{ab + bN} \ys(b,W) [a_\Lambda c] + e^{W \Lambda}
		\ys([a_\Lambda b], W)c.
		\label{}
	\end{equation*}
	Applying $F^\Gamma_W$ to both sides of this formula we get the Jacobi identity. 
\end{proof}
\begin{thm}
	Let $V$ be a $N_W=N$ SUSY vertex algebra. The following identity called
	``quasi-associativity'' of the normally ordered product holds for every
	 $a, b, c \in V$:
	\begin{equation*}
		::ab:c: - :a:bc:: = \sum_{j \geq 0} a_{(-2-j|N)}\left(b_{(j|N)}c\right) +
		(-1)^{ab}  \sum_{j \geq 0} b_{(-2-j|N)} \left(a_{(j|N)}c\right).
	\end{equation*}
	Equivalently
	\begin{equation*}
		: :ab:c: - :a :bc: : =  \left( \int_0^\nabla d\Lambda a \right)[b_\Lambda c]  +
		(-1)^{ab} \left( \int_0^\nabla d\Lambda b \right) [a_\Lambda c],
	\end{equation*}
	where the integral is computed as follows: expand the $\Lambda$-bracket, put the powers
	of $\Lambda$ on the left, under the sign of integral, then take the definite integral
	by the usual rules inside the parenthesis. 
	\label{thm:quasi-associativity1}
\end{thm}
\begin{proof}
	Applying both sides of Theorem \ref{thm:structure.2} (2) to $c$ and taking the constant
	coefficient, 
	 the LHS is $::ab:c:$.  By (\ref{eq:normal_defin}), the RHS of
	Theorem \ref{thm:structure.2} (2) applied to $c$ is
	\begin{multline}
		\sum_{\stackrel{j < 0, J}{k, K \cup J = N}} (-1)^{(N-K)(a + N - J)}
		\sigma(N\setminus J, N\setminus K) Z^{-2-j-k|N\setminus (J \cap K)}
		a_{(j|J)} \left( b_{(k|K)}c \right) + \\ + \sum_{\stackrel{j \geq 0, J}{k, K \cup J
		= N}} (-1)^{(N-J)(b+ N - K) + ab} \sigma(N\setminus K, N\setminus J)
		Z^{-2-j-k|N\setminus (J \cap K)} b_{(k|K)}\left( a_{(j|J)}c \right).
		\label{}
	\end{multline}
	To compute the constant coefficient in the last formula, we let $K = J = N$, and $k =
	-2-j$, to get
	\begin{equation*}
		\sum_{j \geq -1} a_{(-2-j|N)}\left( b_{(j|N)}c \right) + (-1)^{ab} \sum_{j \geq
		0} b_{(-2-j|N)} \left( a_{(j|N)}c \right).
		\label{}
	\end{equation*}
	Noting that the term with $j = -1$ in the first summand in the last
	formula is
	$:a:bc: :$, the theorem follows.
\end{proof}
We thus arrive to the following equivalent definition of an $N_W=N$ SUSY vertex algebra
(cf. \cite{kacbakalov2}):
\begin{defn}
	An $N_W=N$ SUSY vertex algebra is a tuple $(V$, $T$, $S^i$,
	$[\cdot_\Lambda\cdot]$, $\vac$, $::)$, 
	where
	\begin{itemize}
		\item $(V, T, S^i, [\cdot_\Lambda \cdot] )$ is an $N_W=N$ SUSY Lie conformal
			algebra,
		\item $(V, \vac,T, S^i, : :)$ is a unital quasicommutative quasiassociative
			differential superalgebra (i.e. $T$ is an even derivation
			of $: :$ and $S^i$ ($i = 1, \dots, N$) are odd derivations of $: :$),
		\item the $\Lambda$-bracket and the product $: :$ are related by the
			non-commutative Wick formula (\ref{eq:wick.proof.8}). 
	\end{itemize}
	\label{defn:bakalov}
\end{defn}
\begin{proof}
	We have shown that this definition follows from Definition
	\ref{defn:w_defin}. For the converse, we refer the reader to \cite{kacbakalov2}. The
	proof carries over to the SUSY case with minor modifications. 
\end{proof}
Removing the ``quantum corrections'' we arrive to the following definition:
\begin{defn}
	An $N_W = N$  \emph{Poisson} SUSY vertex algebra is tuple
	 $(V$, $\vac$, $T$, $S^i$, $\{\cdot_\Lambda \cdot\}$, $\cdot)$, where
	 \begin{itemize}
		 \item $(V, T, S^i, \{\cdot_\Lambda \cdot\})$ is an $N_W=N$ SUSY Lie conformal
			 algebra,
		 \item $(V, \vac, T, S^i, \cdot)$ is an unital commutative associative
			 differential superalgebra,
		\item the following \emph{Leibniz rule} is satisfied:
			\begin{equation*}
				\{a_\Lambda bc\} = \{a_\Lambda b\} c + (-1)^{(a+N)b} b
				\{a_\Lambda c\}.
			\end{equation*}
	\end{itemize}
	\label{defn:poisson_va}
\end{defn}
\begin{thm}
	Let $V$ be an $N_W=N$ SUSY vertex algebra. For each  $a, b \in V$, $k \in \mathbb{Z}$
	and $K \subset \{1,\dots, N\}$, the following
	identity, called \emph{Borcherds identity}, holds:
	\begin{multline}
		\left(i_{z,w} (Z-W)^{k|K} \right)\ys(a,Z) \ys(b,W) - (-1)^{ab} \left(i_{w,z}
		(Z-W)^{k|K}\right) \ys(b,W) \ys(a,Z) = \\ = \sum_{j \geq 0, J} \sigma(J,K)
		\sigma(J \cup K) \left(\partial_W^{(j|J)}
		\delta(Z,W) \right) \ys(a_{(k+j|K\cup
		J)}b, W).
		\label{eq:borcherds.identity.w.1}
	\end{multline}
	\label{thm:borcherds.identity.w}
\end{thm}
\begin{proof}
	The LHS of (\ref{eq:borcherds.identity.w.1}) is local since multiplied by $(z-w)^n$ for
	$n \gg 0$ it
	is equal to
	\begin{equation*}
		(Z-W)^{n+k|K} [\ys(a,Z), \ys(b,W)] = 0,
	\end{equation*}
	by the locality axiom. Therefore we can apply the decomposition Lemma \ref{lem:decomp}
	to the LHS of (\ref{eq:borcherds.identity.w.1}). We have
	\begin{multline*}
		c_{j|J}(W) = \sigma(J, K) \res_Z \left( \left( i_{z,w} (Z-W)^{k+J|K \cup J} \right)
		\ys(a, Z) \ys(b,W) - \right. \\ \left. - (-1)^{ab} \sigma(J,K) \left( i_{w,z}
		(Z-W)^{k+j|K \cup J}
		\right) \ys(b,W) \ys(a,Z) \right),
	\end{multline*}
	therefore the theorem follows from (\ref{eq:negative.4}) and Theorem
	\ref{thm:structure.2} (1).
\end{proof}
\begin{prop}
	Let $V$ be a $N_W=N$ SUSY vertex algebra. Then 
	\begin{multline}
		[a_{(n|I)}, \ys(b,W)] = \sum_{(j|J): j\geq 0}
		(-1)^{JN + IN + IJ}  \sigma(J) \sigma(I)\times \\ \times
		\bigl(\partial_W^{(j|J)} W^{n|I}\bigr)
		\ys\left( a_{(j|J)}b,W \right). \label{eq:commutator.1}
	\end{multline}
	If, moreover, $n \geq 0$, this becomes:
	\begin{equation}
		[a_{(n|I)}, \ys(b,W)] = \ys(e^{-W\nabla}a_{(n|I)}e^{W\nabla}b,W).
		\label{eq:commutator.2}
	\end{equation}	
	\label{prop:commutator_w.1}
\end{prop}
\begin{proof}
	Multiplying the OPE formula (\ref{eq:ope_w.1}) by $Z^{n|I}$ and taking
	residues we obtain in the left hand side 
	$	\sigma(I) [a_{(n|I)}, \ys(b,W)]$,  
	while the right hand side is
	\begin{multline*}
		\res_Z \sum_{(j|J): j \geq 0} (-1)^{I(N-J)} 
		\sigma(J)  (\partial_W^{(j|J)} \delta(Z,W) Z^{n|I})
		\ys(a_{(j|J)}b,W) = \\ = \res_Z \sum_{(j|J): j \geq 0}
		(-1)^{I(N-J)}
		\sigma(J)  (\partial_W^{(j|J)} \delta(Z,W) W^{n|I})
		\ys(a_{(j|J)}b,W) = \\ = \sum_{(j|J): j \geq 0}
		(-1)^{I(N-J)+JN}
		\sigma(J)  (\partial_W^{(j|J)} W^{n|I})
		\ys(a_{(j|J)}b,W)   
	\end{multline*}
	hence (\ref{eq:commutator.1}) follows. Note that when $n \geq 0$, the RHS of
	(\ref{eq:commutator.1}) is a finite sum of fields of $V$ times monomials $W^{k|K}$ with
	$k \geq 0$.
	 This easily implies that this field is local
	with respect to all fields of $V$, hence the LHS of (\ref{eq:commutator.2}) is local
	with respect to all fields of $V$. On the other hand, the adjoint action of
	$\nabla$ on $a_{(n|I)}$ either decreases $n$ or $\sharp I$ (cf. (\ref{eq:product.5})).
	Using the the formula $\mathrm{Ad} e^X = e^{\mathrm{ad}X}$ for an even element $X$ of a Lie
	superalgebra, we get
	that $e^{-W\nabla} a_{(n|I)} e^{W\nabla}$ is a finite sum, involving
	only positive powers of $w$, hence the
	RHS of (\ref{eq:commutator.2}) is also local with respect to all fields of $V$.
	
	To apply the uniqueness theorem, we need to check that
	both sides agree when valuated at the vacuum vector. The left hand side is
	given by
	\begin{equation*}
		[a_{(n|I)}, \ys(b,W)]\vac = a_{(n|I)} e^{W \nabla} b,
	\end{equation*}
	where we used the fact that $a_{(n|I)}\vac = 0$ and Proposition
	\ref{prop:structure_thm_comienzo} (1). On the other hand, by the same
	proposition the left hand side is
	\begin{equation*}
		\ys(e^{-W \nabla} a_{(n|I)} e^{W\nabla}, W)\vac =
		e^{W\nabla}e^{-W\nabla} a_{(n|I)} e^{W\nabla}\vac,
	\end{equation*} 
	and (\ref{eq:commutator.2}) follows.
\end{proof}
\begin{rem}
	As a consequence of (\ref{eq:commutator.1}) we see that by taking the
	coefficient of $W^{-1-k|N\setminus K}$ we obtain the commutator
	$[a_{(n|I)}, b_{(k|K)}]$ as a linear combination of Fourier modes of fields
	in $V$. This rather complicated formula says that the linear span of
	Fourier modes of $\End(V)$-valued fields is a Lie superalgebra. In order to compute
	explicitly the Lie bracket, we compute the coefficient of
	$W^{-1-k|N\setminus K}$ on the left hand side of (\ref{eq:commutator.1}) to
	obtain:
	\begin{equation}
		(-1)^{(a+N-I)(N-K)} [a_{(n|I)},b_{(k|K)}].
		\label{eq:lie_bracket.1}
	\end{equation}
	To compute this coefficient on the right hand side we first expand:
	\begin{multline}
		(\partial_{W}^{j|J} W^{n|I})W^{-1-l|N\setminus L} =
		\frac{(-1)^{\frac{J(J-1)
		}{2}} n!}{(n-j)!} \times \\ \times \sigma(J, I\setminus J)
		\sigma(I\setminus J, N\setminus L) W^{n-j-1-l|(I
		\setminus
		J)\cup (N\setminus L)}.
		\label{eq:lie_bracket.2}
	\end{multline}
	Note that in order for the corresponding term in (\ref{eq:commutator.1})
	not to vanish, we must have $J \subset I$ and in order for the coefficient
	of $W^{-1-k|N\setminus K}$ not to be zero in (\ref{eq:lie_bracket.2})  we
	must have $(K \cap I) \subset J$. Now we set $n-j-l-1 = -1-k$ and
	$(I\setminus J)\cup (N\setminus L) = N\setminus K$ to obtain $l =
	n+k-j$ and $L = K \cup (I \setminus J)$. We get then for the right hand
	side
	\begin{multline*}
		\sum_{(j|J): j \geq 0} (-1)^{(J + I)(N-J)} \binom{n}{j} \sigma(J)
		\sigma(I) \times \\ \times \sigma(J, I\setminus J)
		\sigma(I\setminus J, (N\setminus
		K)\setminus (I\setminus J)) \bigl(a_{(j|J)}b \bigr)_{(n+k-j|K\cup
		(I\setminus J))}
	\end{multline*}
	Combining with (\ref{eq:lie_bracket.1}), we obtain:
	\begin{multline}
		[a_{(n|I)},b_{(k|K)}] = (-1)^{(a+N-I)(N-K)} \sum_{(j|J): j \geq 0}
		(-1)^{(I-J)(N-J)} \binom{n}{j} \sigma(J) \times \\ \times
		\sigma(I)
		 \sigma(J,I\setminus J) \sigma\left(I\setminus J, (N\setminus K) \setminus
		(I\setminus J)\right)
		\bigl(a_{(j|J)}b\bigr)_{(n+k-j|K\cup (I\setminus
		J))}.
		\label{eq:lie_bracket.4}
	\end{multline}
	\label{rem:lie_bracket}
\end{rem}
\begin{nolabel}
	We can define the \emph{tensor product} of two $N_W=N$ SUSY vertex algebras in
	the usual way, namely, let $V$ and $W$ be two $N_W=N$ SUSY vertex algebras.
	The
	space of states is the vector superspace $V\otimes W$. The vacuum vector is
	$\vac_V \otimes \vac_W$. Let us denote $\ys_V$ and $\ys_W$ the
	corresponding state-field correspondences. We define the state field
	correspondence $\ys$ for $V\otimes W$ as 
	\begin{multline}
		\ys(a\otimes b, Z) = \ys_V(a,Z) \otimes \ys_W(b, Z) = \\
		=\sum_{(j|J), (k|K)}
		(-1)^{a(N-K)} \sigma(N\setminus
		K,N\setminus J) Z^{-2-j-k|(N \setminus (J\cap K))} a_{(j|J)}\otimes
		b_{(k|K)},
		\label{eq:tensor.1}
	\end{multline}
	where the endomorphism $a_{(j|J)}\otimes b_{(k|K)}$ is defined to be
	\begin{equation*}
		(a_{(j|J)} \otimes b_{(k|K)}) (v \otimes w) = (-1)^{(b+N-K) v} a_{(j|J)}
		v \otimes b_{(k|K)} w.
	\end{equation*}
	Note that in order for $\sigma$ not to vanish in (\ref{eq:tensor.1}) we
	must have $J \cup K = N$. Finally, we let the translation operators be $T =
	T_V \otimes Id + Id \otimes T_W$ and $S^i = S^i_V \otimes Id + Id \otimes
	S^i_W$. All the axioms of  SUSY vertex algebra are straightforward to check.

	\label{no:tensor}
\end{nolabel}
\begin{thm}[Existence] Let $V$ be a vector superspace, $\vac \in V$ an even vector, $T$
	an even endomorphism of $V$ and $S^i$, $i = 1, \dots, N$, odd endomorphisms
	of $V$, pairwise anticommuting between themselves and commuting with $T$. Suppose
	moreover that
	$T\vac = S^i \vac =0$. Let $\cF$ be a family
	of $\End(V)$-valued fields 
	\begin{equation*}
		a^\alpha(Z) = \sum_{j \in \mathbb{Z}, J} Z^{-1-j|N\setminus J} a^\alpha_{(j|J)}
	\end{equation*}
	indexed by $\alpha \in A$, such that 
	\begin{enumerate}
		\item $a^\alpha(Z)\vac|_{Z = 0} = a^\alpha \in V$,
		\item $[T, a^\alpha(Z)] = \partial_z a^\alpha(Z)$ and $[S^i,
			a^\alpha(Z)] = \partial_{\theta^i} a^\alpha(Z)$, 
		\item all pairs $(a^\alpha(Z), a^\beta(Z))$ are local,
		\item the vectors 
				$a^{\alpha_s}_{(j_s|J_s)} \dots
				a^{\alpha_1}_{(j_1|J_1)} \vac$ span $V$.
	\end{enumerate}
	Then the formula
	\begin{multline}
		\ys(a^{\alpha_s}_{(j_s|J_s)} \dots a^{\alpha_1}_{(j_1|J_1)} \vac,
		Z) = \\ =\prod \sigma(J_i)
		a^{\alpha_s}(Z)_{(j_s|J_s)} \Bigl( \dots a^{\alpha_2}_{(j_2|J_2)}
		\bigl(a^{\alpha_1}_{(j_1|J_1)} \Id\bigr) \dots \Bigr)
		\label{eq:existence.3}
	\end{multline}
	gives a well defined structure of an $N_W=N$ SUSY vertex algebra on $V$, with vacuum vector
	$\vac$, translation operators $T, \, S^i$, and such that
	$	\ys(a^\alpha, Z) = a^\alpha(Z)$.

	Such a structure is unique.
	\label{thm:existence}
\end{thm}
\begin{proof}
	Let $\bar{\cF}$ be the minimal family of  $\End(V)$-valued
	fields containing
	$\cF$, closed under all $(j|J)$-products and under the derivations
	$\partial_z$ and $\partial_{\theta^i}$, and subject to the conditions (1)-(3) of the
	Theorem. By Theorem \ref{thm:fields.1},
	$\bar{\cF}$ is an $N_W=N$ SUSY vertex algebra. Define a map $\varphi:
	\bar{\cF} \rightarrow V$ by $a(Z) \mapsto a(Z)\vac_{Z = 0}$. This map is
	surjective by (4). Let $a(Z) \in \ker \varphi$. It follows easily by (2) and the fact
	that $T\vac = S^i\vac = 0$ that $a(Z)\vac
	= 0$. Since $\varphi$ is surjective, for each
	$b \in V$, there exists $b(W)  \in \varphi^{-1} (b)$. By $(3)$ there
	exists  $j \in \mathbb{Z}_+$ such that
	\begin{equation*}
		(z-w)^{j} a(Z) b(W) \vac = (-1)^{ab} (z-w)^j b(W) a(Z) \vac = 0
		\label{}
	\end{equation*}
	Letting $W=0$ and canceling $z^j$ we obtain $a(Z)b = 0$ for every $b \in V$, hence
	$a(Z) = 0$ and $\varphi$ is also injective. 	
	 We obtain thus a \emph{state-field correspondence}
	$\ys : a \mapsto \ys(a,Z)$. Formula (\ref{eq:existence.3}) follows from the
	$(j|J)$-product identity in Theorem \ref{thm:structure.2} (1). 
\end{proof}
\subsection{The universal enveloping SUSY vertex algebra} \label{sub:universal}
In this section we construct maps $\varphi$ and $\varphi'$, used in \cite{heluani4} to define
the conformal blocks, and we
construct an $N_W = N$ SUSY vertex algebra attached to each $N_W = N$ SUSY
Lie conformal algebra.
\begin{nolabel} Let $V$ be an $N_W = N$ SUSY vertex algebra. According to Theorem
	\ref{thm:conformal_vertex.1} it is a SUSY Lie conformal algebra. 
	It follows by
	Proposition \ref{prop:formal_distrib_from_conf} that the pair
	$(\lie (V), V)$ is an $N_W = N$ formal distribution Lie superalgebra. 
	
	Recall from Lemma \ref{lem:affinization_1} and \ref{no:formal_from_R_1} the
	construction of the Lie superalgebra $\lie(V) = \tilde{V}/\tilde{\nabla}
	\tilde{V}$, where $\tilde{V} = V
	\otimes_{\mathbb{C}}
	\mathbb{C}[X, X^{-1}]$ and $\tilde{\nabla} \tilde{V}$ is the space spanned by vectors
	of the form:
	\begin{equation}
			T a \otimes f(X) + a \otimes \partial_x f(X), \quad
			S^i a \otimes f(X) + (-1)^{aN} a \otimes \partial_{\eta^i} f(X),
		\label{eq:acabodeagregarte}
	\end{equation}
	for $a \in V$, $f(X) \in \mathbb{C}[X, X^{-1}]$, and we change the
	parity if $N$ is odd. 
	
	Let $\varphi : \lie (V) \rightarrow \End(V)$ be the linear map
	defined by 
	\begin{equation}
		a_{<n|I>} = a \otimes X^{n|I} \mapsto (-1)^{aI}\sigma(I)
		a_{(n|I)}, \qquad a \in V.
		\label{eq:nuevamentealvicio}
	\end{equation}
	Similarly, we construct $V \otimes_{\mathbb{C}} \mathbb{C}( (X))$ and consider its
	quotient $\lie'(V)$ by the vector space generated by vectors of the form
	(\ref{eq:acabodeagregarte}), with reversed parity if $N$ is odd. Then (\ref{eq:nuevamentealvicio}) defines a map $\varphi' : \lie'(V) \rightarrow
	\End(V)$.  Comparing (\ref{eq:lie_bracket.4aa}) and (\ref{eq:lie_bracket.4}) and noting the extra factor $\sigma(J)$ in
	(\ref{eq:conformal_vertex.1}) we obtain
	the following 
\end{nolabel}
\begin{thm}
	The maps $\varphi$, and $\varphi'$ are Lie superalgebra homomorphisms.
	\label{thm:lie_fourier.vertex.1}
\end{thm}
\begin{nolabel}
	Let $\cR$ be an $N_W = N$ SUSY Lie conformal algebra, and let $(\lie(\cR), \cR)$ be the
	corresponding $N_W = N$ formal distribution Lie superalgebra (cf. Proposition
	\ref{prop:formal_distrib_from_conf}). The Lie bracket in $\lie(\cR)$ is
	given by
	(\ref{eq:lie_bracket.4aa}).  Recall from \ref{no:regular_lie_R} that 
	$\lie(\cR)$ is a regular $N_W=N$ formal distribution Lie superalgebra.
	In particular, it carries an even derivation $T$ and $N$ odd derivations
	$S^i$, $i= 1, \dots, N$ defined by (\ref{eq:derivations_lie_R}). Moreover, 
	the annihilation subalgebra $\lie(\cR)_\leq$ is invariant by these
	derivations.
	\label{no:lie_from.conf.1}
\end{nolabel}
\begin{thm}
	Let $\cR$ be an $N_W = N$ SUSY Lie conformal algebra. Let $V=V(\cR)$ be the quotient of
	$U(\lie(\cR))$ by the left ideal generated by $\lie(\cR)_\leq$. Then $V$
	admits an $N_W=N$ SUSY vertex algebra structure whose vacuum vector 
	$\vac$ is the image
	of $1$ in $V$, 	and the translation operators $T$, $S^i$ $(i=1,\dots,
	N)$, are obtained by extending
	the corresponding derivations on $\lie(\cR)$ by the Leibniz rule. 
	 This vertex algebra is called the \emph{universal enveloping
	vertex algebra of $\cR$}. 
	\label{thm:universal.1}
\end{thm}
\begin{proof}
	Let $\cF$ be the family of $\lie(\cR)$-valued formal distributions 
\begin{equation*}
	\cF = \left\{ a(Z) \Bigl| a \in \mathbb{\cR} \right\},
\end{equation*}
	where $a(Z)$ was defined in 
	(\ref{eq:distributions_in_lie_R}). Note that this family defines a family of
	$\End(V)$-valued formal distributions, where the action is by left multiplication. This
	family satisfies (1)-(4) of Theorem \ref{thm:existence}.  Indeed, 
	(2) follows from (\ref{eq:lie_R_is_regular}), (3) follows from Proposition
	\ref{prop:distributions_are_local} and (4) follows since the vectors $a_{(j|J)}$ with
	$a \in \cR$, $j \in \mathbb{Z}$ and $J \subset \{1,\dots, N\}$ span $\lie(\cR)$. The
	theorem will follow from the existence Theorem 
	\ref{thm:existence} if we show that these distributions are in fact $\End(V)$-valued
	fields.

	For that, given $a^{\alpha_1}, \dots, a^{\alpha_s} \in \cR$, we have to
	prove that for any $a \in \cR$ and $J \subset \{1, \dots, N\}$, we have:
	$	a_{(n|J)} a^{\alpha_1}_{(j_1|J_1)} \dots
		a^{\alpha_s}_{(j_s|J_s)}\vac = 0$ for  $n \gg 0$.	
	This is proved by induction on $s$, using (\ref{eq:lie_bracket.4aa}). 
\end{proof}
\numberwithin{thm}{section}
\section{Structure theory of $N_K=N$ SUSY vertex algebras} \label{sub:k.formal}
\begin{nolabel}
	In this section we develop the structure theory of \emph{
	$N_K=N$ SUSY vertex algebras}, where $N$ is a positive integer. This kind of structures has been studied, in some
	particular cases, in the
	physics literature. Roughly speaking an $N_K=N$ SUSY vertex algebra is an
	$N_W=N$ SUSY vertex algebra, where the differential operators
	$\partial_{\theta^i}$ are replaced by the
	differential operators 	
	\begin{equation*}
		D^i_Z = D^{e_i}_Z = \partial_{\theta^i} + \theta^i \partial_z.
	\end{equation*}

	To describe the corresponding \emph{SUSY Lie conformal superalgebras}, perhaps
	the language of \emph{$H$-pseudoalgebras} is more convenient
	\cite{dandrea:seudo}.
	On the other hand, we are interested in their universal
	enveloping vertex algebras and in particular we want a description along the lines of
	the previous sections. 

	In order to have a uniform notation between this section and the previous
	ones, given two sets of coordinates $Z = (z,\theta^i)$ and $W = (w, \zeta^i)$ we
	will denote
	\begin{equation}
		\begin{gathered}
		Z-W = \left( z-w - \sum_{i=1}^N \theta^i \zeta^i, \theta^j -
		\zeta^j \right), \\ (\theta - \zeta)^J = \prod_{i \in J} (\theta^i
		- \zeta^i), \quad (Z-W)^{j|J} = \left( z-w-\sum_{i=1}^N \theta^i\zeta^i
		\right)^j (\theta - \zeta)^J,
	\end{gathered}
		\label{eq:k.intro.2}
	\end{equation}
	where $j \in \mathbb{Z}$ and $J$ is an ordered subset of $\{1, \dots, N\}$. 
	As before, we define \[Z^{j|J} = z^j \theta^J.\] Note that
	\begin{equation}
		(Z-W)^{-1|0} = \sum_{k = 0}^N \frac{\left( \sum_{i=1}^N \theta^i
		\zeta^i \right)^{k}}{(z-w)^{k+1}},
		\label{eq:k.intro.3}
	\end{equation}
	therefore $(Z-W)^{-1|N}$ coincides with that in the $N_W=N$ case:
	\begin{equation}
		(Z-W)^{-1|N} = \frac{(\theta-\zeta)^N}{z-w}.
		\label{eq:k.intro.4}
	\end{equation}

	The differential operators $D^i_Z$ satisfy the commutation relations
	\begin{equation}
		{[}D^i_Z, D^j_Z ] = 2 \delta_{i,j} \partial_z,
		\label{eq:k.intro.5}
	\end{equation}
	and, as before, we denote for $J = (j_1, \dots, j_k)$:
	\begin{equation}
		D_Z = (\partial_z, D^1_Z, \dots, D^N_Z), \quad D^{j|J}_Z = \partial_z^j
		D^{j_1}_Z \dots 
		D^{j_k}_Z, \quad D^{(j|J)}_Z = \frac{(-1)^{\frac{J(J+1)}{2}}}{j!}
		D^{j|J}_Z.
		\label{eq:k.intro.6}
	\end{equation}
	Occasionally,  when $j=0$, we will write $D^{0|J}_Z =
	D^J_Z$.

	Finally, in this section we will consider not necessarily disjoint subsets $I, J
	\subset \{1, \dots, N\}$ as in the
	$N_W=N$ case. Given $I$ and $J$, ordered subsets of $\{1, \dots, N\}$, we will write $I
	\triangle J = (I\setminus J) \cup (J \setminus I)$. Note, however, that
	still $(Z-W)^{j|J} (Z-W)^{k|K} = 0$ if $J \cap K \neq \emptyset$.  We will use the same formal
	$\delta$-function $\delta(Z,W)$ as before. Remarkably, the new binomial $(Z-W)^{j|J}$,
	given by (\ref{eq:k.intro.2}) ``behaves'' with respect to the operators $D^{j|J}_W$, in
	the same way as the old binomials (\ref{eq:w.intro.recontra}) with respect to
	$\partial_W^{j|J}$. 
	 \label{no:k.intro}
\end{nolabel}
\begin{lem} The following identity is true:
	\begin{equation}
		D^{(j|J)}_W \delta(Z,W) =  \sigma(J)
	 (i_{z,w} - i_{w,z}) (Z-W)^{-1-j|N\setminus J}.
		\label{eq:k.deriv.1}
	\end{equation}
	\label{lem:k.deriv}
\end{lem}
\begin{proof}
	Let us assume for simplicity that $j=0$, the general case follows easily
	from this, differentiating by $w$.  We will prove the lemma by induction on $\sharp J$. When $J=
	\emptyset$, it follows from (\ref{eq:k.intro.4}) that (\ref{eq:k.deriv.1})
	coincides with the formula (\ref{eq:intro.3.1}) for $\delta(Z,W)$. When $J = e_i = \{i\}$, the left
	hand side of (\ref{eq:k.deriv.1}) is given by
	\begin{equation*}
		\begin{aligned}
			-D^{i}_W \delta(Z,W) &= -D^i_W (i_{z,w} - i_{w,z})
			\frac{(\theta-\zeta)^N}{z-w}
			\\ 
			&= - (i_{z,w} - i_{w,z}) \bigl( - \sigma(e_i)
			\frac{(\theta -
			\zeta)^{N\setminus e_i}}{z-w} + \zeta^i
			\frac{(\theta - \zeta)^N}{(z-w)^2} \bigr)
		\end{aligned}
	\end{equation*}
	On the other hand, using (\ref{eq:k.intro.3}), we get:
	\begin{equation*}
		\begin{aligned}
			(Z-W)^{-1|N\setminus e_i} &= \sum_{k \geq 0}
			\frac{\left( \sum \theta^i \zeta^i \right)^k}{(z-w)^{k+1}}
			(\theta - \zeta)^{N\setminus e_i} \\ &=
			\frac{(\theta - \zeta)^{N\setminus e_i}}{z-w} +
			\frac{\theta^i \zeta^i}{(z-w)^2} (\theta -
			\zeta)^{N\setminus e_i} \\ 
			&= \frac{(\theta - \zeta)^{N\setminus e_i}}{{z-w}} -
			\sigma(e_i, N\setminus e_i) \frac{\zeta^i}{(z-w)^2}
			(\theta - \zeta)^N,
		\end{aligned}
	\end{equation*}
	hence (\ref{eq:k.deriv.1}) follows when $J=e_i$. To prove the general
	case, assume that the lemma is valid for $J = I\setminus e_i$. Since
	$D^I_W = \sigma(e_i, I\setminus e_i) D^i_W D^{I\setminus e_i}_W$ we have by
	the induction hypothesis
	\begin{multline}
		D^I_W \delta(Z,W) = \sigma(e_i, I\setminus e_i) \sigma(I\setminus
		e_i, N\setminus (I\setminus e_i)) (-1)^{\frac{(I-1)I}{2}} \times \\
		\times D^i_W
		(i_{z,w} - i_{w,z}) (Z-W)^{-1|N\setminus(I\setminus e_i)}.
		\label{eq:deriv.proof.3}
	\end{multline}
	We expand the last factor as:
	\begin{multline}
		D^i (Z-W)^{-1|N\setminus (I\setminus e_i)} = - \sum_{k \geq 1} k
		\zeta^i \frac{\left( \sum \theta^j \zeta^j
		\right)^{k-1}}{(z-w)^{k+1}} (\theta - \zeta)^{N\setminus (I
		\setminus e_i)} - \\ - \sigma(e_i, N \setminus I) \sum_{k \geq 0}
		\frac{\left( \sum \theta^j \zeta^j \right)^k}{(z-w)^{k+1}} (\theta
		- \zeta)^{N\setminus I} + \sum_{k \geq 0} (k+1) \zeta^i
		\frac{\left( \sum \theta^j \zeta^j \right)^{k}}{(z-w)^{k+2}}
		(\theta - \zeta)^{N\setminus (I\setminus e_i)}.
		\label{eq:deriv.proof.4}
	\end{multline}
	Relabeling the indexes we see that the first and last term cancel. Finally
	we note that, by (\ref{eq:not.2b}):
	\begin{equation}
		\sigma(e_i, I \setminus e_i) \sigma(e_i, N\setminus I)
		\sigma(I\setminus e_i) = (-1)^I
		\sigma(I).
		\label{eq:deriv.proof.5}
	\end{equation}
	Combining (\ref{eq:deriv.proof.5}), (\ref{eq:deriv.proof.4}) and
	(\ref{eq:deriv.proof.3}) we obtain the lemma for $j = 0$.
\end{proof}
\begin{nolabel}
	Most of the results proved in section \ref{sec:structure}  for the $N_W =
	N$ situation carry
over to the $N_K=N$ setting with the following modifications. 
\begin{itemize}
	\item replace $\partial_Z=(\partial_z, \partial_{\theta^1}, \dots,
		\partial_{\theta^N})$ by $D_Z = (\partial_z, D^1_Z, \dots, D^N_Z)$, 
	\item replace $Z-W = (z-w, \theta^i - \zeta^i)$ by $Z-W = \left( z - w - \sum_{i =
		1}^N \theta^i \zeta^i, \theta^j - \zeta^j \right)$, 
	\item replace $(Z-W)^{j|J} = (z-w)^j \prod_{i \in J} (\theta^i - \zeta^i)$ by \[
		(Z-W)^{j|J} = \left( z-w-\sum_{i = 1}^N \theta^i \zeta^i \right)^j \prod_{i
		\in J} (\theta^i - \zeta^i),\]
	\item replace the commutative associative ``translation'' superalgebra $\mathbb{C}[T,
		S^i]$ by the non-commutative associative ``translation'' superalgebra $\cH$
		generated by the set $\nabla=(T, S^1, \dots, S^N)$, where $T$ is an even
		generator and $S^i$ are odd generators, subject to the relations:
		\begin{equation}
			[T, S^i] = 0, \quad [S^i, S^j]=2 \delta_{ij} T;
			\label{eq:alpedo_generan}
		\end{equation}
	\item replace the commutative associative ``parameter'' superalgebra
		$\mathbb{C}[\lambda, \chi^i]$ by the non-commutative associative ``parameter''
		superalgebra $\cL$, generated by the set $\Lambda = (\lambda, \chi^1, \dots,
		\chi^N)$, where $\lambda$ is an even generator and $\chi^i$ are odd generators,
		subject to the relations:
		\begin{equation}
			[\lambda, \chi^i] = 0, \quad [\chi^i, \chi^j] = - 2 \delta_{ij}
			\lambda;
			\label{eq:alpedo_generan.2}
		\end{equation}
		Note that we have an isomorphism $\cH \rightarrow \cL$ given by $\nabla \mapsto
		- \Lambda$. 
\end{itemize} \label{no:whattochange}
\end{nolabel}
\begin{lem}
	The formal $\delta$-function satisfies the properties (1)-(7) of \ref{no:delta}
	after replacing $\partial_W$ by $D_W$ and writing $\Lambda + D_W = (\lambda +
	\partial_w, \chi^i + D^i)$. 
	\label{lem:k.delta}
\end{lem}
\begin{proof}
	(1) is clear from Lemma \ref{lem:k.deriv}. In order to prove (2) we use
	Lemma \ref{lem:k.deriv} to write:
	\begin{multline*}
		(Z-W)^{j|J} D^{(n|I)}_W \delta(Z,W) = \sigma(I)
		  \sigma(J, N\setminus I) \times \\
		\times (i_{z,w} - i_{w,z}) (Z-W)^{-1-n+j| N \setminus (I\setminus
		J)}.
	\end{multline*}
	Applying Lemma \ref{lem:k.deriv}  to 
		$D^{(n-j|I\setminus J)}_W \delta(Z,W)$
		the result follows from the following property of $\sigma$, which follows from
		(\ref{eq:not.2b}):
	\begin{equation*}
		 \sigma(J, N\setminus I) \sigma(I\setminus
		J) = \sigma(I)\sigma(I\setminus J, J).
	\end{equation*}
	Properties (3)-(7) are proved as in 
	\ref{no:delta}.
\end{proof}
\begin{lem}
	$D^i_Z (Z-W)^{j|J} = \sigma(e_i, J\setminus e_i) (Z-W)^{j|J\setminus e_i}
	+ j \sigma(e_i, J) (Z-W)^{j-1|J\cup e_i}$.
	\label{lem:k.deriv.2}
\end{lem}
\begin{proof}
	We prove the lemma by direct computation when $j \geq 0$:
	\begin{multline*}
			D^i_Z (Z-W)^{j|J} = (\partial_{\theta^i} + \theta^i
			\partial_z) \left( z - w - \sum \theta^i \zeta^i
			\right)^j (\theta - \zeta)^J \\
			= - j \zeta^i (Z-W)^{j-1 | J} + \sigma(e_i, J\setminus
			e_i) (Z-W)^{j|J\setminus e_i} +  \theta^i j
			(Z-W)^{j-1|J} \\
			= \sigma(e_i, J\setminus e_i) (Z-W)^{j|J\setminus e_i} +
			j \sigma(e_i, J) (Z-W)^{j-1|J \cup e_i}.
	\end{multline*}
	When $j = -1$ we have
	\begin{gather*}
		D^i_Z (Z-W)^{-1|J} = \left( \partial_{\theta^i} + \theta^i
		\partial_z \right) \sum_{k \geq 0} \frac{\left( \sum_{i} \theta^i
		\zeta^i\right)^k (\theta - \zeta)^J}{(z-w)^{k+1}} \\
		= \sum_{k \geq 0} k \zeta^i \frac{\left( \sum_{i} \theta^i
		\zeta^i \right)^{k-1}(\theta - \zeta)^J}{(z-w)^{k+1}}  +
		\sigma(e_i, J \setminus e_i) \times \\  \times \sum_{k
		\geq 0} \frac{\left( \sum_{i} \theta^i \zeta^i
		\right)^k (\theta-\zeta)^{J\setminus e_i}}{(z-w)^{k+1}}  -
		\theta^i \sum_{k \geq 0} (k+1) \frac{\left( \sum_i \theta^i
		\zeta^i \right)^k (\theta-\zeta)^J}{(z-w)^{k+2}} \\
		= \sigma(e_i, J\setminus e_i) (Z-W)^{-1|J\setminus e_i} 
		 -\sigma(e_i, J) \sum_{k \geq 0} (k+1) \frac{\left( \sum_i \theta^i
		\zeta^i
		\right)^{k}(\theta - \zeta)^{J \cup e_i}}{(z-w)^{k+2}} \\
		= \sigma(e_i, J \setminus e_i) (Z-W)^{-1|J\setminus e_i} -
		\sigma(e_i, J) (Z-W)^{-2|J\cup e_i}.
	\end{gather*}
	The general case follows from these by noting that $D^{1|0}_Z=(D^i_Z)^2 =
	\partial_z$, hence
	\begin{equation*}
		(Z-W)^{-j-1|J} = \frac{1}{j!} (D^i_Z)^{2j} (Z-W)^{-1|J}.
	\end{equation*}
	Therefore we get for $j \geq 0$:
	\begin{gather*}
			D^i_Z (Z-W)^{-j-1|J} = \frac{1}{j!} (D^i_Z)^{2j+1}
			(Z-W)^{-1|J} \\ = \frac{1}{j!} (D^i_Z)^{2j} \left(
			\sigma(e_i, J\setminus e_i) (Z-W)^{-1|J} - 
			 -  \sigma(e_i,
			J) (Z-W)^{-2|J \cup e_i}
			\right) \\ = \sigma(e_i, J \setminus e_i) (Z-W)^{-1-j|J
			\setminus e_i}  - (j+1) \sigma(e_i, J) (Z-W)^{-j-2|J \cup
			e_i}.
	\end{gather*}
\end{proof}
The following decomposition lemma is now proved in the same manner as
Lemma \ref{lem:decomp}:
\begin{lem}
	Let $a(Z,W)$ be a local distribution in two variables. Then $a(Z,W)$ can be
	uniquely decomposed in the following finite sum:
	\begin{equation*}
		a(Z,W) = \sum_{(j|J): j \geq 0} 
		\left(  D^{(j|J)}_W \delta(Z,W) \right) c_{j|J} (W).
	\end{equation*}
	The coefficients $c_{j|J}$ are given by
	\begin{equation*}
		c_{j|J}(W) = \res_Z (Z-W)^{j|J} a(Z,W).
	\end{equation*}
	\label{lem:k.decomp}
\end{lem}
\begin{rem} Let $(Z-W)\Lambda = \left( z - w- \sum_{i-1}^N \theta^i \zeta^i \right) \lambda + \sum_{i=1}^N
(\theta^i - \zeta^i) \chi^i$. Note that $-\partial_w$ and $-D^i_W$ satisfy the same commutation
relations (\ref{eq:alpedo_generan.2}) as $\lambda$, $\chi^i$, therefore $- \partial_w$, $-D^i_W$ generate an
associative superalgebra isomorphic to $\cL$. This allows us to consider $\cL$ as a
module over $\cH$, by letting $D^i_W$ and $\partial_w$ act as the following
derivations of the superalgebra $\cL$:
\begin{equation}
	[D^i_W, \chi^j] = 2 \delta_{ij} \lambda, \quad [\partial_w, \chi^i] =
	[\partial_w, \lambda] = [D^i_W, \lambda] = 0.
	\label{eq:k.intro.8}
\end{equation}
\end{rem}
\begin{lem} 
	 \[D^i_Z \exp\left( (Z-W)\Lambda \right) = \chi^i \exp\left(
			(Z-W)\Lambda
			\right) = - [D^i_W, \exp\left(  (Z-W)\Lambda \right)].\]
	\label{lem:k.expderiv}
\end{lem}
\begin{proof}
	Since the exponent is a sum of non-commuting terms, the derivative of the
	exponential is not as obvious as in the $N_W=N$ case. Let
	$	A = \sum_{j = 1}^N (\theta^j - \zeta^j) \chi^j$.
	We have:
	\begin{equation*}
		\begin{aligned}
			\exp\left( (Z-W)\Lambda \right) &= \exp\left( \left(z-w-\sum_{j = 1}^N
			\theta^j \zeta^j \right) \lambda 
			\right) \exp(A), \\
			\partial_{\theta^i} A^k &= \sum_{j = 0}^{k-1} A^j \chi^i
			A^{k-j-1}.
		\end{aligned}
	\end{equation*}
	Since $[A, \chi^i] = - 2 (\theta^i - \zeta^i) \lambda$ we obtain
	\begin{equation*}
			\partial_{\theta^i} A^k = \sum_{j=0}^{k-1} \chi^i A^{k-1}
			- 2 j \lambda (\theta^i - \zeta^i) A^{k-2} = 
			 k \chi^i A^{k-1} - k(k-1) \lambda (\theta^i - \zeta^i)
			A^{k-2},
	\end{equation*}
	therefore
	\begin{equation}
		\partial_{\theta^i} \exp(A) = \left( \chi^i - \lambda(\theta^i - \zeta^i)
		\right) \exp(A)
		\label{eq:k.exp.proof.4}
	\end{equation}
	from which the first equality of the lemma follows easily. The proof of the second
	equality of the lemma is similar. Note that from (\ref{eq:k.intro.8}) we have:
	\begin{equation*}
		[D^i_W, A] = - \chi^i - 2 \lambda (\theta^i - \zeta^i),
	\end{equation*}
	from where it follows as in (\ref{eq:k.exp.proof.4}) that
	\begin{equation*}
		[D^i_W, \exp(A)] = - (\chi^i + \lambda (\theta^i - \zeta^i))
		\exp(A),
	\end{equation*}
	and the lemma follows by a straightforward computation.
\end{proof}
\begin{nolabel}
	Now we are in position to define the \emph{formal Fourier transform} and 
	$N_K=N$ \emph{SUSY Lie conformal algebras} as we did in
	\ref{no:fourier}. We put
	\begin{equation*}
		\cF^\Lambda_{Z,W} a(Z,W) = \res_Z \exp\left( (Z-W)\Lambda \right)
		a(Z,W),
	\end{equation*}
	which formally looks exactly like (\ref{eq:fourier.1}) but in this expression
	the variables $\chi^i$ in $\Lambda = (\lambda, \chi^1, \dots, \chi^N)$ do
	not commute, but rather satisfy (\ref{eq:alpedo_generan.2}), and $(Z-W)$ is
	given by 
	  (\ref{eq:k.intro.2}) instead of  (\ref{eq:w.intro.recontra}). Using this formal
	Fourier transform, we define the $\Lambda$-bracket of two formal
	distributions $a(W)$ and $b(W)$ as in (\ref{eq:lie.2}).
	The $N_K=N$ version of Proposition \ref{prop:fourier}, on the properties
	of the Fourier transform, is proved in the same way as in the $N_W=N$ case with
	the aid of
	Lemma \ref{lem:k.expderiv}. There is only one subtlety involved in
	stating and proving (3) of Proposition \ref{prop:fourier}, and
	consequently,
	(3) in Proposition \ref{prop:lie}. Since the exponentials involved in this case
	do not commute, the argument in \ref{eq:fourierproof.6} is no
	longer valid.  
		We define 
	\begin{equation*}
		\cF^{\Lambda + \Gamma}_{X, W} = \cF^\Psi_{X,W}|_{\Psi = \Lambda + \Gamma},
	\end{equation*}
	where the Fourier transform on the RHS is computed as follows. First compute
	$\cF^\Psi_{X,W}$, and then replace $\Psi$ by $\Lambda + \Gamma = (\lambda + \gamma,
	\chi^1 + \eta^1, \dots, \chi^N + \eta^N)$. Here $\Psi$ denotes another
	set of indeterminates
	$\Psi = (\psi,
	\upsilon^1, \dots, \upsilon^N)$, where $\psi$ is an even indeterminate and $\upsilon^i$ are
	odd indeterminates, subject to the relations:
	\begin{equation*}
		\begin{gathered}
			{[\psi}, \upsilon^i] = 0, \quad  [\upsilon^i,
			\upsilon^j] = -2 \delta_{ij}\psi.
	\end{gathered}
	\end{equation*}
	We need to check that
	$\cF^\Lambda_{Z,W} \cF^\Gamma_{X,W} = (-1)^N
	\cF^{\Lambda + \Gamma}_{X,W} \cF^\Lambda_{Z,X}$, or, equivalently
	\begin{equation}
		\exp ( (Z-W) \Lambda) \exp ( (X-W) \Gamma) = \exp( (X-W) \Psi) \exp(
		(Z-W)\Lambda ) |_{\Psi = \Lambda + \Gamma}.
		\label{eq:invento.1}
	\end{equation}
	In order to compare both sides in (\ref{eq:invento.1}), we assume that 
	$\Psi$ and $\Lambda$ commute :
	\begin{equation}
		[\lambda, \psi] = [\lambda, \upsilon^i] = 0, \quad [\chi^i,
		\psi] = [\chi^i,
		\upsilon^j] = 0.
		\label{eq:unamasynojodemosmas}
	\end{equation}
	 Note that the
	RHS of (\ref{eq:invento.1}) can be also computed as 
	\begin{equation}
		\exp( (X-W)(\Lambda + \Gamma) ) \exp ( (Z-W) \Lambda),
		\label{eq:invento.2}
	\end{equation}
	where we have to use the following commutation relations between $\Lambda$
	and $\Gamma$:
	\begin{equation} [\eta^i, \chi^j]
	= 2 \lambda \delta_{i,j}, \quad [\gamma, \chi^i] = [\gamma, \lambda] = [\lambda,
	\eta^i] = 0, \label{eq:etacommutes_with_lambda}
	\end{equation}
	which follow from (\ref{eq:unamasynojodemosmas}) after replacing $\Psi$ by
	$\Lambda +
	\Gamma$.

	In order to check (\ref{eq:invento.1}), recall that, 
	given two operators $A, B$, such that their commutator $[A,B]=C$
	commutes with both $A$ and $B$, we have 
	\begin{equation} 
		e^A e^B = e^C e^B e^A.
		\label{eq:how.exp.commute}
	\end{equation}
	Let $Z=(z,\theta^i)$, $W=(w,\zeta^i)$ and $X=(x,\pi^i)$. Now we expand:
	\begin{multline}
		\exp\left( (Z-W)\Lambda \right) \exp \left( (X-W) \Gamma \right)
		= \exp\left( (z - w - \sum \theta^i \zeta^i) \lambda \right)
		\times \\ \times \exp\left(  \sum (\theta^i - \zeta^i) \chi^i  \right) 
		 \exp \left(
		(x-w-\sum \pi^i \zeta^i) \gamma \right) \exp\left(  \sum (\pi^i -
		\zeta^i) \eta^i
		\right).
		\label{eq:k.jacobi.subtle.1}
	\end{multline}
	Note also that we have
	\begin{equation}
		\begin{split}
			\exp\left( \sum(\theta^i - \zeta^i)\chi^i \right) =
			\prod \exp\left(  (\theta^i - \zeta^i) \chi^i \right) =
			\prod (1 + (\theta^i - \zeta^i) \chi^i) \\ = \prod \left(
			(1+\theta^i \chi^i)(1-\zeta^i \chi^i)(1 + \theta^i \zeta^i
			\lambda) \right) = 
			\prod e^{\theta^i \chi^i} e^{-\zeta^i \chi^i} e^{\theta^i
			\zeta^i \lambda} \\
			= \exp\left( \sum \theta^i \chi^i  \right) \exp\left(
			 - \sum \zeta^i \chi^i \right) \exp\left( \sum \theta^i
			 \zeta^i \lambda
			 \right),
		\end{split}
		\label{eq:k.jacobi.subtle.2}
	\end{equation}
	therefore (\ref{eq:k.jacobi.subtle.1}) reads:
	\begin{multline}
		\exp\left( (Z-W)\Lambda \right) \exp \left( (X-W) \Gamma \right)
		= \exp\left( (z - w) \lambda \right)
		\times \\ \times \exp\left(  \sum \theta^i \chi^i \right) \exp
		\left( - \sum \zeta^i \chi^i  \right) 
		 \exp \left(
		(x-w) \gamma \right) \exp\left(  \sum \pi^i \eta^i \right) \exp
		\left( -
		\sum \zeta^i  \eta^i
		\right).
		\label{eq:k.jacobi.subtle.3}
	\end{multline}
	Commuting the exponentials using (\ref{eq:how.exp.commute}) and
	(\ref{eq:etacommutes_with_lambda}),  (\ref{eq:k.jacobi.subtle.3}) can be
	expressed as:
	\begin{multline}
		\exp\left( (z-w)\lambda + (x-w)\gamma \right) \exp\left(  - \sum
		\zeta^i \chi^i\right) \exp \left( \sum \pi^i \eta^i \right)
		\times \\ \times
		\exp\left(- \sum \zeta^i \eta^i\right)  \exp\left(\sum \theta^i \chi^i
		\right) \exp \left(- 2 \sum \theta^i \pi^i \lambda \right) .
		\label{eq:k.jacobi.subtle.4}
	\end{multline}
	Multiplying and dividing by $\exp(\sum \pi^i \chi^i)$ and using
	(\ref{eq:k.jacobi.subtle.2}) we can express (\ref{eq:k.jacobi.subtle.4})
	as 
	\begin{multline*}
		\exp\left( (z-w)\lambda + (x-w)\gamma \right) \exp\left(  - \sum
		\zeta^i \chi^i\right) \exp \left( \sum \pi^i \eta^i \right)
		\times \\ \times
		\exp\left(- \sum \zeta^i \eta^i\right) \exp\left( \sum \pi^i
		\chi^i
		\right)  \exp\left(\sum (\theta^i
		- \pi^i) \chi^i
		\right)  \exp \left( - \sum \theta^i \pi^i \lambda \right) .
	\end{multline*}
	Combining again the exponentials it is easy to express this as
	\begin{multline*}
		\exp\left( (z-w)\lambda + (x-w)\gamma \right) \exp\left(
		\sum(\pi^i - \zeta^i)(\chi^i + \eta^i)
		\right) \times \\ \times\exp \left( -\sum \pi^i \zeta^i (\lambda + \gamma) \right)
		 \exp\left(\sum (\theta^i
		- \pi^i) \chi^i
		\right)  \exp \left( - \sum \theta^i \pi^i \lambda \right),
	\end{multline*}
	which is equal to (\ref{eq:invento.2}). 
	\label{no:k.conformal.1}
\end{nolabel}
\begin{nolabel}
	The definition of $(j|J)$-th products for $j \geq 0$ and the definition of an
	$N_K=N$ formal distribution Lie superalgebra generalizes in a
	straightforward way the corresponding definitions in \ref{no:lie}.
	The $N_K=N$ version	
	 of Proposition 
	\ref{prop:lie}, on the properties of the $\Lambda$-bracket, is now proved
	in the same way as in the
	$N_W=N$ case, with the aid of (\ref{eq:invento.1}).
	\label{no:jajajajaja.ja}
\end{nolabel}
\begin{defn} An $N_K = N$ SUSY Lie conformal algebra is a $\mathbb{Z}/2 \mathbb{Z}$-graded
	$\cH$-module $\cR$, endowed with a parity $N \mod 2$ $\mathbb{C}$-bilinear map $\cR \otimes_\mathbb{C}
	\cR \rightarrow \cL \otimes_{\mathbb{C}} \cR$ denoted (as before, we omit
	the symbol $\otimes$ in the $\Lambda$-bracket)
	\begin{equation*}
		a \otimes b \mapsto [a_\Lambda b] = \sum_{\stackrel{j \geq 0, J}{\text{finite}}}
		(-1)^{JN} \Lambda^{(j|J)} a_{(j|J)}b,
	\end{equation*}
	where $a_{(j|J)} b \in \cR$. This data should satisfy the following axioms:
	\begin{enumerate}
		\item sesquilinearity (this is an equality in $\cL \otimes \cR$):
			\begin{equation}
				[S^i a_\Lambda b] = - (-1)^N \chi^i [a_\Lambda
				b],
				\qquad [a_\Lambda S^i b] = (-1)^{a+N} \left(S^i
				+ \chi^i
				\right) [a_\Lambda b],
				\label{eq:k.sesqui.1}
			\end{equation}
			where in the RHS of the second equation, to obtain an element of $\cL
			\otimes \cR$, we first compute the $\Lambda$-bracket, and then we
			commute $S^i$ to the right using the relations $[S^i, \chi^j] = 2 \delta_{ij}
			\lambda$ (cf. (\ref{eq:k.intro.8})).
		\item skew-symmetry (this is an equality in $\cL \otimes \cR$):
			\begin{equation}
				[a_\Lambda b] = - (-1)^{a b + N} [b_{-\Lambda -
				\nabla} a],
				\label{eq:k.skew.1}
			\end{equation}
			where the commutator on the right hand side is computed as
			follows: first compute $[b_\Gamma a] = \sum_{j \geq 0 , J} \Gamma^{j|J}
			c_{j|J} \in \cL' \otimes \cR$, where $ \cL'$ is another copy of $\cL$
			generated by the set $\Gamma =
			(\gamma , \eta^1, \dots, \eta^N)$, where $\gamma$ is an even generator,
			$\eta^i$ are odd generators, subject to the relations
			 \[ [\gamma, \eta^i] = 0, \quad [\eta^i, \eta^j]
			= -2 \delta_{ij} \gamma.\]  Then replace $\Gamma$ by  $ - \nabla -
			\Lambda = (- T -
			\lambda, - S^1- \chi^1, \dots, - S^N - \chi^N)$
			 and apply $T$ and $S^i$ to $c_{j|J} \in \cR$.
		\item Jacobi identity (this is an equality in $\cL \otimes \cL' \otimes \cR$):
			\begin{equation*}
				[a_\Lambda [b_\Gamma c]] = (-1)^{aN + N} \left[
				[ a_\Lambda b]_{\Gamma + \Lambda} c \right] +
				(-1)^{(a+N)(b+N)} [b_\Gamma [a_\Lambda c]],
			\end{equation*}
			where $[ [a_\Lambda b]_{\Lambda + \Gamma} c]$ is computed as follows,
			first compute $[ [a_\Lambda b]_\Psi c] \in \cL \otimes \cL''
			\otimes \cR$, where $\cL''$ is another copy of $\cL$ generated by the
			set $\Psi = (\psi, \upsilon^1, \dots, \upsilon^N)$, where $\psi$ is an
			even generator, $\upsilon^i$ are odd generators, subject to the
			relations \[ [\psi, \upsilon^i] = 0, \quad [\upsilon^i, \upsilon^j] =
			-2 \delta_{ij} \psi.\]
			Then replace $\Psi$ by $\Lambda + \Gamma = (\lambda +
			\gamma, \theta^1 + \eta^1, \dots, \theta^N + \eta^N)$ to obtain an
			element of $\cL \otimes \cL' \otimes \cR$. 
	\end{enumerate}
	\label{defn:k.conformal.1}
\end{defn}
It follows that given $N_K=N$ formal distribution Lie superalgebra $(\fg, \cR)$,
the space $\cR$ is an $N_K=N$ SUSY Lie conformal algebra with $T= \partial_w$ and $S^i
= D^i_W$.
\begin{rem} We want to give an explanation for the commutation relations  $[S^i,
	\chi^j] = 2 \delta_{ij} \lambda$ 
	 appearing in sesquilinearity. 
	For this, we give an abstract descriptions of the axioms of a Lie conformal
	algebra as follows. Let $\cH$ be a co-commutative Hopf superalgebra with
	commultiplication $\Delta: \cH \rightarrow \cH \otimes \cH$ and antipode
	$S$ 
	(note that this is the case in definition \ref{defn:k.conformal.1}). Let
	$\cR$ be a (left) $\cH$-module. The spaces $\cR \otimes \cR$ and $\cH
	\otimes \cR$ are canonically $\cH$ modules with $h \mapsto \Delta h$ and
	we consider $\cH$ as a $\cH$-module with the adjoint action. An
	\emph{$\cH$-Lie conformal algebra} structure in $\cR$
	is a linear map
	$\phi: \cR \otimes \cR \rightarrow \cH \otimes \cR$ satisfying the
	following axioms (see \cite{dandrea:seudo}):
	\begin{itemize}
		\item $\phi$ is an homomorphism of $\cH$-modules,
			namely, the following diagram is commutative for any $h
			\in \cH$:
			\begin{equation*}
				\xymatrix{ \cR \otimes \cR
				\ar[r]^{\phi} \ar[d]_{\Delta h} &  \cH
				\otimes \cR
				\ar[d]^{\Delta h} \\ \cR \otimes \cR
				\ar[r]_{
				\phi}& \cH \otimes
				\cR}
        		\end{equation*}
		\item (Sesquilinearity) Let $L_h$ be the operator of left multiplication by
			$h$ in $H$. The following diagram is commutative:
			\begin{equation*}
				\xymatrix{ \cR \otimes \cR \ar[r]^{\phi}
				\ar[d]_{h \otimes 1} & \cH \otimes \cR
				\ar[d]^{L_h \otimes 1} \\ \cR \otimes \cR
				\ar[r]_{ \phi} & \cH \otimes \cR}
			\end{equation*}
		\item (Skew-symmetry) Let $A$ and $B$ be two $\cH$-modules. Let
			$\sigma_{12}$ be the permutation isomorphism $A \otimes B \simeq B
			\otimes A$. Let $\mu : \cH \otimes \cR \rightarrow \cR$
			be the natural multiplication coming from the
			$\cH$-module structure in $\cR$. The following diagram is
			commutative:
			\begin{equation*}
				\xymatrix{ \cR \otimes \cR \ar[r]^{\phi}
				\ar[d]_{\sigma_{12}} & \cH \otimes \cR
				\ar[d]^{(S\otimes \mu) \circ (\Delta \otimes 1)}
				\\ \cR \otimes \cR \ar[r]_{- \phi} & \cH
				\otimes \cR}
			\end{equation*}
		\item (Jacobi identity) Define three morphisms
			$\cR^{\otimes 3}
			\rightarrow \cH^{\otimes 2} \otimes \cR$ corresponding to
			the three terms in the Jacobi identity. First, let
			$\mu_{1\{23\}}$ be the composition
			\begin{equation*}
					\cR^{\otimes 3} \xrightarrow{1 \otimes
					\phi} \cR \otimes
					\cH \otimes \cR \xrightarrow{\sigma_{12}
					(1 \otimes \phi) \sigma_{12}}\cH^{\otimes
					2} \otimes \cR.
			\end{equation*}
			Similarly, we define $\mu_{2\{13\}}$ to be the
			composition:
			\begin{equation*}
				\cR^{\otimes 3} \xrightarrow{\sigma_{12}
				(1\otimes \phi) \sigma_{12}} \cH \otimes
				\cR^{\otimes 2} \xrightarrow{(1 \otimes \phi)}
				\cH^{\otimes 2} \otimes \cR
			\end{equation*}
			Finally, let $\nu : \cH \otimes \cH \rightarrow \cH$ be
			the multiplication map. We define $\mu_{\{12\}3}$ to be
			the composition:
			\begin{equation*}
				\cR^{\otimes 3} \xrightarrow{\phi \otimes 1} \cH
				\otimes \cR^{\otimes 2} \xrightarrow{1 \otimes
				\phi} \cH^{\otimes 2} \otimes \cR
				\xrightarrow{(\nu \otimes 1 \otimes 1)(1 \otimes \Delta \otimes 1)}
				\cH^{\otimes 2} \otimes \cR 
			\end{equation*}
			The Jacobi identity is the following
			axiom:
			\begin{equation*}
				\mu_{1\{23\}} = \mu_{\{12\}3} + \mu_{2\{13\}}
			\end{equation*}
	\end{itemize}
	In the $N_K=N$ case, identifying:
	\begin{equation*}
		S^i \mapsto -\chi^i, \qquad T \mapsto - \lambda, \qquad \gamma
		\mapsto \lambda, \qquad \eta^i \mapsto \chi^i,
	\end{equation*}
	and, as in \ref{lem:lie.from.pseudoaa}, changing the parity of $\cR$ if $N$
	is odd and defining 
	\begin{equation*}
		\phi(a \otimes b) = (-1)^{aN+N} [a_\Lambda b],
	\end{equation*}
	it is straightforward to check that the axioms of an $N_K=N$ SUSY Lie
	conformal algebra,
	as in Definition \ref{defn:k.conformal.1}, get transformed into the axioms
	of
	 an $\cH$-Lie conformal algebra. \label{rem:remarketacommute}
\end{rem}
\begin{nolabel}
Lemmas \ref{lem:Lie_R_1} and \ref{lem:affinization_1} hold in the $N_K=N$ setting replacing
$\partial_W$ with $D_W$ in (\ref{eq:lambda_bracket_affine}). This allows us to
construct a Lie
superalgebra of degree $N \mod 2$ $L(R) = \tilde\cR/ \tilde\nabla \tilde\cR$, and
the corresponding Lie superalgebra $\lie(\cR)$, from any $N_K=N$ SUSY Lie conformal
algebra $\cR$. In order to prove the $N_K=N$ versions of Propositions
\ref{prop:distributions_are_local} and \ref{prop:formal_distrib_from_conf}, we note
that for
$J=(j_1, \dots, j_k)$, we have
\begin{equation}
	\begin{aligned}
		D_W^{j|J} &= \partial_w^j (\partial_{\zeta^{j_1}} + \zeta^{j_1} \partial_w)
		\dots (\partial_{\zeta^{j_k}} + \zeta^{j_k} \partial_w) \\
		&= \sum_{K \subset J} \sigma(K, J\setminus K) \zeta^K \partial_W^{j + \sharp K|
		J \setminus K}.
	\end{aligned}
	\label{eq:expansion_of_D}
\end{equation}
Let now $a_{<n|I>} = a \otimes W^{n|I} \in L(\cR)$ for each $a \in \cR$. 
Using (\ref{eq:expansion_of_D}) and (\ref{eq:lambda_bracket_affine}) with $f = W^{n|I}$, $g =
W^{k|K}$ and letting $\Lambda = 0$, we compute the Lie bracket (of parity $N \mod 2$) in
$L(\cR)$:
\begin{multline}
	\{a_{<n|I>}, b_{<k|K>}\} = \sum_{j \geq 0, J} (-1)^{a J + b(I - J) +
	\frac{(J \cap I)(J \cap I - 1)}{2} + \frac{J(J-1)}{2}}
	\frac{(n)_{n-j-\sharp(J\setminus I)}}{j!} \times \\ \times \sigma(J\setminus I, J \cap I) \sigma(J \cap I, I \setminus J) \sigma(J
	\setminus I, I \setminus J) \sigma(I \triangle J, K) \left( a_{(j|J)}b \right)_{<n+k -j
	-\sharp (J\setminus I)| K \cup (I\triangle J)>},
	\label{eq:horrible.1.bb}
\end{multline}
where $(n)_k = n (n-1) \dots (n-k+1)$.
Defining $a_{(n|I)}$ as the image of $(-1)^{aI} \sigma(I) a_{<n|I>}$ in
$\lie(\cR)$ and using (\ref{eq:horrible.1.bb}) and Lemma \ref{lem:lie.from.pseudoaa} we
compute the Lie bracket in $\lie(\cR)$:
\begin{multline}
	{[}a_{(n|I)}, b_{(k|K)}] = (-1)^{(a+N-I)(N-K)}\sum_{j \geq 0, K} (-1)^{J(N-I) + IN +\frac{(J
	\cap I)(J \cap I - 1)}{2} + \frac{J(J+1)}{2}}  \times \\ \times
	\frac{(n)_{n-j-\sharp(J\setminus
	I)}}{j!} \sigma(I \triangle J, N \setminus (K \cup I \triangle J)) 
	 \sigma(I)
	 \sigma(J \setminus I, J \cap I) \sigma(J \cap I, I
	\setminus J) \times \\ \times \sigma(J \setminus I, I \setminus J) \left( a_{(j|J)}b \right)_{(n + k - j
	- \sharp (J\setminus I)| K \cup (I\triangle J))}.
	\label{eq:horripilante.1aa}
\end{multline}
Substituting  (\ref{eq:delta_expanded}) in (\ref{eq:expansion_of_D}) we find:
\begin{multline}
	D^{(j|J)}_W \delta(Z,W) = \sum_{n \in \mathbb{Z}, I} (-1)^{\frac{(J\cap I)(J\cap I
	-1)}{2} + \frac{J(J+1)}{2} + I + (N-I)(J-I)} \times \\ \times
	\frac{(n)_{n-j-\sharp(J\setminus I)}}{j!} \sigma(J \setminus I, J \cap I) \sigma(J \cap I) \times \\
	\times \sigma(I \setminus J, N \setminus I) \sigma(J \setminus I, I \setminus J)
	Z^{-1-n|N\setminus I} W^{n-j-\sharp(J \setminus I)| I \triangle J}.
	\label{eq:calculo_horrendo}
\end{multline}
For each $a \in \cR$ define the following $\lie(\cR)$-valued formal distribution:
\begin{equation}
	a(Z) = \sum_{n \in \mathbb{Z}, I} Z^{-1-n|N\setminus I} a_{(n|I)}.
	\label{eq:formal_k_distr}
\end{equation}
Using (\ref{eq:horripilante.1aa}) and (\ref{eq:calculo_horrendo}) we obtain the
$N_K=N$ analog of Proposition \ref{prop:distributions_are_local}:
\begin{equation*}
	{[a(Z)},b(W)] = \sum_{j \geq 0, J} \left( D^{(j|J)}_W \delta(Z,W) \right) \left(
	a_{(j|J)}b
	\right) (W).
\end{equation*}
To prove the $N_K=N$ analog of Proposition \ref{prop:formal_distrib_from_conf} we
need the following identity which is straightforward to check:
	\begin{equation*}
			(S^i a)_{(j|J)} = 
			\begin{cases}
				\sigma(e_i, N\setminus J) a_{(j|J\setminus e_i)}
				& \quad \text{for } e_i \in J \\
				-j \sigma(e_i, N\setminus (J\cup e_i))
				a_{(j-1|J\cup e_i)} & \quad \text {for } e_i \notin
				J.
			\end{cases}
	\end{equation*}	
\end{nolabel}
\begin{nolabel}
	As in \ref{no:regular_lie_R}, the $N_K=N$ formal distribution Lie
	superalgebra $\lie(\cR)$ carries an even derivation $T$ and $N$ odd
	derivations $S^i$ $(i= 1, \dots, N)$, given by:
	\begin{equation*}
		\begin{aligned}
			T \bigl(a_{(j|J)} \bigr) &= - j a_{(j-1|J)}, \\
		S^i \bigl( a_{(j,J)} \bigr) &=\begin{cases}
			\sigma(N\setminus J, e_i) a_{(j, J \setminus e_i)} &
			\quad e_i \in J \\
			j \sigma(N\setminus (J \cup e_i), e_i) a_{(j-1|J \cup
			e_i)} & \quad e_i \notin J.
		\end{cases} 
	\end{aligned}
	\end{equation*}
	It follows easily that the formal distributions (\ref{eq:formal_k_distr})
	satisfy:
	\begin{equation*}
		Ta(Z) = \partial_z a(Z), \qquad S^ia(Z) = (\partial_{\theta^i} -
		\theta^i \partial_z) a(Z), \quad i=1, \dots, N,
		\label{eq:formality_of_k}
	\end{equation*}
	and therefore $(\lie(\cR), \cR)$ is a \emph{regular} $N_K=N$ formal
	distribution Lie superalgebra.
\end{nolabel}
	We define the normally ordered product of fields 
	by the same formula (\ref{eq:normal_defin}) as in the $N_W=N$ case, and all the other products by using the
	derivations $D_W$ instead of $\partial_W$.  
	Lemma \ref{lem:normal_lem} is still valid in the $N_K=N$ setting (recall that $(Z-W)^{-1|N}$
	is the same in both situations). The $N_K=N$ version of
	Proposition \ref{prop:sesquilinearity} is:
\begin{prop}
	The following identities analogous to sesquilinearity for all pairs
	$(j|J)$ are true:
	\begin{equation}
		\begin{aligned}
			\left( D^i_W a(W) \right)_{(j|J)}b(W) &= -(-1)^J \left(
			\sigma(e_i, J) a(W)_{(j|J\setminus e_i)} b(W) + \right.
			\\ & \quad \left. + j \sigma(e_i, J) a(W)_{(j-1|J \cup
			e_i)} b(W) \right) \\
			D^i_W \left( a(W)_{(j|J)}b(W) \right) &= (-1)^{N-J}\left(
			\left( D^i_W a(W) \right)_{(j|J)} b(W) + \right. \\ &
			\quad \left. + (-1)^a a(W)_{(j|J)} \left( D^i_W b(W)
			\right) \right).
		\end{aligned}
		\label{eq:k.sesqui.prop.1}
	\end{equation}
	\label{prop:k.sesqui.prop}
\end{prop}
\begin{proof}
	According to Lemma \ref{lem:k.deriv.2} we have:
	\begin{multline}
		\res_Z i_{z,w} (Z-W)^{j|J} D^i_Z a(Z)b(W)= \\= - (-1)^{J} \res_Z
		\left( D^i_Z i_{z,w} (Z-W)^{j|J} \right) a(Z) b(W) =\\ = - (-1)^J
		 \res_Z \left( \sigma(e_i, J\setminus e_i) i_{z,w}(Z-W)^{j|J\setminus
		 e_i} + \right. \\ + \left. j \sigma(e_i, J) i_{z,w} (Z-W)^{j-1|J\cup e_i} \right)
		 a(Z) b(W).
		\label{eq:k.sesqui.prop.proof.1}
	\end{multline}
	Similarly we have:
	\begin{multline}
		- (-1)^{(a+1)b} \res_Z i_{w,z} (Z-W)^{j|J}b(W) \left( D^i_Z a(Z)
		\right) = \\ = (-1)^{ab + J} \res_Z \left( D^i_Z
		i_{w,z}(Z-W)^{j|J} \right) b(W)a(Z) = \\ = (-1)^{ab + J} \res_Z
		\left( \sigma(e_i, J\setminus e_i) i_{w,z}(Z-W)^{j|J\setminus
		e_i} + \right. \\ + \left. j \sigma(e_i, J) i_{w,z} (Z-W)^{j-1|J
		\cup e_i} \right) b(W)a(Z).
		\label{eq:k.sesqui.prop.proof.2}
	\end{multline}
	Adding (\ref{eq:k.sesqui.prop.proof.1}) and
	(\ref{eq:k.sesqui.prop.proof.2}) we obtain:
	\begin{multline*}
		\left( D^i_W a(W) \right)_{(j|J)} b(W) = - (-1)^J \left( \sigma(e_i, J
		\setminus e_i) a(W)_{(j|J\setminus e_i)}b(W) + \right. \\ +
		\left. j \sigma(e_i, J) a(W)_{(j-1|J\cup e_i)}b(W) \right).
	\end{multline*}
	The fact that $D^i_W$ is a derivation of all $(j|J)$-products is proved
	in the same way as in (\ref{eq:sesqui_proof.4}): 
	\begin{multline*}
		D^i_W \left( a(W)_{(j|J)}b(W) \right) =
		D^i_W \res_Z \left( i_{z,w} (Z-W)^{j|J} a(Z) b(W) -
		\right. \\ - \left. (-1)^{ab} i_{w,z} (Z-W)^{j|J} b(W)a(Z) \right)
		= \\ (-1)^N  \res_Z \left( \left( - \sigma(e_i,J\setminus e_i) i_{z,w}
		(Z-W)^{j|J \setminus e_i} - j \sigma(e_i, J) i_{z,w}
		(Z-W)^{j-1|J\cup e_i} \right)\times \right. \\ \left . \times  a(Z)b(W)  + (-1)^{J+a}
		i_{z,w} (Z-W)^{j|J} a(Z) D^i_W b(W) + \right. \\ \left.
		+ (-1)^{ab} \left( \sigma(e_i, J\setminus e_i) i_{w,z} (Z-W)^{j|J
		\setminus e_i} + j  \sigma(e_i, J) i_{w,z} (Z-W)^{j-1|J \cup
		e_i} \right) b(W) a(Z) - \right. \\ \left. - (-1)^{ab + J}
		i_{w,z} (Z-W)^{j|J} D^i_W b(W) a(Z)
		\right) = \\ = - (-1)^N \sigma(e_i, J\setminus e_i)
		a(W)_{(j|J\setminus e_i)} b(W) - (-1)^N j \sigma(e_i, J)
		a(W)_{(j-1|J\cup e_i)} b(W)+ \\ + (-1)^{N + J+a}
		a(W)_{(j|J)}D^i_W b(W) = \\ = (-1)^{N-J} \left( \left(
		D^i_W a(W)
		\right)_{(j|J)} b(W) + (-1)^{a} a(W)_{(j|J)} D^i_Wb(W)
		\right)
	\end{multline*}
\end{proof}
\begin{nolabel}
	Even though the general Proposition \ref{prop:wick.1} is no longer valid
	in the $N_K=N$ setting, we easily check that its proof works in the
	particular case $(j|J)=(-1|N)$. Therefore, the non-commutative Wick formula
	(\ref{eq:wick.proof.8}) is still valid in the $N_K=N$ case. 
	The $N_K=N$ version of 
	Dong's Lemma \ref{lem:dong} is proved as in the usual vertex algebra case.
\end{nolabel}
\begin{defn}
	An $N_K=N$ SUSY vertex algebra is the data
	consisting of a super
	vector space $V$, an even vector $\vac \in V$, $N$ odd endomorphisms
	$S^i$ and a parity preserving linear map $\ys$ from $V$ to $\End(V)$-valued
	superfields $a \mapsto \ys(a,Z)$, satisfying the following axioms:
\begin{itemize}
	\item vacuum axioms:
		\begin{equation*}
			\begin{aligned}
			\ys(a, Z) \vac = a + O(Z), \quad
			S^i \vac = 0, \qquad i = 1, \dots, N,
			\end{aligned}
		\end{equation*}
	\item translation invariance:
		\begin{equation*}
			{[} S^i, \ys(a,Z)] = \bar{D}^i_Z \ys(a,Z),
		\end{equation*}
		where $\bar{D}^i_Z = \partial_{\theta^i} - \theta^i \partial_z$,
	\item locality:
		\begin{equation*}
			(z-w)^n [\ys(a,Z), \ys(b,W)] = 0 \qquad \text{for some } n
			\in \mathbb{Z}_+.
		\end{equation*}
\end{itemize}
	\label{no:k.vertex.defin}
\end{defn}
\begin{nolabel}
	We define the $(j|J)$-products for a $N_K=N$ SUSY vertex algebra, as in the
	$N_W=N$ case,  by (\ref{eq:product.3}). 

	As in \ref{no:product_axioms} we see easily that the vacuum axioms may be
	formulated as (\ref{eq:product.4}) and translation invariance is
	equivalent to:
	\begin{equation}
		[S^i, a_{(j,J)}] =\begin{cases}
			\sigma(N\setminus J, e_i) a_{(j, J \setminus e_i)} &
			\quad e_i \in J \\
			j \sigma(N\setminus (J \cup e_i), e_i) a_{(j-1|J \cup
			e_i)} & \quad e_i \notin J 
		\end{cases}
		\label{eq:k.product.1}
	\end{equation}
	\label{no:k.product}
\end{nolabel}
\begin{nolabel}
	It follows easily from (\ref{eq:k.product.1}) and the vacuum axioms that 
	\begin{equation*}
		[S^i, S^j] = 2 \delta_{i,j} T, \qquad [S^i, T] = 0,
	\end{equation*}
	where $T$ is an even operator satisfying:
	\begin{equation*}
		[T, a_{(j|J)}] = - j a_{(j-1|J)} \qquad \forall a, (j|J),
	\end{equation*}
	or equivalently:
	\begin{equation*}
		[T, \ys(a,Z)] = \partial_z \ys(a,Z).
	\end{equation*}
	\label{no:k.t.defin}
\end{nolabel}
With these results we can prove the $N_K=N$ version of Theorem
\ref{thm:fields.1}. 
\begin{thm}
	Let $\cU$ be a vector superspace and $V$ a space of pairwise local $\End(\cU)$-valued
	fields such that $V$ contains the constant field $\Id$, it is invariant
	under the derivations $D^i_Z= \partial_{\theta^i}+ \theta^i \partial_z$
	and closed under
	all $(j|J)$-th products. Then $V$ is an $N_K=N$ SUSY vertex algebra with
	vacuum vector $\Id$, the translation operators are $S^i a(Z) = D^i_Z
	a(Z)$,
	the $(j|J)$ product is the one for
	distributions multiplied by $\sigma(J)$.
	\label{thm:k.fields.1}
\end{thm}
\begin{proof}
	The proof goes like the proof of \ref{thm:fields.1}. To check translation
	invariance we see that $D^i_Z 1 = 0$ and that
	\begin{multline*}
		\sigma(J) D^i_Z \left( a(Z)_{(j|J)}b(Z) \right) - (-1)^{a+N-J}
		a(Z)_{(j|J)}D^i_Z b(Z) = \\ = (-1)^{N-J}\sigma(J) \left( D^i_Z a(Z)
		\right)_{(j|J)} b(Z).
	\end{multline*}
	But in view of (\ref{eq:k.sesqui.prop.1}) this is:
	\begin{multline*}
		- (-1)^N \sigma(J) \left( \sigma(e_i, J)
		a(Z)_{(j|J\setminus e_i)}b(Z) +  j \sigma(e_i, J) a(W)_{(j-1|J
		\cup e_i)} b(Z) \right)= \\ = \sigma(N\setminus J, e_i)
		\sigma(J\setminus e_i)
		a(Z)_{(j|J\setminus e_i)}b(Z) + j \sigma(N\setminus (J \cup e_i),
		e_i) \sigma(J\cup e_i) a(Z)_{(j-1|J\cup
		e_i)} b(Z),
	\end{multline*}
	proving equation (\ref{eq:k.product.1}). 
	Locality is proved in the same way as in \ref{thm:fields.1}. 
\end{proof}
Lemma \ref{lem:constant_term} is still valid for $N_K=N$ SUSY vertex algebras. Its proof
parallels the proof for $N_W=N$ SUSY vertex algebras. Lemma
\ref{lem:differential}, on the existence and uniqueness of solutions to a system of
differential equations, is straightforward to generalize to the $N_K=N$ setting. The proof of Proposition
\ref{prop:structure_thm_comienzo} in this context is more subtle:
\begin{prop}
	Let $V$ be a $N_K=N$ SUSY vertex algebra. Then for every $a, b \in V$ we
	have:
	\begin{enumerate}
		\item $\ys(a,Z)\vac = \exp(Z \nabla) a$,
		\item $\exp(Z \nabla)\ys(a,W) \exp(-Z \nabla) = i_{w,z} \ys(a,
			W+Z)
			$,
		\item $\ys(a,Z)_{(j|J)}\ys(b,Z)\vac = \sigma(J)
			\ys(a_{(j|J)}b,Z)\vac$.
	\end{enumerate}
	where $\nabla = (T, S^1, \dots, S^N)$, $Z
	\nabla = zT + \sum \theta^i S^i$, and we define $W+Z = W - (-Z) =
	(z+w+ \sum \zeta^i \theta^i, \theta^j + \zeta^j)$\footnote{Note that $Z+W
	\neq W+Z$}.  
	\label{prop:k.structure.com}
\end{prop}
\begin{proof}
	As in the proof of Proposition \ref{prop:structure_thm_comienzo} we note
	that both sides of (1) and (3) are elements of $V[ [Z]]$, whereas both
	sides of (2) are elements of $\End(V)[ [W, W^{-1}] ] [ [Z] ]$. By
	evaluating at $Z = 0$ we get equalities in all three cases. Indeed (1)
	and (2) are trivial, and (3) follows from the $N_K=N$ version of
	Lemma \ref{lem:constant_term}. We need to show that both sides in
	each equation satisfy the same system of differential equations.

		(1) Similarly to the proof of Lemma
			\ref{lem:k.expderiv}, we expand:
			\begin{equation}
				\begin{split}
					\bar{D}^i_Z \exp(Z \nabla) = \left(
					\partial_{\theta^i} - \theta^i \partial_z 
					\right) \exp(z T) \sum_{k \geq 0}
					\frac{(\sum_i \theta^i S^i)^k}{k!} \\
					= (S^i + T \theta^i) \exp(Z \nabla) -
					\theta^i  T \exp(Z \nabla) 
					= S^i \exp(Z \nabla),
				\end{split}
				\label{eq:k.structure.proof.1}
			\end{equation}
			from where the RHS $X(Z)$ of (1) satisfies
			the system of differential equations:
			\begin{equation*}
				\bar{D}^i_Z X(Z) = S^i X(Z).
			\end{equation*}
			Similarly by translation invariance we have for the 
			LHS of (1):
			\begin{equation*}
				\bar{D}^i_Z \ys(a,Z) \vac = [S^i, \ys(a,Z)]\vac =
				S^i \ys(a,Z) \vac.
			\end{equation*}

			We also point out that a computation similar to
			(\ref{eq:k.structure.proof.1}) shows that 
			\begin{equation}
				\bar{D}^i_Z \exp(-Z \nabla) = - \exp(-Z\nabla)
				S^i,
				\label{eq:k.structure.proof.3.b}
			\end{equation}
			which is not entirely obvious since $S^i$ does not commute
			with the exponential.

		(2) By translation invariance we have:
			\begin{multline*}
				\bar{D}^i_Z \ys(a, W+Z) = \left(- \zeta^i
				\partial_{w+z+\sum \zeta^i \theta^i}
				 +  \partial_{\zeta^i +
				 \theta^i} - \theta^i
				 \partial_{w+z+\sum \zeta^i \theta^i}
				 \right) \ys(a,W+Z) = \\ = \bar{D}^i_{W+Z}
				 \ys(a,Z+W) = [S^i, \ys(a,Z+W)]. 
			\end{multline*}
			On the other hand, letting $Y(Z) = e^{Z\nabla} \ys(a,W)
			e^{-Z \nabla}$ we have (cf.
			(\ref{eq:k.structure.proof.3.b})):
			\begin{equation*}
					\bar{D}^i_Z Y(Z) = S^i 
					 Y(Z) - (-1)^a Y(Z) S^i 
					= [S^i, Y(Z)].
			\end{equation*}
			
		(3) For the RHS we have by translation invariance and
			the vacuum axioms:
			\begin{equation*}
				S^i \ys(a_{(j|J)}b,Z)\vac =
				[S^i,\ys(a_{(j|J)}b,Z)]\vac = \bar{D}^i_Z
				\ys(a_{(j|J)}b,Z)\vac .
			\end{equation*}
			To prove that the LHS satisfies the same
			differential equation, we proceed exactly in the same way as
			in the proof of Proposition
			\ref{prop:structure_thm_comienzo}. We only need the fact
			that $\bar{D}^i_Z$ is a derivation of all $(j|J)$-products.
			But $\partial_z=(D^i_Z)^2$ is a derivation since:
			\begin{multline*}
				\partial_z a(Z)_{(j|J)}b(Z) = (-1)^{N-J} D^i_Z \left( (D^i_Z
				a(Z))_{(j|J)}b(Z) + (-1)^a a_{(j|J)} D^i_Z b(Z)
				\right) = \\ (\partial_z a(Z))_{(j|J)} b(Z) +
				(-1)^{a+1} (D^i_Z a(Z))_{(j|J)} D^i_Zb(Z) + \\ + (-1)^a
				(D^i_Z a(Z))_{(j|J)} D^i_Z b(Z) +
				a(Z)_{(j|J)}\partial_z b(Z)
			\end{multline*}
			therefore $\bar{D}^i_Z = D^i_Z - 2 \theta^i \partial_z$ is
			a derivation of all $(j|J)$-products.
\end{proof}
The uniqueness Proposition \ref{prop:uniqueness.1} is still valid in the $N_K=N$
setting,
As its corollary, 
 we obtain an analogous version of Theorem \ref{thm:structure.2}, namely
\begin{thm}
	On an $N_K=N$ SUSY vertex algebra the following identities hold
	\begin{enumerate}
		\item $\ys(a_{(j|J)}b,Z) = \sigma(J)
			\ys(a,Z)_{(j|J)}\ys(b,Z)$ .
		\item $\ys(a_{(-1|N)}b,Z) = :\ys(a,Z)\ys(b,Z):$.
		\item $\ys(S^ia, Z) = D^i_Z \ys(a,Z)$.
		\item We have the following OPE formula:
			\begin{equation*}
				{[}\ys(a,Z), \ys(b,W)] = \sum_{(j|J): j \geq 0}
				\sigma(J)
				(D_W^{(j|J)} \delta(Z,W) )
				\ys(a_{(j|J)}b,W) 
			\end{equation*}
	\end{enumerate}
	\label{thm:k.structure.2}
\end{thm}
\begin{rem}
	Note that as a consequence of (3) we obtain
	\begin{equation*}
		[S^i,\ys(a,Z)] \neq \ys(S^i a, Z),
	\end{equation*}
	in contrast to the $N_W=N$ and, in particular, the ordinary vertex algebra case.
\end{rem}
\begin{cor}
	\begin{equation*}
		\begin{aligned}
			(S^i a)_{(j|J)} &= \sigma(e_i, N\setminus J)
			a_{(j|J\setminus e_i)} - j \sigma(e_i,N\setminus (J\cup
			e_i)) a_{(j-1|J\cup e_i)} \\
			&=\begin{cases}
				\sigma(e_i, N\setminus J) a_{(j|J\setminus e_i)}
				& \quad \text{for } e_i \in J \\
				-j \sigma(e_i, N\setminus (J\cup e_i))
				a_{(j-1|J\cup e_i)} & \quad \text {for } e_i \notin
				J,
			\end{cases}\\
			S^i\left( a_{(j|J)}b \right) &= (-1)^{N-J}\left( (S^i
			a)_{(j|J)}b + (-1)^a a_{(j|J)}S^i b.
			\right)
		\end{aligned}
	\end{equation*}
	\label{cor:k.cor.1}
\end{cor}
\begin{nolabel} In order to prove the $N_K=N$ version of the \emph{associativity}
	formulas 
	(\ref{eq:quasi-assoc.vertex.6}) and (\ref{eq:quasi-assoc.vertex.7}), we
	proceed as in \ref{no:quas-assoc.vertex}, by taking the generating series
	of \ref{thm:k.structure.2} (1) and using the following $N_K=N$ version of the
	\emph{Taylor expansion}. For a 
	formal power series $a(Z) \in \mathbb{C}[ [Z]]$ we have:
	\begin{equation}
		a(W+Z) = \sum_{(j|J): j \geq
		0}(-1)^{\frac{J(J-1)}{2}} \frac{W^{j|J}}{j!}
		D^{j|J}_Z a(Z) = e^{W D_Z} a(Z).
		\label{eq:k.no.3.3}
	\end{equation}
	Indeed, the usual Taylor expansion is:
	\begin{equation*}
		a(W+Z) = a\left( w + z + \sum \zeta^i \theta^i, \zeta^j +
		\theta^j \right) = \sum (-1)^{\frac{J(J-1)}{2}} \frac{\partial_W^{j|J}}{j!}
		a(W+Z)|_{W = 0}.
	\end{equation*}
	In this case:
	\begin{equation*}
		\begin{aligned}
			\partial_W^{1|0} a(W+Z)|_{W=0} &= D_{Z}^{1|0}a(Z), \\
			\partial_W^{0|i} a(W+Z)|_{W=0} &= (\theta^i \partial_z +
			\partial_{\theta^i}) a(Z) = D_Z^i a(Z),
		\end{aligned}
	\end{equation*}
	proving (\ref{eq:k.no.3.3}).
	
	Also, according to our prescription to add coordinates we see
	that
	\begin{equation*}
		(X-Z)-W = X- (W+Z) = X- (W - (-Z)),
	\end{equation*}
	and  note that 
	equation (\ref{eq:fields.proof.6}) is still valid in the $N_K=N$ setting.
	The proof in \ref{no:quas-assoc.vertex} generalizes now easily. Similarly,
	we obtain as a corollary,  
	the $N_K=N$ version of the Cousin property \ref{cor:cousin.w}.
	
	The proofs for \emph{skew-symmetry} in Theorem \ref{thm:skew-symmetry.w}
	and
	 quasi-commutativity for the normally ordered product as in
	 \ref{no:skew.cor} carry over verbatim to the $N_K=N$ case.
\end{nolabel} 
\begin{nolabel}
	Defining as the \emph{Fourier Transform} as in \ref{no:fourier_1variable}
	we obtain an analogous result to Theorem \ref{thm:conformal_vertex.1},
	namely an $N_K=N$ SUSY vertex algebra
	gives rise to an $N_K=N$ SUSY Lie conformal algebra. 

	The $N_K=N$ version of \emph{quasi-associativity} for the normally ordered
	 product is the same ans is proved in the same way as Theorem
	 \ref{thm:quasi-associativity1}. 
\end{nolabel}
As in the $N_W=N$ case, we have the following equivalent definition of $N_K=N$
SUSY vertex algebras:
\begin{defn}
	An $N_K=N$ SUSY vertex algebra is a tuple $(V, T, S^i, [\cdot_\Lambda\cdot], \vac, :
	:)$, $i = 1, \dots, N$,
	where
	\begin{itemize}
		\item $(V, T, S^i, [\cdot_\Lambda \cdot] )$ is an $N_K=N$ SUSY Lie conformal
			algebra,
		\item $(V, \vac,T, S^i, : :)$ is a unital quasicommutative quasiassociative
			differential superalgebra (i.e. $T$ is an even derivation
			of $: :$ and $S^i$ are odd derivations of $: :$),
		\item the $\Lambda$-bracket and the product $: :$ are related by the
			non-commutative Wick formula (\ref{eq:wick.proof.8}). 
	\end{itemize}
	\label{defn:bakalov.k}
\end{defn}
\begin{nolabel}
	The definition of an \emph{$N_K=N$ Poisson SUSY vertex algebra} is
	straightforward to generalize. Similarly,
	the $N_K=N$ version of Borcherds identity in Theorem
	\ref{thm:borcherds.identity.w} and the $N_K=N$ version of the
	commutator formula in  
	Proposition \ref{prop:commutator_w.1}, are the same with $\partial_W$
	replaced by $D_W$ and are proved in the same way as for
	$N_W=N$ SUSY vertex algebras. Following the argument in Remark \ref{rem:lie_bracket}
	and using (\ref{eq:expansion_of_D}), we obtain the formula for the 
	commutator of the Fourier coefficients of
	the $\End(V)$-valued fields of the $N_K=N$ SUSY vertex algebra: it is
	equal to the RHS of 
	(\ref{eq:horripilante.1aa}), multiplied by $\sigma(J)$.
\end{nolabel}
	The rest of section \ref{sub:existence} carries over to the $N_K=N$ case with minor
modifications, in particular we define tensor products of $N_K=N$ SUSY vertex
algebras as in \ref{no:tensor} and we have an \emph{existence theorem} as in
\ref{thm:existence} that we restate here:
\begin{thm}[Existence of $N_K=N$ SUSY vertex algebras] 
	Let $V$ be a vector space, $\vac \in V$ an even vector, $T$
	an even endomorphism of $V$ and $S^i$, $i = 1, \dots, N$ odd endomorphisms
	of $V$, satisfying $[S^i,S^j] = 2 \delta_{i,j} T$. Suppose moreover that
	$S^i \vac =0$. Let $\cF$ be a family
	of fields 
	$a^\alpha(Z) = \sum_{(j|J)} Z^{-1-j|N\setminus J} a^\alpha_{(j|J)}$,
	indexed by $\alpha \in A$, and such that 
	\begin{enumerate}
		\item $a^\alpha(Z)\vac|_{Z = 0} = a^\alpha \in V$,
		\item $[S^i,
			a^\alpha(Z)] = \bar{D}^i_Z
			a^\alpha(Z)$,
		\item all pairs $(a^\alpha(Z), a^\beta(Z))$ are local,
		\item the vectors 
				$a^{\alpha_s}_{(j_s|J_s)} \dots
				a^{\alpha_1}_{(j_1|J_1)} \vac$ span $V$.
	\end{enumerate}
	Then the formula
	\begin{multline*}
		\ys(a^{\alpha_s}_{(j_s|J_s)} \dots a^{\alpha_1}_{(j_1|J_1)} \vac,
		Z) = \\ =  \prod \sigma(J_i)
		a^{\alpha_s}(Z)_{(j_s|J_s)} \Bigl( \dots a^{\alpha_2}_{(j_2|J_2)}(Z)
		\bigl(a^{\alpha_1}_{(j_1|J_1)}(Z) \Id\bigr) \dots \Bigr)
	\end{multline*}
	defines a structure of an $N_K=N$ SUSY vertex algebra on $V$, with vacuum vector
	$\vac$, translation operators  $S^i$ and such that
		$\ys(a^\alpha, Z) = a^\alpha(Z)$, $\alpha \in A$.

	Such a structure is unique.
	\label{thm:k.existence}
\end{thm}
\begin{nolabel}
	The results in section \ref{sub:universal} generalize to this context without
	difficulty. In particular, we obtain 
	 the \emph{universal enveloping $N_K=N$ SUSY vertex
	algebra} of an $N_K=N$ SUSY Lie conformal algebra. 
	
	In the same way as in Theorem
	\ref{thm:lie_fourier.vertex.1}, we obtain:
\end{nolabel}
\begin{thm}
	Let $V$ be an $N_K=N$ SUSY vertex algebra. Let $\cA = \mathbb{C}[X,
	X^{-1}]$, 
	define $L(V)$ to be the
	quotient of $\tilde{V}=\cA \otimes_\mathbb{C} V$ by the span of vectors of the
	form:
	\begin{equation*}
		S^i a \otimes f(X) + (-1)^{a} a \otimes D^i_X f(X).
	\end{equation*}
	and let $L'(V)$ be its completion with respect to the natural topology in
	$\cA$. Then $L(V)$ (resp. $L'(V)$) carries a natural Lie superalgebra
	of degree $N \mod 2$ structure. Let $\lie(V)$ (resp. $\lie'(V)$) be the corresponding
	Lie superalgebra. The map $\varphi: \lie(V) \rightarrow \End(V)$
	(resp. $\varphi': \lie'(V) \rightarrow \End(V)$), given by formula
	(\ref{eq:nuevamentealvicio}), is a Lie superalgebra homomorphism. 
	\label{thm:k.lie.algebra.morphism}
\end{thm}
\section{Examples}\label{sub:examples}
\begin{ex}[$V(\cW_N)$ series] Let $X = (x, \xi^1, \dots, \xi^N)$, where
	 $x$ is even and $\xi^i$ are odd anticommuting variables
	 commuting with $x$. Consider the Lie algebra $\fg=W(1|N)$ of derivations of
	 $\mathbb{C}[X, X^{-1}]$.
	It is spanned by elements of the form $X^{j|J} \partial_x$ and $X^{j|J}
	\partial_{\xi^i}$ (cf. Example \ref{ex:w_n.series.aa}). Define the following $\fg$-valued formal distributions:
	\begin{equation}
		L(Z) = - \delta(Z,X) \partial_x, \qquad Q^i(Z) = - \delta(Z,X)\partial_{\xi^i}, \quad i = 1,
		\dots, N.
		\label{eq:w.ex.1.a}
	\end{equation}
	A straightforward computation shows that these distributions satisfy the
	following commutation relations:
	\begin{equation}
		\begin{aligned}
			{[L(Z)}, L(W)] &= \delta(Z,W) \partial_w L(W) + 2 \left( \partial_w
			\delta(Z,W)
			\right) L(W), \\
			[L(Z), Q^i(W)] &= \delta(Z,W) \partial_w Q^i(W) + \left(
			\partial_{\zeta^i} \delta(Z,W)
			\right) L(W) + \\ &\quad + \left( \partial_w \delta(Z,W) \right) Q^i(W), \\
			[Q^i(Z), Q^j(W)] &= \delta(Z,W) \partial_{\zeta^i} Q^j(W) + (-1)^{N}
			\left( \partial_{\zeta^i} \delta(Z,W) \right) Q^j(W) - \\ & \quad - (-1)^N \left(
			\partial_{\zeta^j} \delta(Z,W) \right) Q^i(W).
		\end{aligned}
		\label{eq:w.ex.2.b}
	\end{equation}
	In particular, the distributions (\ref{eq:w.ex.1.a}) are pairwise local. Let
	$\cF$ be the family of $\fg$-valued formal distributions
	\[
	\cF = \left\{ \partial_Z^{j|J} L(Z),\, \partial^{j|J}_{Z} Q^i(Z)\bigr| \quad j\geq 0,\, J
	\subset \{1,\dots, N\}, \, i = 1,
	\dots, N \right\}.\]  Then  
	 $(\fg, \cF)$ is an $N_W = N$ SUSY formal distribution Lie superalgebra.

	Let $\cW_N$ be the corresponding $N_W=N$ SUSY Lie conformal algebra. It is
	generated as a $\mathbb{C}[T, S^i]$-module by a vector $L$ of parity $N \mod 2$ and
	$N$-vectors $Q^i$, $i=1, \dots, N$, of parity $N+1 \mod 2$ satisfying the following
	$\Lambda$-brackets (which follow from (\ref{eq:w.ex.2.b}))
	\begin{equation}
		\begin{aligned}
			{[}L_\Lambda L] &= (T + 2 \lambda) L, \quad
			[Q^i_\Lambda Q^j] =  (S^i + \chi^i) Q^j - \chi^j Q^i, \\
			[L_\Lambda Q^i] &= (T + \lambda) Q^i + (-1)^N \chi^i L.
		\end{aligned}
		\label{eq:ex.w.n.1}
	\end{equation}
	When $N=0$, this is the centerless Virasoro conformal algebra, which admits
	a central extension defined by: 
	\begin{equation*}
		[L_\lambda L] = (T + 2 \lambda) L + \frac{\lambda^3}{12} C
	\end{equation*}
	where $C$ is even, central, and satisfies $T C = 0$.
	
	Translating the formulas in \cite{kacfattori}, it follows that when $N=1$,
	$\cW_1$ admits a central extension of the form:
	\begin{equation}
		\begin{aligned}
			{[}L_\Lambda L] &= (T + 2 \lambda) L, \quad
			[Q_\Lambda Q] =  S Q + 
			\frac{\lambda \chi}{3} C,\\
			[L_\Lambda Q] &= (T + \lambda) Q - \chi L +
			 \frac{\lambda^2}{6} C,
		\end{aligned}
		\label{eq:ex.w.n.3}
	\end{equation}
	where $C$ is even, central, and satisfies $T C = S C = 0$.
	
	When $N=2$, $\cW_2$  admits a central extension given by:
	\begin{equation*}
		\begin{aligned}
				{[}L_\Lambda L] &= (T + 2 \lambda) L,  &
			[L_\Lambda Q^i] &= (T + \lambda) Q^i +  \chi^i L, \\
			[Q^i_\Lambda Q^i] &=  S^i Q^i, &
			[Q^1_\Lambda Q^2] &= (S^1 + \chi^1) Q^2 - \chi^2 Q^1 +
			\frac{\lambda}{6} C, \\ 
		\end{aligned}
	\end{equation*}
	where $C$ is as above. It follows from \cite{kacleur},
	\cite{kacfattori} that these algebras do not admit central extensions for
	$N \geq 3$. 
	
	If $N \leq 2$, we let $V(\cW_N)$ be the universal enveloping $N_W=N$ SUSY vertex
	algebra of the central extension of $\cW_N$ as given by Theorem
	\ref{thm:universal.1}, and let $V(\cW_N)^c$ be the quotient of $V(\cW_N)$ by the ideal
	$(C - c \vac)_{(-1|N)} V(\cW_N)$, where $c \in \mathbb{C}$ is called the
	\emph{central charge}. 

	When $N=1$, expanding the superfields as
	\begin{equation*}
		Q (Z) = -J(z) + \theta G^+ (z), \qquad L(Z) = G^-(z) + \theta
		\left( L(z) + \frac{1}{2} \partial_z J(z) \right),
	\end{equation*}
	we check that the fields $L, J, G^\pm$ satisfy the commutation relations
	of the $N=2$ vertex algebra as defined in Example 
	\ref{ex:n=2}. 

	When $N \geq 3$, we let $V(\cW_N)$ be the universal enveloping $N_W = N$ SUSY vertex
	algebra of $\cW_N$. It follows from the definitions that $\lie(\cW_N)
	= W(1|N)$, the Lie
	superalgebra of derivations of $\mathbb{C}[X, X^{-1}]$. Also, 
	$\lie(\cW_N)_\leq = W(1|N)_\leq$ is the Lie superalgebra of derivations of
	$\mathbb{C}[X]$. Denote by $W(1|N)_<=\lie(\cW_N)_< \subset
	\lie(\cW_N)_\leq$ the Lie
	superalgebra of vector fields vanishing at the origin; it is
	spanned by vectors of the form $X^{j|J} \partial_x$ 
	 and $X^{j|J} \partial_{\xi^i}$, with $j + \sharp J
	> 0$.  
	\label{ex:w_n.conformal}
\end{ex}
\begin{defn}
	An $N_W=N$ SUSY vertex algebra $V$ is called \emph{conformal}  if there exists
	$N+1$ vectors $\nu, \tau^1,\dots, \tau^N$ in $V$ such that their associated superfields
	$L(Z) = \ys(\nu, Z)$ and $Q^i(Z) = \ys(\tau^i, Z)$ satisfy (\ref{eq:ex.w.n.1}) (or
	possibly a central extension) and moreover:
	\begin{itemize}
		\item $\nu_{(0|0)} = T$, 
		 $\tau^i_{(0|0)} =  S^i$,
		\item the operator $\nu_{(1|0)}$ acts diagonally with eigenvalues bounded 
			below and with finite dimensional eigenspaces.
	\end{itemize}
	
	If moreover, 
	 the action of $\lie(\cW_N)_<$ on $V$ can be
	exponentiated to the group of automorphisms of the disk $D^{1|N}$, we will say that
	$V$ is \emph{strongly conformal}. 
	 This amounts to the following extra condition:
	\begin{itemize}
		\item the operators $\nu_{(1|0)}$ and $\sum_{i=1}^N \sigma(e_i)
			\tau^i_{(0|e_i)}$ have integer
			eigenvalues.
	\end{itemize}

	If $a \in V$ is an eigenvector of the operator $\nu_{(1|0)}$ of eigenvalue $\Delta$, we
	say that $a$
	has
	 \emph{conformal weight} $\Delta$. This happens if $a$ satisfies
	 \begin{equation*}
		 {[}L_\Lambda a] = (T + \Delta \lambda) a + O(\lambda^2) + O(\chi^1, \dots, \chi^N),
	 \end{equation*}
	 If, moreover, $a$ satisfies $[L_\Lambda a] = (T+ \Delta \lambda) a$ we say that $a$ is
	 \emph{primary}. 
	
	 As in the ordinary vertex algebra case, the conformal weight $\Delta(a)$
	 is an important \emph{book-keeping device}: \[\Delta(Ta)= \Delta(a) + 1, \quad
	 \Delta(S^i a) = \Delta(a), \quad \Delta(a_{(j|J)}b) = \Delta(a) +
	 \Delta(b) - j - 1.\]
	 (Note that $\Delta(:ab:) = \Delta(a) + \Delta(b)$.) Furthermore,
	 letting $\Delta(\chi^i) = 0$ and $\Delta(\lambda) = 1$, all
	 terms in $[a_\Lambda b]$ have conformal weight $\Delta(a) + \Delta(b) - 1$. 
	\label{defn:w.conformal.definition}
\end{defn}
\begin{rem}
	It is clear from this definition that the $N_W=N$ SUSY vertex algebra
	$V(\cW_N)^c$
	defined in Example \ref{ex:w_n.conformal} is strongly conformal. 
	\label{rem:w_n.conf.}
\end{rem}
\begin{ex}[Free Fields] As an example of a strongly conformal $N_W=N$ SUSY vertex algebra,
	we will compute explicitly the free fields case, namely, let
	$\alpha, \,C$ be even vectors and let $\varphi$ be an odd vector. Consider
	the $N_W=N$ SUSY Lie conformal algebra generated by these three vectors, where $C$ is
	central and annihilated by $\nabla$, and the other commutation
	relations are:
	\begin{equation*}
		{[}\alpha_\Lambda \varphi] = C.
	\end{equation*}
	Let $\tilde{F}_N$ be its universal enveloping $N_W=N$ SUSY vertex algebra and $F_N$
	its quotient by the ideal $(C -  \vac)_{(-1|N)} \tilde{F}_N$.
	
	Expanding the superfields 
	\begin{equation*}
		\alpha(Z) = a(z) + \theta \psi(z), \qquad \varphi(Z) = \phi(z) +
		\theta b(z)
	\end{equation*}
	we find that the fields $a, b, \psi$ and $\phi$ generate the well known
	$bc-\beta\gamma$-system, namely, the non-trivial $\lambda$-brackets
	are (up to skew-symmetry):
	\begin{equation*}
		{[}b_\lambda a] = [\psi_\lambda \phi] = 1.
	\end{equation*}

	When $N=1$, this SUSY vertex algebra
	admits a $N_W = 1$
	 strongly conformal structure with:
	\begin{equation*}
			\nu = \alpha_{(-2|1)}\varphi_{(-1|1)} \vac, \quad
			\tau = - \alpha_{(-1|0)}\varphi_{(-1|1)} \vac,
	\end{equation*}
	and central charge $c = 3$. 
	The associated fields $L = \ys(\nu, Z)$ and $Q = \ys(\tau, Z)$
	are: 
	\begin{equation*}
		L = :(T\alpha)\varphi: \qquad Q = - :(S\alpha)\varphi:.
	\end{equation*}
	In order to check the commutation relations (\ref{eq:ex.w.n.3}) we use the
	non-commutative Wick formula (\ref{eq:wick.proof.8}) to find:
	\begin{equation*} 
		{[}\alpha_\Lambda L ] = T\alpha, \quad  [\varphi_\Lambda L ] =
		\lambda \varphi, \quad
		[\alpha_\Lambda Q] = S \alpha, \quad  [\varphi_\Lambda Q] = - \chi
		\varphi.
	\end{equation*}
	And now by skew-symmetry and sesquilinearity we obtain:
	\begin{xalignat}{2}
		{[}L_\Lambda \alpha] &= T\alpha & [L_\Lambda \varphi] &= (\lambda
		+ T) \varphi \label{eq:conformal.alrepedo.1} \\
		[L_\Lambda T\alpha] &= (T+\lambda)T \alpha & [L_\Lambda S\alpha]
		&= (S+\chi)T\alpha \\
		[Q_\Lambda \alpha] &= S\alpha & [Q_\Lambda \varphi] &= (S +
		\chi) \varphi \\
		[Q_\Lambda S\alpha] &= - \chi S \alpha.
		\label{eq:w.free.ex.1.e}
	\end{xalignat}
	Formula (\ref{eq:conformal.alrepedo.1}) says that $\alpha$ and $\varphi$
	are primary fields of conformal weight $0$ and $1$ respectively.
	With these formulas and using again the Wick formula
	(\ref{eq:wick.proof.8}) we obtain
	\begin{multline*}
		{[}L_\Lambda L] = [L_\Lambda :(T\alpha) \varphi:] = :( (T +\lambda)
		T\alpha)\varphi: + :T\alpha(\lambda + T) \varphi: = \\ = 2 \lambda
		L + :(T(T\alpha))\varphi: + :T\alpha T\varphi: = (T + 2 \lambda) L,
	\end{multline*}
	since the integral term obviously vanishes. For the other commutation
	relations we compute:
	\begin{multline*}
		{[}Q_\Lambda Q] = - [Q_\Lambda :(S\alpha)\varphi:] = \chi :(S \alpha)
		\varphi: + :S\alpha (\chi + S) \varphi: + \int_0^\Lambda[\chi
		S\alpha_\Gamma \varphi] d \Gamma = \\ = :S\alpha S\varphi: +
		\int_0^\Lambda(\eta - \chi) \eta d \Gamma =  S Q + \lambda \chi.
	\end{multline*}
	Finally for the last commutator we find:
	\begin{multline*}
		{[}L_\Lambda Q] = - [L_\Lambda :(S\alpha)\varphi:] = -:( (S +
		\chi)T\alpha) \varphi: - :S\alpha (\lambda +T) \varphi: - \\ -
		\int_0^\Lambda[ (S+\chi)T\alpha_\Gamma \varphi] d \Gamma  =
		TQ - \chi :(T\alpha)\varphi: + \lambda Q + \int_0^\Lambda (\eta -
		\chi) \gamma d \Gamma = \\ = (T+\lambda)Q - \chi L +
		\frac{\lambda^2}{2}.
	\end{multline*}
	According to (\ref{eq:ex.w.n.3}) this is a conformal $N_W=1$ SUSY
	vertex 
	algebra with central charge $3$. It is easy to check that this SUSY vertex
	algebra is indeed strongly conformal.
	\label{ex:w.free.ex.1}
\end{ex}
\begin{ex}[$V(\cK_N)$ series] Consider now the Lie subsuperalgebra $K(1|N) \subset W(1|N)$
	consisting of those derivations of $\mathbb{C}[X, X^{-1}]$ preserving the form
	\[\omega = dx + \sum_{i =1}^N \xi^i d\xi^i,\]up to multiplication by a function
	(cf. Example \ref{ex:k_series.begin}). Define the following $\fg$-valued formal
	distribution:
	\begin{equation*}
		G(Z) = - 2\delta(Z,X) \partial_x - (-1)^N \sum_{i = 1}^N \left( D^i_X
		\delta(Z,X)
		\right) D^i_X.
	\end{equation*}
	It follows from (\ref{eq:fattorikac1}) that its $Z$-coefficients form a
	basis of $K(1|N)$. 
		A straightforward computation shows that this formal distribution satisfies the
	following commutation relation:
	\begin{multline*}
		{[G(Z)},G(W)] = 2 \delta(Z,W) \partial_w G(W) + (4-N) \left( \partial_w
		\delta(Z,W) \right) G(W) + \\ +(-1)^N \sum_{i = 1}^N \left( D^i_W \delta(Z,W)
		\right) D^i_W G(W),
	\end{multline*}
	in particular, the pair of $\fg$-valued formal distributions $(G(Z),G(Z))$ is local.
	Letting $\cF = \left\{ \partial_Z^{j|J} G(Z), \; j \geq 0,\, J \subset \{1, \dots, N\}
	\right\}$, we obtain an $N_K = N$ formal distribution Lie
	superalgebra $(\fg, \cF)$. 

	Let $\cK_N$ be the associated $N_K = N$ SUSY Lie conformal algebra. It is generated
	as a free $\cH$-module by a vector $G$ of parity $N \mod 2$ satisfying the following
	$\Lambda$-bracket (for the definition of the
	algebra $\cH$ see \ref{no:whattochange})
	\begin{equation}
		[G_\Lambda G] = \left( 2 T + (4-N) \lambda + \sum_{i = 1}^N
		\chi^i S^i \right) G.
		\label{eq:ex.k_n.1}
	\end{equation}
	When $N \leq 3$, $\cK_N$ admits a non-trivial central extension, obtained
	by adding the term
	$\tfrac{\lambda^{3-N}\chi^N}{3} C$ to the RHS of (\ref{eq:ex.k_n.1}), where $C$ is even, central
	and satisfies $T C = S^i C = 0$, cf. \cite{kacfattori}. 
	
	When $N = 4$, $\cK_4$ admits a central extension, obtained by
	adding the term $\lambda C$ to the RHS of (\ref{eq:ex.k_n.1}). It follows from
	\cite{kacleur}, \cite{kacfattori} that $\cK_N$ does
	not admit central extensions when $N > 4$.

	When $N \leq 4$, we let $V(\cK_N)$ be the universal enveloping $N_K = N$ SUSY vertex
	algebra of the central extension of $\cK_N$ and define $V(\cK_N)^c$ to be its quotient by
	the ideal $(C - c \vac)_{(-1|N)} V(\cK_N)$, where $c \in \mathbb{C}$ is called the
	\emph{central charge}. When $N \geq 5$, we let $V(\cK_N)$ be the universal enveloping
	$N_K=N$ SUSY vertex algebra of $\cK_N$. 

	In the case $N=1$, if we expand the corresponding superfield 
	 as
	\begin{equation*}
		G(z,\theta) = G(z) + 2 \theta L(z),
	\end{equation*}
	we find that the fields $G(z)$ and $L(z)$ generate a Neveu Schwarz vertex
	algebra of central charge $c$ as in Example \ref{ex:2.11}. 

	When $N=2$ expanding the corresponding superfield as (cf.
	\ref{eq:tau.k.2.field})
	\begin{equation*}
		G(z, \theta^1,\theta^2) = \sqrt{-1} J(z) +\theta^1
		G^{(2)}(z) - \theta^2  G^{(1)}(z) + 2 \theta^1 \theta^2 L(z)
	\end{equation*}
	We find that the corresponding fields $J, L, G^\pm$ satisfy the
	commutation relations of the $N=2$ vertex algebra as in Example
	\ref{ex:n=2}. 

	When $N=4$ the corresponding $N_K=4$  SUSY vertex algebra is not simple. Indeed the
	SUSY Lie conformal superalgebra $\cK'_4 \subset \cK_4$ generated by
	$S^i G$, $i=1,\dots,4$ is an ideal. 
	The central
	extension of $\cK_4$ described above restricts to a central extension of
	$\cK'_4$
        whose cocycle is given by:
	\begin{equation}
		\alpha_1 (S^i G, S^j G) = -  \chi^i \chi^j C_1.
		\label{eq:cocyclederived}
	\end{equation}
	This SUSY Lie conformal algebra admits another central extension given by (cf.
	\cite{kacfattori}):
	\begin{equation}
		\alpha_2 (S^i G, S^j G) =  \chi^1 \chi^2 \chi^3 \chi^4 C_2.
		\label{eq:cocyclederived2}
	\end{equation}
	We let ${V(\cK'_4)^{c_1,c_2}}$ be the corresponding $N_K = 4$ SUSY vertex
	algebra with $C_i = c_i$.

	Note that $\lie(\cK_N) = K(1|N)$ by definition, while $\lie(\cK_N)_\leq
	= K(1|N)_\leq$ is the Lie
	superalgebra of
	 regular vector fields preserving $\omega$ up to multiplication by a
	 function. We will denote by
	$K(1|N)_< = \lie(\cK_N)_< \subset \lie(\cK_N)_\leq$ the Lie superalgebra of regular vector
	fields, preserving $\omega$ up to multiplication by a function, and vanishing at the origin. 

	A field $G$ satisfying the commutation relations (\ref{eq:ex.k_n.1})
	with a central extension,
	will be called a \emph{super Virasoro} field. 
	\label{ex:k.n.series}
\end{ex}
\begin{defn}
	Let $N \leq 4$, an $N_K=N$ SUSY vertex algebra $V$ is called \emph{conformal} if there exists
	a vector $\tau \in V$ (called the conformal vector) such that the corresponding field
	$G(Z) = \ys(\tau, Z)$ satisfies (\ref{eq:ex.k_n.1}) with a central
	extension,
	and moreover
	\begin{itemize}
		\item $\tau_{(0|0)} = 2 T$, 
		 $\tau_{(0|e_i)} = \sigma(N\setminus e_i, e_i) S^i$,
		\item the operator $\tau_{(1|0)}$ acts diagonally with eigenvalues bounded
			 below and finite dimensional eigenspaces.
	\end{itemize}

	If moreover,  the
	representation of $\lie(\cK_N)_<$ can be exponentiated to the group of
	automorphisms of the  disk $D^{1|N}$ preserving the SUSY structure 
	\begin{equation*}
		\omega = d x + \sum_{i = 1}^N \xi^i d \xi^i,
	\end{equation*}
	we will say that $V$ is \emph{strongly conformal}. 
	This amounts to the extra condition \begin{itemize} \item the operator $\tau_{(1|0)}$ has
			integer eigenvalues and, if $N=2$, the operator
			$\sqrt{-1} \tau_{(0|N)}$ has integer eigenvalues. \end{itemize}

	If a vector $a$ in a conformal $N_K=N$ SUSY vertex algebra $V$ is an eigenvector of
	$\tau_{(1|0)}$ with eigenvalue $2 \Delta$, we say that $a$ has \emph{conformal weight}
	$\Delta$. This happens iff $a$ satisfies 
	\begin{equation*}
		{[}G_\Lambda a ] = \left( 2 T + 2 \Delta \lambda + \sum_{i=1}^N \chi^i S^i \right) a +
		O(\Lambda^2),
	\end{equation*}
	where $O(\Lambda^2)$ denotes a polynomial in $\Lambda$ with vanishing constant and
	linear terms. If, moreover, $a$ satisfies \[ [G_\Lambda a] = \left( 2 T + 2
	\Delta \lambda +
	\sum_{i=1}^N \chi^i S^i \right) a, \] we say that $a$ is \emph{primary}.
	For example, formula (\ref{eq:ex.k_n.1}) with a central extension, says that $G$ has conformal weight
	$2 - N/2$, and it is primary if the central extension is trivial.  As in
	the $N_W= N$ case, the conformal weight is an important \emph{book-keeping
	device}
	\[\Delta(Ta) = \Delta(a) + 1, \quad \Delta(S^i a) = \Delta(a) +
	\frac{1}{2}, \quad \Delta(a_{(j|J)}b) = \Delta(a) + \Delta(b) - j - 1 +
	\frac{N - \sharp J}{2}.\](Note that $\Delta(:ab:) = \Delta(a) +
	\Delta(b)$.) Furthermore,
	letting $\Delta(\lambda) = 1$ and $\Delta(\chi^i) = 1/2$, all terms in
	$[a_\Lambda b]$ have conformal weight $\Delta(a) + \Delta(b) - 1 + N/2$.
		\label{defn:k.conf.vertex.def}
\end{defn}
\begin{rem}
	The $N_K=N$ SUSY vertex algebra $V(\cK_N)^c$ defined in Example \ref{ex:k.n.series} is
	strongly conformal when $N < 4$, and ${V(\cK'_4)^{c_1,c_2}}$ is
	strongly conformal as well.
	\label{rem:k.is.stronly}
\end{rem}
\begin{ex}(Free fields) The well-known \emph{boson-fermion} system is an
	$N_K=N$ vertex algebra generated by one superfield. Let $\Psi$ be a
	vector of parity $(-1)^N$, $C$ an even vector, and define a $N_K=N$ SUSY Lie conformal
	algebra generated by $\Psi$ and $C$ where $C$ is central, satisfies $T C = S^i C = 0$
	and the remaining commutation relations are:
	\begin{equation*}
		{[}\Psi_\Lambda \Psi] = \Lambda^{1|N} C,
	\end{equation*}
	when $N$ is even, and 
	\begin{equation*}
		{[}\Psi_\Lambda \Psi] = \Lambda^{0|N} C,
	\end{equation*}
	when $N$ is odd. Skew symmetry is clear and the Jacobi identity is
	obvious since all triple brackets vanish. We let $\tilde{B}_N$ be the
	corresponding universal enveloping $N_K=N$ SUSY vertex algebra, and let $B_N$ be its
	quotient by the ideal $(C - \vac)_{(-1|N)} \tilde{B}_N$.  

	To show an application of the above formalism, as well as the subtleties
	involved in calculations, we will show explicitly that the $N_K=1$ SUSY
	vertex algebra $B_1$ is
	conformal, the corresponding super Virasoro field being
	\begin{equation*}
		G = :(S \Psi) \Psi: + m T \Psi, \qquad m \in \mathbb{C}.
	\end{equation*}
	Indeed, from sesquilinearity (\ref{eq:k.sesqui.1}) and skew-symmetry
	(\ref{eq:k.skew.1}) we find
	\begin{equation*}
		{[}\Psi_\Lambda S\Psi ] = (S+\chi)\chi = \lambda, \qquad [S\Psi_\Lambda \Psi] = -
		\lambda,
	\end{equation*}
	where we used $[S, \chi] = 2 \lambda$ and $\chi^2 = -\lambda$. 
	Using sesquilinearity once more we get:
	\begin{equation*}
		{[}S\Psi_\Lambda S\Psi] = \chi \lambda, \qquad [\Psi_\Lambda T
		\Psi] = \lambda \chi.
	\end{equation*}
	Now we can use the $N_K=1$ version of the non-commutative Wick formula
	(\ref{eq:wick.proof.8}) to find:
	\begin{gather*}
			{[}\Psi_\Lambda G ] = ( \lambda + \chi S)
			\Psi  + m \lambda \chi, \quad 
			[S\Psi_\Lambda G] = \lambda (\chi - S) \Psi -
			m \lambda^2  \\
			[T\Psi_\Lambda G ] = - \lambda (\lambda + 
			 \chi S) \Psi - m \lambda^2 \chi,
	\end{gather*}
	where we note that all the integral terms vanish. Using skew-symmetry
	again we get:
	\begin{equation}
		\begin{gathered}
			{[}G_\Lambda \Psi] = (\lambda + 2 T + \chi S) \Psi
			- m \lambda \chi, \quad
			[G_\Lambda S\Psi ] = (\lambda + T)(  \chi + 2 S
			  ) \Psi - m \lambda^2 \\
			[G_\Lambda T\Psi] = (T + \lambda)(\lambda + 2 T +
			\chi S) \Psi - m \lambda^2 \chi.
		\end{gathered}
		\label{eq:free.fields.7}
	\end{equation}
	With these formulas we can use again (\ref{eq:wick.proof.8}) to get:
	\begin{multline*}
		{[}G_\Lambda G]  = :\left((\lambda + T) (\chi + 2 S) \Psi -
		m \lambda^2 \right) \Psi: + \\ +  :S\Psi \left(  \lambda + 2 T + \chi
		S) \Psi - m \lambda \chi \right): + m (T+ \lambda) (\lambda +
		2T + \chi S) \Psi  - m^2 \lambda^2 \chi
	\end{multline*}
	where again the integral term is easily seen to vanish. Note that from
	quasi-commutativity of the normally ordered product we find $:\Psi \Psi: =
	0$, from where the expression above reduces to:
	\begin{multline*}
		2 \lambda :(S\Psi) \Psi: + \chi :(T\Psi)\Psi:  + 2 :(S^3 \Psi) \Psi: -
		m \lambda^2 \Psi + 
		 \\ +  :S\Psi \left(  (\lambda + 2 T + \chi
		S) \Psi - m \lambda \chi \right): + m (T+ \lambda) (\lambda +
		2T + \chi S) \Psi  - m^2 \lambda^2 \chi.
	\end{multline*}
	Expanding this expression and after a simple cancellation we find
	\begin{equation*}
		{[}G_{\Lambda} G] = \left( 2 T + 3\lambda  + \chi S \right) G -
		m^2 \lambda^2 \chi.
	\end{equation*}
	Therefore $B_1$ is a strongly conformal $N_K=1$ SUSY vertex algebra with
	central charge $-3 m^2$. Note that, by (\ref{eq:free.fields.7}), $\Psi$ has
	conformal weight $1/2$ (but it is not primary if $m \neq 0$). 
	Expanding the superfield
	\begin{equation*}
		\Psi(Z) = \varphi(z) + \theta \alpha(z),
	\end{equation*}
	we find easily that
	\begin{equation*}
		[\varphi_\lambda \varphi ] = 1, \qquad [\alpha_\lambda\alpha] = \lambda,
	\end{equation*}
	hence the name \emph{boson-fermion} system.
	\label{ex:k.free.fields}
\end{ex}
\begin{ex}(Super Currents) Let $\fg$ be a simple or abelian Lie superalgebra with a
	non-degenerate invariant
	supersymmetric bilinear form $(\, , \,)$. If $N$ is even, then  we
	define a SUSY Lie conformal algebra (either $N_K=N$ or $N_W=N$) generated by
	$\fg$ with commutation relations:
	\begin{equation*}
		[a_\Lambda b] = [a,b] + (k+h^\vee) (a,b) \lambda \qquad \forall a,\,b \,
		\in \fg,
	\end{equation*}
	where $2 h^\vee$ is the eigenvalue of the Casimir operator on $\fg$.

	When $N$ is odd we let $\Pi\fg$ be $\fg$ with reversed parity, and for each
	element $a \in \fg$ we let $\bar{a}$ be the same element thought in
	$\Pi\fg$. In this case we define a SUSY Lie conformal algebra generated by
	$\Pi\fg$ with commutation relations:
	\begin{equation*}
		[\bar{a}_\Lambda \bar{b}] = (-1)^a \left( \overline{[a,b]} + (k+h^\vee) (a,b)
		\sum_{i = 1}^N \chi^i \right).
	\end{equation*}
	We let $V^k(\fg)$ be the associated universal enveloping SUSY vertex
	algebra\footnote{Here as before, we are considering a central extension of a SUSY Lie
	conformal algebra, and then we identify the central element with a multiple of the
	vacuum vector in the universal enveloping SUSY vertex algebra.},
	either the $N_K=N$ or the $N_W=N$ vertex algebra, the choice will be
	clear in each context, as well as the value of $N$. 

	When $N=1$, the corresponding $N_K=1$ SUSY vertex algebra is strongly
	conformal,
	 the corresponding conformal vector is 
	\begin{equation*}
		\tau = \frac{1}{k + h^\vee} \left( \sum (-1)^{a^i} :(S\bar{a}^i)
		\bar{b}^i: +
		\frac{1}{3(k + h^\vee)} \sum ([a^i, a^j],a^r) :\bar{b}^i :\bar{b}^j
		\bar{b}^r: :
		\right),
	\end{equation*}
	where $\{a^i\}$ and $\{b^i\}$ are dual bases for $\fg$ with respect to $(\,,\,)$.  This
	is known as the \emph{Kac-Todorov} construction \cite{kactodorov}. The
	super Virasoro field $Y(\tau,Z)$ has central charge \[ c = \frac{k \mathrm{sdim} \fg}{k
	+ h^\vee} + \frac{\mathrm{sdim} \fg}{2},\] and the fields $\bar{a} \in
	\bar{\fg}$ have conformal weight $1/2$.
	\label{ex:ex.currents.2}
\end{ex}
\begin{ex}($N=2$ vertex algebra) As a vertex algebra it
  is generated by $4$ fields (cf. Example \ref{ex:n=2}). As we have seen in Example
  \ref{ex:k.n.series}, this is a $N_K=2$ SUSY vertex algebra generated by one field
  $G$. On the other hand, the $N=2$ vertex algebra admits an embedding of the $N=1$
  vertex algebra. Therefore we can view the $N=2$ vertex algebra as an $N_K=1$ SUSY
  vertex algebra. As such, this algebra is generated by two superfields $G$ and
  $J$, where
  $G$ is a super Virasoro field of central charge $c$ and $J$ is even, primary of
  conformal weight $1$. The remaining
  $\Lambda$-bracket is given by:
  \begin{equation*}
      {[J}_\Lambda J] =  G + \frac{c}{3} \lambda \chi.
  \end{equation*}
  This is computed using the decomposition Lemma \ref{lem:k.decomp}.
  \label{ex:k.n=2.double}
\end{ex}
\begin{ex}($N=4$ vertex algebra) As a vertex
	algebra it is generated by 8 fields: a Virasoro field, three currents (for
	the Lie algebra $\fs\fl_2$) and four fermions \cite[p. 187]{kac:vertex}. This
	vertex algebra admits and embedding of the Neveu Schwarz vertex algebra,
	therefore we can consider a $N_K=1$ SUSY vertex algebra structure on it. As a $N_K=1$
	vertex algebra, it is of rank $3|1$, generated by an $N=1$ conformal vector $G$
	with central charge $c$ and three even vectors $J^i$, $i=1,2,3$. Each pair
	$(G,J^i)$ generates an $N=2$ vertex algebra, viewed as an $N_K=1$ SUSY
	vertex algebra, as in the previous example. The
	remaining commutation relations are:
	\begin{equation*}
		{[ J^i}_\Lambda J^j] =  \sqrt{-1}
			 \varepsilon^{ijk} (S + 2 \chi) J^k, \qquad i \neq j,
	\end{equation*}
	where $\varepsilon$ is the totally antisymmetric tensor. 

	This vertex algebra is the universal enveloping vertex algebra of the
	central extension of the superconformal Lie algebra $S(1|2;0)$ (cf.
	\cite{kacfattori}). 
	\label{ex:ex.n4}
\end{ex}
\begin{ex} ($bc-\beta \gamma$ system) This is a $N_K=1$ SUSY vertex algebra generated
	by $n$ even fields $B^i$ and $n$ odd fields $\Psi^i$. The only
	non-vanishing $\Lambda$-brackets (up to skew-symmetry) are:
	\begin{equation*}
		{[B^i}_\Lambda \Psi^j] = \delta_{ij}.
	\end{equation*}
	This SUSY vertex algebra is
	strongly conformal with super Virasoro field
	\begin{equation*}
		G = \sum_{i = 1}^n \left( :(SB^i) (S\Psi^i): + :(TB^i) \Psi^i:
		\right),
	\end{equation*}
	and central charge $3n$. The fields $B^i$ (resp. $\Psi^i$) are primary of
	conformal weight $0$ (resp. $1/2$). 

	Let $\sigma^s_{ij}$, $s = 1,2,3$, be
	three $n \times n$ matrices satisfying 
	\begin{equation*}
		\sigma^i \sigma^j = \sqrt{-1} \varepsilon^{ijk} \sigma^k, \qquad
		(\sigma^s)^2 =  \Id.  
	\end{equation*}
	The fields
	\begin{equation*}
		J^i = \sum_{j,k=1}^n \sigma^i_{jk} :SB^j \Psi^k:, \quad i =
		1,\,2,\,3,
	\end{equation*}
	together with $G$ generate an $N=4$ vertex algebra as in
	the previous example\footnote{Note however that these fields $J^i$
	differ
        from those used in \cite{heluani2} by a factor of $\sqrt{-1}$.} (cf. \cite{heluani2}).
	\label{ex:bc}
\end{ex}
\begin{ex} Here we explain briefly the construction of the chiral de Rham complex
	of a smooth manifold introduced in \cite{malikov}, using the formalism of
	$N_K=1$ SUSY vertex algebras \cite{heluani2}.
	     Let $U$ be a differentiable manifold.
    Let $\cT$ be the tangent bundle of $U$ and $\cT^*$ be its
    cotangent bundle. We let $T=\Gamma(U, \cT)$ be the space of vector
    fields on $U$ and $A = \Gamma(U, \cT^*)$ be the space of
    differentiable $1$-forms on $U$. We let $\cC=\cC^\infty(U)$ be the
    space of differentiable functions on $U$. Denote by
    \begin{equation*}
        <\,,\,> : A \otimes T \rightarrow \cC
        \label{cdr.2.1}
    \end{equation*}
    the natural pairing, and, as before, by $\Pi$ the functor of
    change of parity.

    Consider now an $N_K = 1$ SUSY Lie conformal algebra $\cR$ generated by the
    vector superspace
    \begin{equation*}
        \cC \oplus \Pi T \oplus A \oplus \Pi A.
        \label{cdr.2.2}
    \end{equation*}
    That is, we consider differentiable functions (to be denoted
    $f,g, \dots$) as even elements, vector fields $X, Y, \dots$ as
    odd elements, and finally we have two copies of the space of
    differential forms. For differential forms, which we consider to be even
    elements, $\alpha, \beta, \dots
    \in A$,
    we will denote the corresponding elements of $\Pi A$ by
    $\bar{\alpha}, \bar\beta, \dots$. The nonvanishing $\Lambda$-brackets
     in
    $\cR$ are given by (up to skew-symmetry):
    \begin{equation*}
        \begin{aligned}
            {[X}_\Lambda f] &= X(f), &
            [X_\Lambda Y] &= [X,Y]_{\mathrm{Lie}}, \\
            [X_\Lambda \alpha] &= \lie_X \alpha + \lambda
            <\alpha, X>,  &
            [X_\Lambda \bar{\alpha}] &= \overline{\lie_X \alpha} +
            \chi <\alpha, X>,
        \end{aligned}
        \label{cdr.2.3}
    \end{equation*}
    where $[,]_{\mathrm{Lie}}$ is the Lie bracket of vector fields
    and $\lie_X$ is the action of $X$ on the space of differential
    forms by the Lie derivative. The fact that (\ref{cdr.2.3})
    satisfies the Jacobi identity is a  straightforward
    computation (cf. \ref{no:intro.4}). 

    We let $V(U)$ be the corresponding universal enveloping $N_K=1$
    SUSY vertex algebra of $\cR$.  This vertex
    algebra is too big. We impose some relations in $V(U)$ as follows.
    Let $1_U$ denote the constant function $1$ in $U$. Let $d: \cC
    \rightarrow A$ be the de Rham differential. Define $I(U) \subset
    V'(U)$ to be the ideal generated by elements of the form:
    \begin{equation*}
	    \begin{aligned}
        :fg: - (fg), && :fX: - (fX), &&
         :f\alpha: - (f\alpha), &&  :f \bar{\alpha}:
         - (\overline{f\alpha}), \\
         1_U - \vac, && Tf - df, && Sf - \overline{df}. &&
 \end{aligned}
    \end{equation*}
    Define the $N_K = 1$ SUSY vertex algebra
    \begin{equation*}
        \Omega^{\mathrm{ch}} (U) := V(U)/I(U).
    \end{equation*}
    The following theorem is a reformulation of the corresponding result in
    \cite{LL}:
\begin{thm}\hfill
    \begin{enumerate}
    \item Let $M \subset \mathbb{R}^n$ be an open submanifold. The
    assignment $U \mapsto
    \Omega^{\mathrm{ch}}(U)$ defines a sheaf of SUSY vertex algebras
    $\Omega^{\mathrm{ch}}_M$ on $M$.
        \item For any diffeomorphism of open sets $M'
        \xrightarrow{\varphi} M $
         we obtain a canonical isomorphism of
        SUSY vertex algebras $\Omega^{\mathrm{ch}}(M)
        \xrightarrow{\Omega^{\mathrm{ch}}(\varphi)}
        \Omega^{\mathrm{ch}}(M')$. Moreover, given
        diffeomorphisms $M'' \xrightarrow{\varphi'} M'
        \xrightarrow{\varphi} M$, we have
        $\Omega^{\mathrm{ch}}(\varphi \circ \varphi') =
        \Omega^{\mathrm{ch}}(\varphi') \circ
        \Omega^{\mathrm{ch}}(\varphi)$.
    \end{enumerate}
    \label{thm:cdr.3}
\end{thm}

    This theorem allows one to construct a sheaf of SUSY vertex
    algebras
    in the Grothendieck topology on $\mathbb{R}^n$ (generated by open
    embeddings). This in turn allows one to attach to any smooth manifold
    $M$, a sheaf of SUSY vertex algebras $\Omega^{\mathrm{ch}}_M$,
    called the \emph{chiral de Rham complex of $M$}.
    \label{ex:chiralderham}
\end{ex}
\begin{ex} (Free Fields) We can generalize Examples \ref{ex:k.free.fields},
	\ref{ex:bc}, and \ref{ex:w.free.ex.1} as follows. Let $A = A_{\bar{0}}
	\oplus A_{\bar{1}}$ be a vector superspace, and let $(\, , \,)$ be a
	non-degenerate bilinear form in $A$. Recall that the bilinear form $(\,,\,)$ is said to
	be of parity $p \in \mathbb{Z}/2\mathbb{Z}$ if $(a,b)=0$ unless $p(a)+p(b)=p$, and it
	is supersymmetric (resp. skew-supersymmetric) if $(a,b)=(-1)^{ab}(b,a)$ (resp $(a,b) =
	-(-1)^{ab}(b,a)$). 

	Let $\cH = \mathbb{C}[T, S^i]$ in the $N_W=N$ case, and let $\cH$ be defined as in
	\ref{no:whattochange} in the $N_K=N$ case. Let \[\cR = \cH \otimes A \oplus
	\mathbb{C}C,\]where $C$ is an even element such that $TC = S^iC = 0$. Given a non-zero
	homogeneous polynomial $Q(\Lambda)$ of degree $s$ (in PBW basis) and parity $p$, define
	the following $\Lambda$-bracket on $A \oplus \mathbb{C}C$:
	\begin{equation*}
		{[a}_\Lambda b] = Q(\Lambda) (a,b)C, \quad a, b \in A,\; \text{and } C\;
		\text{central},
		\label{lamda-pedito}
	\end{equation*}
	and extend it to $\cR$ by sesquilinearity. Then the Jacobi identity automatically holds
	since all triple brackets are zero. Skewsymmetry holds if and only if
	\begin{equation*}
		(a,b) = -(-1)^{ab} (-1)^{N+s} (b,a).
		\label{lambda-pedaso}
	\end{equation*}
	Thus, (\ref{lamda-pedito}) defines a structure of a SUSY Lie conformal algebra,
	provided that (\ref{lambda-pedaso}) holds together with the following parity condition:
	\begin{equation*}
		p + p( (\, ,\,)) = N \mod 2.
		\label{}
	\end{equation*}
	Thus, $\cR$ is a SUSY Lie conformal algebra if and only if $N+s$ is even (resp. odd)
	and the bilinear form $(\, ,\,)$ is supersymmetric (resp. skew-supersymmetric) of
	parity $(N - p) \mod 2$. 
	The corresponding \emph{free field} SUSY vertex algebra
	$F(A,Q)$ is the quotient of the universal enveloping vertex algebra $V(\cR)$ by the
	ideal $(C - \vac)_{(-1|N)} V(\cR)$. 
\end{ex}
\begin{ex} ($\mathrm{Spin}_7$ vertex algebra) In \cite{vafa}, Shatashvili and Vafa constructed a vertex algebra
	associated to any manifold with $\mathrm{Spin}_7$ holonomy. This vertex
	algebra is generated by four fields and comes equipped with an
	$N=1$ superconformal vector, therefore we can view it as an $N_K = 1$ SUSY vertex
	algebra. As such, it is generated by a super Virasoro field $G$ of central charge
	$1/2$, and a (non-primary) even field $X$ of
	conformal weight $2$. The corresponding $\Lambda$-brackets, derived
	from the OPE in \cite{vafa}, using Lemma \ref{lem:k.decomp}, are:
\begin{equation*}
	\begin{aligned}
		{[}G_\Lambda X] &= \left( 2T + \chi S + 4 \lambda \right) X +
		\frac{\chi \lambda}{2} G + \frac{2}{3} \lambda^3,\\
		[X_\Lambda X] &= \left( \frac{5}{2} TS X + \frac{5}{4} T^2 G + 6 :G X: \right)
		+ \\ & \quad +8 \left( \chi T + \lambda S + 2 \lambda \chi\right) X +  \frac{15}{4} \lambda (T + \lambda) G
		+ \frac{8}{3} \lambda^3 \chi. 
	\end{aligned}
	\label{}
\end{equation*}
	Note that this $\Lambda$-bracket is quadratic in the generating fields.
	Thus, 
	 this SUSY vertex algebra is not the universal enveloping SUSY
	vertex algebra of a SUSY Lie conformal algebra.
	Expanding these superfields as:
\begin{equation*}
	G(Z) = G(z) + 2 \theta T(z), \qquad X(Z) = \tilde{X}(z) + \theta
	\tilde{M}(z),
	\label{}
\end{equation*}
we obtain the generating fields as in \cite{vafa}. 
	\label{ex:vafa}
\end{ex}
\begin{ex} (Odake's vertex algebra) In \cite{odake}, Odake constructed a vertex
	algebra, generated by eight fields, as
	an extension of the $N=2$ vertex superalgebra, associated to manifolds with
	$\mathrm{SU}_3$ holonomy. It carries therefore an
	$N=1$ superconformal vector, so we can view this algebra as an $N_K=1$ SUSY
	vertex algebra. As such, it is generated by two superfields $G, J$ forming
	an $N=2$ vertex algebra of central charge $9$, as in Example
	\ref{ex:k.n=2.double}, and two odd
	superfields $X^\pm$, primary of conformal weight $3/2$. The remaining
	$\Lambda$-brackets are as follows:
\begin{equation*}
	\begin{aligned}
		{[J}_\Lambda X^\pm] &= \pm (S + 3\chi) X^\pm, \qquad
		{[X^\pm}_\Lambda X^\pm] = 0,\\
	[{X^+}_\Lambda X^-] &= \left( :J G: + :J SJ: + T G + TS J \right) + \\ & \quad +\chi \left( :JJ: +
	T J
	\right) + \lambda (G + SJ) +2 \lambda \chi J + \lambda^2 \chi.
\end{aligned}
\end{equation*}
	Note that this relations are also quadratic in the generators as in the
	previous example. Expanding the generating superfieds as 
	\begin{equation*}
		\begin{aligned}
			J(Z) &= I(z) + \theta \frac{1}{\sqrt{2}} (\bar{G}(z) -
			G(z)), & G(Z) &= \frac{1}{\sqrt{2}} ( G(z) + \bar{G}(z) ) +
			2
     			\theta T(z) , \\
		X^+ (Z) &= X(z) + \theta \sqrt{2} Y(z), & X^-(Z) &= \bar{X}(z) +
		\theta \sqrt{2} \bar{Y}(z),
	\end{aligned}
	\end{equation*}
	we obtain the generating fields as in \cite{odake}.  
\end{ex}

\bibliographystyle{alpha}
\bibliography{refs}
\end{document}